\documentclass[10pt,twoside,english,a4paper]{article}
\usepackage{etex}
\usepackage{caption}
\usepackage[latin1]{inputenc}
\usepackage[intlimits] {amsmath}   % muss for defjs# stehen
\usepackage{graphicx}
\usepackage{psfrag}
\usepackage{ifthen}
\usepackage{fancyhdr}
\usepackage{rotating}
\usepackage{multirow}
\usepackage{booktabs}              % Dicke Tabellenlinien
\usepackage{amssymb}
\usepackage{amsbsy}
\usepackage{bm}                    % fette Schrift in Formeln
\usepackage{bbm}
\usepackage{babel}
\usepackage{theorem}
\usepackage{upgreek}  % gerade (roman) griechische Buchstaben z.B. \upalpha
\usepackage{pifont}   % zusaetzliche Schriften, z.B. eingekreiste Zahlen
\usepackage[active]{srcltx}   % aktive suche im xdv
\usepackage{subcaption}
\usepackage{suetterl}  % Suetterlin Schrift
\newfont{\suetdbl}{suet14 scaled 2000}  % Suetterlin skaliert um Faktor 2,000
\usepackage{comment}
\usepackage{yfonts}  % altdeutsche Schrifte und Initialen
\newfont{\gothdbl}{ygoth scaled 2000} 
\newfont{\frakdbl}{yfrak scaled 2000}
\newfont{\swabdbl}{yswab scaled 2000}
\usepackage{caption}
\usepackage{subcaption}
\usepackage{placeins} 
\usepackage[normalem]{ulem}  % zusaetzliche Unterstreichungsformen wie z.B.:
\usepackage[textsize=footnotesize,textwidth=1.7cm]{todonotes}
\usepackage[]{defjs3}      % eigene Definitionen: Option BoldVarsRoman
\usepackage{stmaryrd}
\usepackage{subcaption}
\usepackage[]{defjs3}  
\usepackage{amsfonts}
\usepackage{graphicx}
\usepackage{amsmath}
\usepackage{comment}

\usepackage[numbers,sort&compress]{natbib}
\bibpunct{[}{]}{,}{n}{}{;}
\fancypagestyle{plain}{%
  \fancyhead{}%
  \fancyfoot[c]{\sffamily\thepage}%
}
\makeatletter                           % Leere Fuellseiten
\def\cleardoublepage{\clearpage\if@twoside \ifodd\c@page\else
  \hbox{}
  \vspace*{\fill}
  \thispagestyle{empty}
  \newpage
  \if@twocolumn\hbox{}\newpage\fi\fi\fi}
\makeatother
\def\Curl{\operatorname{Curl}}      % Curl  ohne Argument
\def\curl{\operatorname{curl}}       % curl  ohne Argument
\newcommand{\Ce}{\mathbb{C}_{\mathrm{e}}}
\newcommand{\Cc}{\mathbb{C}_{\mathrm{c}}}
\newcommand{\Cmicro}{\mathbb{C}_{\mathrm{micro}}}
\newcommand{\Cmacro}{\mathbb{C}_{\mathrm{macro}}}
\newcommand{\CVoigt}{\mathbb{C}_{\mathrm{Voigt}}}
\newcommand{\Cmatrix}{\mathbb{C}_{\mathrm{matrix}}}
\newcommand{\Cinclusion}{\mathbb{C}_{\mathrm{inclusion}}}
\newcommand{\Lc}{L_{\mathrm{c}}}
\newcommand{\IL}{\mathbb{L}}
\newcommand{\R}{\mathbb{R}}
\newcommand{\Bdis}{\bP}
\newcommand{\dis}{P}
\begin{document}
\unitlength1.0cm
\frenchspacing

%=== cover page
%\input{fgeneralized_continua_01}

%=== sections ===
\thispagestyle{empty}

\vspace{1mm}
\ce{\bf 
Size-effects of metamaterial beams subjected to pure bending: on boundary} \ce{\bf conditions and parameter identification in the relaxed micromorphic model}

\vspace{3mm}
\ce{Mohammad Sarhil$^{1,\ast}$, Lisa Scheunemann$^2$, J\"org Schr\"oder$^1$ and Patrizio Neff$^3$}

\vspace{3mm}
\ce{$^1$Institute of Mechanics, University of Duisburg-Essen}
\ce{Universit\"atsstr. 15, 45141 Essen, Germany}
\ce{\small e-mail: 
mohammad.sarhil@uni-due.de
}    

\vspace{2mm}
\ce{$^2$Chair of Applied Mechanics, Gottlieb-Daimler-Str., TU Kaiserslautern,}
\ce{ 67663 Kaiserslautern, Germany}

\vspace{2mm}
\ce{$^3$Chair of Nonlinear Analysis and Modeling, Faculty of Mathematics, }
\ce{Thea-Leymann-Str. 9, University of Duisburg-Essen,}
\ce{ 45141 Essen, Germany}

\vspace{2mm}
\begin{center}
{\bf \large Abstract}
\bigskip

{\footnotesize
\begin{minipage}{14.5cm}
\noindent
In this paper we model the size-effects of metamaterial beams under bending with the aid of the relaxed micromorphic continuum. We analyze first the size-dependent bending stiffness of heterogeneous fully discretized  metamaterial beams subjected to pure bending loads. Two equivalent loading schemes are introduced which lead to a constant moment along the beam length with no shear force. The relaxed micromorphic model is employed then to retrieve the size-effects. We present a procedure for the determination of the material parameters of the relaxed micromorphic model based on the fact that the model operates between two well-defined scales. These scales are given by linear elasticity with micro and macro elasticity tensors which bound the relaxed micromorphic continuum from above and below, respectively.  The micro elasticity tensor is specified as the maximum possible stiffness that is exhibited by the assumed metamaterial while the macro elasticity tensor is given by standard periodic first-order homogenization. For the identification of the micro elasticity tensor, two different approaches are shown which rely on affine and non-affine Dirichlet boundary conditions of candidate unit cell variants with the possible stiffest response. The consistent coupling condition is shown to allow the model to act on the whole intended range between macro and micro elasticity tensors for both loading cases. We fit the relaxed micromorphic model against the fully resolved metamaterial solution by controlling the curvature magnitude after linking it with the specimen's size. The obtained parameters of the relaxed micromorphic model are tested for two additional loading scenarios. 

\end{minipage}
}
\end{center}


%=========================================================================
\sect{\hspace{-5mm}. Introduction}
\label{sec:intor}
Mechanical metamaterials are unconventional materials with exotic mechanical properties that are governed by the geometry of the complex underlying microstructure rather than by the properties of the constituting materials \cite{FisHilEbe:2020:mmo,GolGryGioLauCos:2019:mwr,LeeSinTho:2012:mnm,YaZhoLiaJiaWu:2018:mma,Zad:2016:mmm}. They can be optimized to obtain the intended mechanical properties to fit the wanted functionality \cite{SurGaoDuLiXioFanLu:2019:AEM}. However, mechanical metamaterials typically reveal size-effect phenomena and therefore the classical Cauchy-Boltzmann theory and first-order homogenization methods are incapable to describe such mechanical behavior. Generalized continua are enhanced continua that can model these size-effects as a homogeneous continuum without accounting for the detailed microstructure. The enhancement can be achieved by expanding the kinematics to contain additional degrees of freedom, e.g. the classical micromorphic theory \cite{Eri:1968:mom,EriSub:1964:nto,JuMahLiaXu:2021:CS,Min:1964:msi,SuhEri:1964:nto,NefFor:2007:age} and the Cosserat theory \cite{AlaGanSad:2022:cch,CosCos:1909:tof,LeiMah:2015:coh,Nef:2006:tcc,NefJeoMueRam:2010:lce,TroPau:2014:dom}, or by accounting for higher-grade differential operators in the energy functional, e.g. gradient elasticity models \cite{Aif:2011:otg,AltAif:1997:osa,AskAif:2011:gei,AskMetPicBen:2008:fsg,EreCazDel:2021:ond,FisKlaMerSteMue:2011:iao,GodGan:2016:cof,MinEsh:1968:ofsg,SheAbaBer:2022:abs,YanTimAbaLiMue:2021:vos}. However, the identification of the material parameters of these models is not trivial and in general remains unsolved. Different schemes were presented for the homogenization of  the heterogeneous fully resolved microstructures into the Cosserat continuum in \cite{AlaGanSad:2022:cch,ForSab:1998:com,Hue:2019:otm,RedAlaNasGan:2021:htc},
different variants of the gradient elasticity continuum in \cite{AbaBar:2021:ami,AbaYanPap:2019:aca,BacPagDalBig:2018:iof,KhaNii:2020:asg,LahGodGan:2022:sis,SchKruKeiHes:2022:cho,SkrEre:2020:ote,Wee:2021:nho,YanMue:2021:seo,YanAbaTimMue:2020:dom,YanAbaMueetal:2022:voa} and the
classical Eringen-Mindlin micromorphic continuum in \cite{AlaGanRedSad:2021:com,BisPoh:2017:amc,For:2002:hma,Hue:2017:hoa,RokAmePeeGee:2019:mch,RokAmePeeGee:2020:emc,RokZemDosKry:2020:ris,ZhiPohTayTan:2022:dfm}, however, without leading to a universally accepted answer. Mainly two approaches are employed for the determination of higher-order homogenized properties, which are asymptotic expansion methods, see e.\,g. \cite{Bou:1996:mei,BacGam:2010:soc} (also in combination with fast Fourier transform methods, see \cite{TraMonBon:2012:amb,LiZha:2013:ana}) and heuristic approaches relying upon the ad-hoc definition of modified kinematic boundary conditions on the microscale compared to first-order problems, see \cite{GodGan:2016:cof,BerDeoGodPicGan:2017:cof}. Among the latter, quadratic boundary conditions have been applied and analyzed to a large extent, see e.\,g. \cite{For:2002:hma,ForTri:2011:gca,TriJaeAufDieFor:2012:eog,For:2016:nro,KouGeeBre:2002:msc,KouGeeBre:2004:mss,AufBouBre:2010:sge}, in the field of homogenization towards second gradient continua and classical micromorphic continua. However, several problematic issues are described in the literature for this choice. Indeed, this natural extension does not lead to vanishing effective higher-order moduli when a homogeneous RVE or unit cell is homogenized. Moreover, when scale separation holds, Cauchy theory is also not a priori recovered, cf. \cite{AlaGanRedSad:2021:com}. To correct some of these spurious effects, additional microstructure-dependent body forces along with quadratic boundary conditions have been introduced \cite{YvoAufMon:2020:cso,MonAufYvo:2020:sgh}. Even though the presented results agree well with results from asymptotic homogenization, there remain artifacts for the special case of soft inclusions in a hard matrix material, for which the higher-order properties diverge. An alternative formulation has been proposed in \cite{Hue:2019:otm} by averaging solely over microheterogeneities in a homogenization scheme from a Cauchy continuum to micromorphic media. Note that \cite{TriJaeAufDieFor:2012:eog} stated that a quartic boundary condition would be necessary for fully describing the moduli in a micromorphic theory. A harmonic decomposition is recently applied to the governing equations and used to interpret the related modes on the micro-scale in  \cite{Hue:2022:iom}. 
 A variational approach is presented for the homogenization from a Cauchy continuum on the lower scale towards a second gradient or micromorphic continuum on the macro-scale in  \cite{GanRed:2021:ava} and \cite{ AlaGanRedSad:2021:com}, respectively. In their approach, as in other theories, the microscopic displacement is decomposed into a homogeneous and a fluctuation part. However, in contrast to other theories involving a heuristic definition of boundary conditions, the homogeneous part of the deformation arises here from a variational approach.   
With an enhancement of the description on the microlevel \cite{Wee:2021:nho} presents a homogenization procedure from metamaterial unit cell structures modeled using beam-lattice  structures on the micro-scale to a second gradient linear elastic model on the macro-scale. Similar to \cite{ForTri:2011:gca}, zero energy modes are observed for the higher-order moduli. In \cite{SchKruKeiHes:2022:cho}, a homogenization strategy for higher-order continua is presented which scales from a second-order continuum on the meso-scale to second- and third-order continua on the macro-scale under application of Isogeometric Analysis (IGA) and thereby also enhances the micro-scale continuum. Further developments in the context of multiphysical applications are e.g. discussed in \cite{WasHeuStaGeeKou:2020:ecf}.   
 
 In the field of asymptotic expansion homogenization and especially homogenization of metamaterials, the authors in \cite{AbaVazNew:2022:iom} have used asymptotic homogenization for the analysis of different unit cells in the framework of metamaterials. They exploited insight on the beam bending problem with a focus on the observable size-effects. In \cite{AbaBar:2021:ami}, metamaterials with honeycomb microstructure are analyzed in the framework of asymptotic expansion homogenization. The work \cite{AbaYanPap:2019:aca} presents a straight forward computational scheme for the determination of effective moduli through comparison with microstructure simulations. Here, the model is chosen apriori and does not originate from a homogenization strategy.
 Numerical and analytical solutions have been compared on a 3D structure for different deformation modes in  \cite{YanTimAbaLiMue:2021:vos} and pointed out the necessity of wedge and double traction forces for a correct overlap of both solutions. In \cite{YanMue:2021:seo}, mechanical metamaterials are analyzed by means of asymptotic expansion with an eye on appearing size-effects, which could only be detected for shear and torsion modes.

The relaxed micromorphic model considered by us is a generalized continuum model that allows in principle to capture size-effects and to describe band gaps phenomena in the dynamical case, see for example \cite{AivTalAgoDaoNefMad:2020:fan,BarTalDagAivNefMad:2019:rmm,DagBarGhiEidNefMad:2020:edo,DemRizColNefMad:2022:uem,MadNefGhiPlaRos:2015:bgi,MadNefBar:2016:cbg,MadNefGhiRos:2016:rat,MadNefBarGhi:2017:aro,RizdAgNefMad:2022:mam,RizNefMad:2022:msf,RizTalNefMad:2022:ttc,RizNefMad:2022:msf}. This model has been introduced in \cite{GhiNefMadPlaRos:2015:trl,NefGhiMadPlaRos:2014:aup} and its well-posedness for the static and dynamic problems has been proved in \cite{NefGhiLazMad:2015:trl,OwcGhidAgNef:2019:nmi}.  In \cite{KneOwcNef:2022:alr} the regularity of the model was investigated.  Being a micromorphic model, it features the classical translational degrees of freedom $\bu: \B \in \R^3 \rightarrow \R^3$ as well as a non-symmetric micro-distortion field $\Bdis: \B \in \R^3 \rightarrow \R^{3\times3}$. Compared to the classical micromorphic approach, the assumed strain energy is drastically simplified; notably, the curvature part (derivatives of $\Bdis$) intervenes only through $\Curl \Bdis$, so that solutions are found in $H^1(\B) \times H(\curl,\B)$ for the pair $(\bu,\Bdis)$. Using only the $\Curl$ of $\Bdis$ has some decisive advantages. It generates "bounded stiffness" \cite{RizHueKhaGhiMadNef:2021:aso1,RizHueMadNef:2021:aso3,RizHueMadNef:2021:aso4,RizKhaGhiMadNef:2021:aso2} for arbitrary large characteristic length (arbitrary small samples), in opposition to all strain gradient, Cosserat-micropolar or classical micromorphic approaches. Moreover, the appearing length-scale independent elasticity tensors $\Ce$ and $\Cmicro$ are related by a Reuss-like homogenization formula as function of the uniquely known elasticity tensor $\Cmacro$ from classical periodic homogenization. It remains therefore to determine $\Cmicro$, which happens to be the largest observable stiffness in the model (such an identification does not exist for the classical micromorphic model or other variants of it). As it turns out, the relaxed micromorphic model interpolates between two well-defined scales: the classical continuum scales of macroscopic elasticity, whose stiffness is given by $\Cmacro$ and a microscopic scale, with stiffness $\Cmicro$. The role of the characteristic length $\Lc > 0$ is then to scale correctly with the size of the specimen and to describe the interaction between the two scales. For $\Lc \rightarrow 0$ we recover macroscopic elasticity (complete scale separation, stiffness $\Cmacro$) and for $\Lc \rightarrow \infty$ (zoom into the microstructure) we obtain the microscopic scale (stiffness $\Cmicro$).

 In this contribution, we want to explore the possibilities that this unique interpretation of the relaxed micromorphic model provides. We consider an architected  material (hard matrix with soft inclusions). The determination of $\Cmacro$ is a standard identification in periodic homogenization theory. The identification of $\Cmicro$ will be guided by the largest stiffness idea alluded to above. Therefore, we consider a bending test of slender metamaterial beams.  The size-dependent bending was analyzed by means of other enriched models such as strain gradient, Cosserat-micropolar and other continua in \cite{AbaVazNew:2022:iom,AlAbdMah:2015:sdb,HosNii:2022:3sg,KhaNii:2020:asg,KhaBalNii:2018:msd,LiWanSonCheSuZhoWan:2022:osg,Lak:2022:cse,LieMue:2016:cof,KhaNii:2019:lsa,YinXiaDenetal:2021:iao}.
Modeling the mechanical behavior of many metamaterials was achieved for a variety of applications using generalized continua in \cite{Abi:2019:rtp,CarDelEspPul:2015:mdo,DelGio:2014:dpo,Dhaba:2020:rmm,PlaBarBat:2017:aim,GlaBasConetal:2021:ctm,RugHaLak:2019:cel,SheAbaBarBer:2021:iao,ShiFanTroLiWei:2022:sfe,SirLiuKouGee:2018:hec,SirKouGee:2016:hol}. 

 In this work the size-effects of metamaterial beams with fully discretized microstructure are analyzed.  Afterward, we employ the relaxed micromorphic continuum to describe these size-effects without accounting for the detailed microstructure. The material parameters and adequate boundary conditions of the micro-distortion field ${\Bdis}$ should be identified in order to establish a simplified fitting procedure on the fully resolved metamaterial beams.  The so-called consistent coupling condition (applied on the Dirichlet boundary for $\bu$) allows the relaxed micromorphic to operate on the whole scale between $\Cmacro$ and $\Cmicro$ which is of pivotal importance for a correct identification of its material parameters. However, an alternative loading by a normal linear traction (applied moment), which delivers exactly the same results for the fully resolved metamaterial, achieves consistent results as well for the relaxed micromorphic model when the consistent coupling condition is imposed via the penalty approach on the part of the boundary where the traction is set. 

In a previous  attempt \cite{NefEidMad:2019:ios} $\Cmicro$ was supposed  to be given by the L\"owner matrix supremum  ${\mathbb{C}}^\textrm{L\"owner}_\textrm{micro}$ of elasticity tensors appearing under affine Dirichlet conditions on the unit cell level. From the results in the present paper it inspires that ${\mathbb{C}}^\textrm{L\"owner}_\textrm{micro}$ is too soft, when compared with the appearing stiffness in the bending regime. Here, we extend our understanding of $\Cmicro$ towards all scenarios, notably including non-affine Dirichlet conditions. We limit our consideration to the planar case, in which the isotropic curvature energy in terms of $\Curl \Bdis$ has only one free parameter. 

The outline of the paper is as follows: in  Section \ref{sec:model} we recall the energy functional of the relaxed micromorphic model, define the material parameters, and introduce the strong forms with the associated boundary conditions obtained by the energy minimization. We present briefly in Section \ref{sec:fem} the main aspects of the construction of $H(\curl,\B)$-conforming finite elements. The size-effects of the heterogeneous microstructured metamaterial beams are investigated in Section \ref{sec:reference} for two loading cases which lead to the same results. In Section \ref{sec:iden} we determine the material parameters of the relaxed micromorphic model and discuss the boundary condition for symmetric and non-symmetric force stresses. We then fit the relaxed micromorphic model solution to the microstructured metamaterial solution by calibrating the curvature in Section  \ref{sec:fit}. In Section \ref{sec:validation}, the relaxed micromorphic model is shown to be capable of handling two loading scenarios in addition to pure bending. Finally, we provide our conclusions and outlook in Section \ref{sec:con}. 

%=========================================================================

\sect{\hspace{-5mm}. The relaxed micromorphic model and its discretization}
\subsection{The relaxed micromorphic model}
\label{sec:model}  
The relaxed micromorphic model (RMM) is an enriched continuum model. The kinematics of each material point is determined, similar to the general micromorphic theory \cite{EriSub:1964:nto,Min:1964:msi,SuhEri:1964:nto}, by a displacement vector $\bu\colon\B\subseteq\R^3\to\R^3$ and a non-symmetric micro-distortion field  $\Bdis\colon\B\subseteq\R^3\to\R^{3\times3}$. The displacement and the micro-distortion fields are defined for the static case by minimizing the energy functional  
\begin{equation}
\label{eq:pot}
\Pi (\bu,\Bdis)= \int_\B W\left(\nabla \bu,\Bdis,\Curl \Bdis \right) \, - { \overline\bbf}\cdot{\bu} \, \,\textrm{d}V\ - \int_{\partial \B_t}  \overline\bt \cdot \bu  \,\textrm{d}A  \longrightarrow\ \min\,,
\end{equation}
with $(\bu,\Bdis)\in H^1(\B)\times H(\curl,\B)$. The vector $ \overline\bbf$  describes the applied body force. The vector $ \overline\bt$ is the traction vector acting on the boundary $\partial \B_t \subset \partial \B$. The elastic energy density $W$ reads 
\begin{equation}
\label{eq:W}
\begin{aligned}
W\left(\nabla \bu,\Bdis,\Curl \Bdis \right) = & \frac{1}{2} ( \symb{ \nabla \bu - \Bdis} : \Ce : \symb{ \nabla \bu - \Bdis}  +   \sym \Bdis : \Cmicro: \sym \Bdis  \\
& + \skewb{ \nabla \bu - \Bdis} : \Cc : \skewb{ \nabla \bu - \Bdis}  +  \mu \, \Lc^2 \, \textrm{Curl} \Bdis : \IL :\textrm{Curl} \Bdis  )\,.
\end{aligned}
\end{equation}
Here, $\Cmicro,\Ce > \bzero $ are fourth-order positive definite standard elasticity tensors, $\Cc \ge \bzero$ is a fourth-order positive semi-definite rotational coupling tensor, $\IL$ is a positive definite fourth-order tensor acting on non-symmetric arguments, $\Lc \ge 0 $ is the characteristic length parameter and $\mu$ is a shear modulus for dimensional consistency. The  characteristic length parameter is related to the size of the microstructure  and determines its influence on the macroscopic mechanical behavior.  The characteristic length allows to scale the number of considered unit cells keeping all remaining parameters of the model scale-independent where the macro-scale with $\Cmacro$ and the micro-scale with $\Cmicro$ are retrieved for $\Lc \rightarrow 0$ and  $\Lc \rightarrow \infty$, respectively, if suitable boundary conditions are applied, see \cite{NefEidMad:2019:ios,SchSarSchNef:2022:lhb}. The macro-scale elasticity tensor $\Cmacro$ associated with $\Lc \rightarrow 0$ can be defined by the standard first-order periodic homogenization (the scale separation holds) while the micro-scale elasticity tensor $\Cmacro$ associated with $\Lc \rightarrow \infty$ represents the stiffest extrapolated response (zooming in the microstructure). The constitutive coefficients are assumed constant with the following symmetries 
\begin{equation}
\begin{aligned}
&(\Cmicro)_{ijkl} = (\Cmicro)_{klij} = (\Cmicro)_{jikl} \,, \qquad  &(\Cc)_{ijkl} = (\Cc)_{klij} \,, \\
&(\Ce)_{ijkl} = (\Ce)_{klij} = (\Ce)_{jikl} \,, \qquad   &{ (\IL)_{ijkl} = (\IL)_{klij}  \,,}
\end{aligned}
\end{equation} 
where $\Cmicro$ and $\Ce$ are connected to $\Cmacro$ through a Reuss-like homogenization relation { \cite{BabMadDagAbrGhiNeff:2017:taf}}
\begin{equation}
\Cmacro^{-1} = \Cmicro^{-1} + \Ce^{-1} \quad \Rightarrow \quad \Ce = \Cmicro (\Cmicro-\Cmacro)^{-1} \Cmacro\,.
\end{equation}

The variation of the potential with respect to the displacement yields the weak form
\begin{equation}
\begin{split}
\delta_{\bu} \Pi =& \int_{\B}\{  \underbrace{ \Ce : \symb{ \nabla \bu - \Bdis}+  \Cc : \skewb{ \nabla \bu - \Bdis}}_{\textstyle=:\Bsigma} \} :  {\nabla \delta \bu} -  \overline\bbf \cdot \delta \bu  \, \textrm{d}V - \int_{\partial \B_t}  \overline\bt \cdot \delta \bu  \, \textrm{d}A = 0\,,
\end{split}
\end{equation} 
which leads, using  integration by parts and employing the divergence theorem, to
 \begin{equation}
\delta_{\bu} \Pi =  \int_{\B}   \{ \div \Bsigma +  \overline\bbf \} \cdot \delta \bu \, \textrm{d}V  = 0 \, ,
\end{equation}  
where  $\Bsigma$ is the non-symmetric force stress tensor (symmetric if $\Cc  \equiv \bzero$ which is permitted). In a similar way, the variation of the potential with respect to the micro-distortion field $\Bdis$ leads to the weak form
\begin{equation}
 \delta_{\Bdis} \Pi =   \int_{\B} \{ \Bsigma  - \underbrace{\Cmicro : \sym \Bdis}_{\textstyle =:\Bsigma_\textrm{micro} } +  \overline\bM\} :  \delta {\Bdis}  - \underbrace{\mu \, \Lc^2  (\IL : \Curl \Bdis)}_{\textstyle =:\bbm} : \Curl \delta {\Bdis} \, \textrm{d} V  = 0\,, 
\end{equation}
which can be rewritten, using integration by parts and applying Stokes' theorem, as 
\begin{equation}
 \delta_{\Bdis} \Pi =   \int_{\B}  \{ \Bsigma - \Bsigma_\textrm{micro}  - \Curl \bbm \} :  \delta {\Bdis} \, \textrm{d} V +  \int_{ \partial \B}  \{ \sum_{i=1}^3  \left( \bbm^i  \times \delta {\Bdis}^i  \right) \cdot \bn  \}  \;\; \textrm{d}A  
 = 0\,, 
\end{equation}
where the stress measurements $\Bsigma_\textrm{micro}  $ and $\bbm$ are the micro- and moment stresses, respectively, $\bn$ is the outward unit normal vector on the boundary, and $\bbm^i$ and $ \delta {\Bdis}^i$ are the row vectors of the related second-order tensors. 
The strong form of the relaxed micromorphic model with the associated boundary conditions read 
    \begin{subequations}
      \label{eq:sfs1}
 \begin{align}
  \div \Bsigma + \overline\bbf &= \bzero &\quad  \textrm{on} &\quad  \B , \\ 
  \bu &= \overline{\bu} \quad& \textrm{on} &\quad  \partial \B_u\,,  \\
  \overline\bt &= \Bsigma \cdot \bn \quad& \textrm{on} &\quad  \partial \B_t \,,  \\
 \Bsigma - \Bsigma_\textrm{micro}  - \Curl \bbm  &= 0  \quad& \textrm{on} &\quad  \B \,, \\  
  \sum_{i=1}^3   \Bdis^i \times \bn &=  \overline\bt_{P}  \quad& \textrm{on} &\quad  \partial \B_{P} \\  
  \sum_{i=1}^3   \bbm \bm^i \times \bn &= \bzero \quad& \textrm{on} &\quad  \partial \B_m\,, 
   \end{align}
 \end{subequations}
 where $\partial \B_\dis \cap \partial \B_m = \partial \B_u \cap \partial \B_t  = \emptyset $ and  $\partial  \B_\dis \cup \partial \B_m = \partial \B_u \cup \partial \B_t = \partial \B $.  The strong form represents a generalized balance of linear momentum (force balance) and a generalized balance of angular momentum (moment balance). For more details regarding derivations of the boundary conditions, the reader is referred to \cite{SchSarSchNef:2022:lhb}. \\
An additional dependence between the displacement field and the micro-distortion field on the boundary was proposed in \cite{NefEidMad:2019:ios} and subsequently considered in \cite{SkyNeuMueSchNef:2021:CM,RizHueMadNef:2021:aso3,RizHueMadNef:2021:aso4,DagRizKhaLewMadNef:2021:tcc}. This so-called {\bf consistent coupling condition} is defined by
 \begin{equation}
  \Bdis \cdot \Btau = \nabla \bu \cdot \Btau \, \Leftrightarrow \, {\Bdis}^i  \times \bn  =   {\nabla \bu}^i  \times \bn \quad \textrm{for} \quad i=1,2,3  \quad \textrm{on}\quad  \partial \B_\dis = \partial \B_u  \,,
 \end{equation}

  where $\Btau$ is the tangential vector on the boundary and ${\Bdis_i}$ and $\nabla {\bu}^i$ are the row-vectors of the associated tensors.  However, we can extend this relative boundary condition to
 parts of ${\partial} \B_m$ by enforcing the consistent coupling condition on $\partial \B_{\widehat{m}} \subseteq \partial \B_{{m}}$ via a penalty approach as
 \begin{equation}
 \Pi \Leftarrow \Pi + \int_{ \partial \B_{\widehat{m}}}   \frac{\kappa_1}{2} \sum_{i=1}^3 ||(\Bdis^i - \nabla {\bu}^i) \times \bn )||^2 \, \textrm{d}A  \,,
 \end{equation}
 where $\kappa_1$ is the penalty parameter. 
 
 The micro-distortion field has the following general form for the three-dimensional case
 \begin{equation}
\Bdis = \left( \begin{array}{c}
(\Bdis^{1})^T \\
(\Bdis^{2})^T \\
(\Bdis^{3})^T \\
\end{array}\right) = \left( \begin{array}{c c c}
 \dis_{11}  &  \dis_{12}  & \dis_{13}  \\   
 \dis_{21}  &  \dis_{22}  & \dis_{23}  \\ 
 \dis_{31}  &  \dis_{32}  & \dis_{33}  \\ 
\end{array}\right) \quad \textrm{with} \quad \Bdis^{i} = \left( \begin{array}{c}
\dis_{i1} \\
\dis_{i2} \\ 
\dis_{i3}
\end{array} \right)\,  \textrm{for } i=1,2,3\,.
 \end{equation} 
We let the Curl operator act on the row vectors of the micro-distortion field $\Bdis$ as
 \begin{equation}
\Curl \Bdis = \left( \begin{array}{c}
(\curl \Bdis^{1})^T \\
(\curl \Bdis^{2})^T \\
(\curl \Bdis^{3})^T \\
\end{array}\right) = \left( \begin{array}{c|c|c}
 \dis_{13,2} - \dis_{12,3}	& \dis_{11,3} - \dis_{13,1}  &  \dis_{12,1} - \dis_{11,2} \\
 \dis_{23,2} - \dis_{22,3}	& \dis_{21,3} - \dis_{23,1}  &  \dis_{22,1} - \dis_{21,2} \\
 \dis_{33,2} - \dis_{32,3}	& \dis_{31,3} - \dis_{33,1}  &  \dis_{32,1} - \dis_{31,2} 
\end{array}\right) \,.
 \end{equation}
\subsection{$H^1(\B)\times H(\curl,\B)$-conforming finite element in 2D}
\label{sec:fem} 
Different finite element formulations of the relaxed micromorphic model were introduced for the plane strain case in \cite{SchSarSchNef:2022:lhb,SarSchNefSch:2021:oat,SarSchSchNef:2023:mts}, antiplane shear in \cite{SkyNeuMueSchNef:2021:CM} and 3D case in \cite{SkyNeuMueSchNef:2022:pam,SkyMueRizNeff:2023:hob}. For the two-dimensional case, the micro-distortion field has only four non-vanishing components, which are in the plane, and its Curl operator is reduced to only two components out of the plane, namely $(\Curl \Bdis)_{13}$ and $(\Curl \Bdis)_{23}$, 
 \begin{equation}
\Bdis = \left( \begin{array}{c}
(\Bdis^{1})^T \\
(\Bdis^{2})^T \\
\bzero^T 
\end{array}\right) = \left(\begin{array}{c c c }
 \dis_{11}  &  \dis_{12}   & 0\\   
 \dis_{21}  &  \dis_{22}   & 0\\
 0 & 0 & 0 
\end{array}\right) \quad \textrm{and} \quad 
\textrm{Curl} \Bdis = 
\left(\begin{array}{c|c|c}
0	&   0  &  \dis_{12,1} - \dis_{11,2} \\
0   &   0  &  \dis_{22,1} - \dis_{21,2} \\
0 & 0 & 0
\end{array}\right)\,.
 \end{equation}

It has been shown in \cite{SchSarSchNef:2022:lhb} that  $H^1(\B)\times H(\curl,\B)$ elements obtain the discontinuous solution of the micro-distortion field while the standard nodal $ H^1(\B)\times H^1(\B)$ elements are unable to capture the jumps. Therefore, transition zones emerge for $ H^1(\B)\times H^1(\B)$ elements which need to be resolved by distinctly refining the mesh in contrast to $H^1(\B)\times H(\curl,\B)$ elements which exhibit faster convergences rates. 

\begin{figure}[ht]
	\unitlength=1mm
	\center
	\begin{picture}(50,40)
	\put(0,0){\def\svgwidth{5 cm}{\small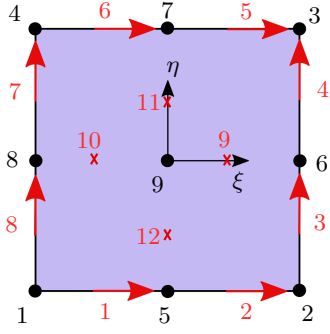}}
	\end{picture}
		  	\caption{Q2NQ2 Element. Black dots represent the displacement nodes. Red arrows and crosses indicate the edge and inner vectorial dofs, respectively, of the micro-distortion field used in N\'ed\'elec formulation. }
\label{Fig:Q2NQ2}
\end{figure}

We demonstrate briefly the main aspects of the finite element formulation of a quadrilateral element $(\bu,\Bdis) \in H^1(\B) \times H(\curl,\B)$ shown in Fig. \ref{Fig:Q2NQ2}. The finite element, denoted as Q2NQ2, utilizes Lagrange-type shape functions of the second-order for the displacement field, denoted as Q2. The suitable finite element space for the micro-distortion field is known as N\'ed\'elec space, see \cite{Ned:1980:mfe,Ned:1986:anf}. In this work, we choose the N\'ed\'elec space of first-kind and second-order, denoted as NQ2.  N\'ed\'elec formulation uses vectorial shape functions that satisfy the tangential continuity at element interfaces. General reviews about the edge elements are available in \cite{KirLogRogTer:2012:cau} and \cite{RogKirAnd:2009:eao}.  For more details regarding the derivation  of shape functions and the FEM-implementation aspects, the reader is referred  to \cite{SchSarSchNef:2022:lhb}. 

The Q2NQ2 element uses 9 nodes for the discretization of the displacement field $\bu$. The geometry and the displacement field are approximated employing the related quadratic scalar shape functions $N^u_I$ defined in the parameter space with natural coordinates $\Bxi = \{\xi,\eta\}$ by
\begin{equation}
\label{eq:t:approx_chi_u}
 \bX_h = \sum_{I=1}^{9} N^u_I ( \Bxi) \bX_I\,, \qquad \qquad \bu_h = \sum_{I=1}^{9} N^u_I ( \Bxi) \bd^u_I\,, 
\end{equation}
where $\bX_I$ are the coordinates of the displacement node $I$ and $\bd^u_I$ are its displacement degrees of freedom. The deformation gradient is obtained then in physical space by 
\begin{equation}
\nabla \bu_h = \sum_{I=1}^{9} \bd^u_I \otimes \nabla N^u_I ( \Bxi)  \,, \qquad \textrm{with} \qquad 
\nabla {N^u_I} ( \Bxi)  = \bJ^{-T} \cdot \nabla_\Bxi N^u_{I}\,,
\end{equation}
where $\bJ = \frac{\partial \bX}{\partial \Bxi  }$ is the Jacobian, $\nabla$ and $\nabla_\Bxi$ are the gradient operators to $\bX$ and $\Bxi$, respectively.  
The micro-distortion field $\Bdis$ is  approximated by the vectorial dofs $\bd^\dis_I$ presenting its tangential components at the location $I=1,...,12$. The micro-distortion field and its Curl operator are interpolated as 
\begin{equation}
\label{eq:t:dis_restoring}
\Bdis_h =   \sum_{I=1}^{12}   \bd^\dis_I \otimes \Bpsi^2_I \, , \qquad  \qquad 
\Curl \Bdis_h  =  \sum_{I=1}^{12}   \bd^\dis_I  \otimes \curl{{\Bpsi}^2_I} \,.
\end{equation}
The non-vanishing components of the Curl operator of the micro-distortion field for the 2D case are obtained by 
\begin{equation}
\left[\begin{array}{c}
\curl^{2D}{\Bdis_h^1} \\
  \curl^{2D}{\Bdis_h^2} \\
\end{array}\right]\,
 =  \sum_{I=1}^{12}   \bd^\dis_I  \curl^{2D}{{\Bpsi}^2_I}  =  \left[\begin{array}{c}
\sum_{I=1}^{12}   (d^\dis_I)_1  \curl^{2D}{{\Bpsi}^2_I} \\
 \sum_{I=1}^{12}   (d^\dis_I)_2  \curl^{2D}{{\Bpsi}^2_I}  \\
\end{array}\right]. 
\end{equation}

The simulations presented in this paper are performed within AceGen and AceFEM programs. The interested reader is referred to \cite{KorWri:2016:aofem,Kor:2009:aof}.

{\section{Reference study: size-effects of metamaterial specimens subjected \\ to bending}
\label{sec:reference}
We investigate here the size-effect phenomena of an assumed metamaterial with fully resolved microstructure.  The size-effect phenomena will be analyzed via the effective bending stiffness  of beams subjected to pure bending.  According to the elementary beam theory, the moment is linked to the curvature by 
$             
M (x) = D(x) \kappa(x), 
$
where $D(x)$ and $\kappa(x)$ are the bending stiffness and the curvature at a position $x$ along the beam. For a constant bending moment  $\overline M$ along the beam length, we assume an effective  flexural rigidity $\overline{D}$ and an effective curvature $\overline{\kappa}$ so that we obtain 
\begin{equation}
\label{Eq:bending_stiffness}
\overline{D} = \frac{\overline M}{ \overline{\kappa}}.
\end{equation}

We design in the following two beams subjected to a vanishing shear force and a constant moment along the length $L$, see Fig. \ref{Figure:beams}. For the first loading case a rotation $\theta$ is applied on the right end while a moment load is enforced for the second loading case instead. A deflection equation $ \overline w(x)$, which will be fitted later to the heterogeneous beams, featuring an effective constant curvature reads

\begin{equation}
\label{Eq:deflection}
 \overline w(x) = \frac{\overline \kappa}{2} (x^2-L^2 )\,
 \quad \quad \, \textrm{satisfying}  \quad \overline{w}(L) = 0, \quad \quad \, \textrm{and} \quad  \frac{\textrm{d}\overline{w}(0)}{\textrm{d}x} = 0 \, , 
\end{equation}

  \begin{figure}[ht]
\center
	\unitlength=1mm
	\begin{picture}(140,70)
	\put(0,5){\def\svgwidth{14 cm}{\small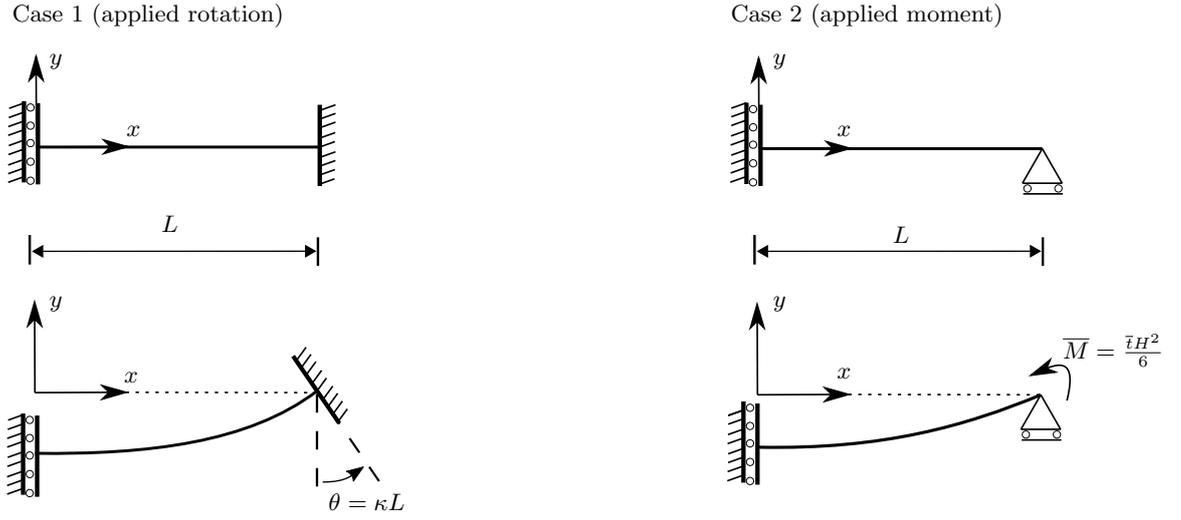}}
	\end{picture}
	\caption{The beam models, compare Fig. \ref{Figure:beam_metamaterials_BCs}.}
	\label{Figure:beams}
\end{figure}

A 2D metamaterial is considered with a unit cell consisting of a square with an edge length $l=1.9 \cdot 10^{-2} \, \textrm{m}$ and a circular inclusion at its center with a diameter of $d=1.2 \cdot 10^{-2} \, \textrm{m}$, see Fig. \ref{Figure:unit_cell}. Both matrix and inclusion are isotropic linear elastic with the material parameters shown in Table \ref{Tab:matpar1}.  The inclusion is 20 times softer than the matrix.  A standard triangular finite element with quadratic shape functions (T2) is used for this analysis. 
The specimens are considered with dimensions $H \times L = n \, l \times 12 \, n \, l $  so { that} the length is always twelve times the height where $n$ is the number of unit cells in the height direction, see Fig. \ref{Figure:unit_cell}. 

\begin{table}[ht]
\caption{Material parameters of the assumed metamaterial.}
\label{Tab:matpar1}
\center
\begin{tabular}{ |c|c|c|c|c| }
\hline
& Young's modulus: $E$ & Poisson's ratio: $\nu$ & $\lambda$ & $\mu$  \\ \hline 
Matrix &  $70 \, \textrm{GPa}$ & $0.333$ & $52.35 \, \textrm{GPa}$ &  $26.25 \, \textrm{GPa}$ \\ \hline 
Inclusion &  $3.5 \, \textrm{GPa}$ & $0.333$ & $2.62 \, \textrm{GPa}$ &  $1.31 \, \textrm{GPa}$  \\ \hline 
\end{tabular}
\end{table}

\begin{figure}[ht]
\center
	\unitlength=1mm
	\begin{picture}(140,110)
	\put(0,5){\def\svgwidth{14 cm}{\small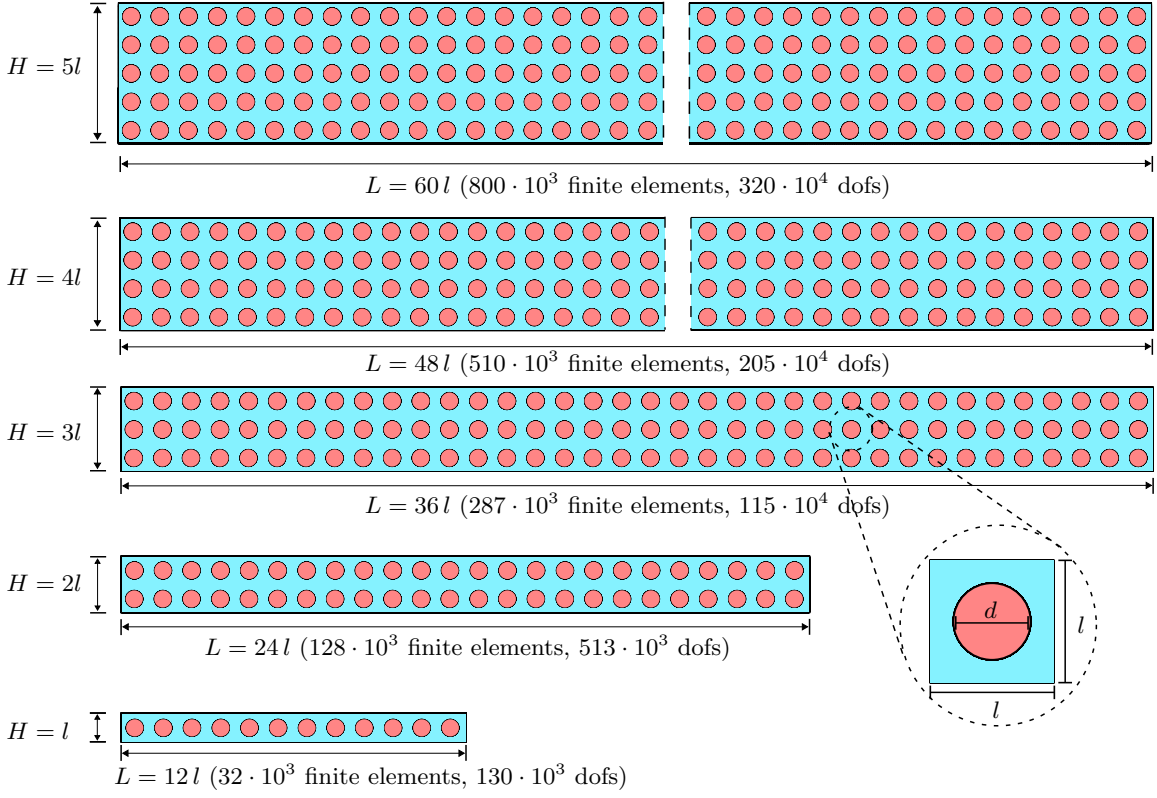}}
	\end{picture}
	\caption{Illustration shows the geometry of the specimens for $n=1,2,3,4,5$ with the assumed unit cell. The number of finite elements with degrees of freedom (dofs) are shown in parentheses.  }
	\label{Figure:unit_cell}
\end{figure} 
\FloatBarrier

The boundary condition of the beam models in Fig. \ref{Figure:beams} are passed on the 2D metamaterial as shown in Fig. \ref{Figure:beam_metamaterials_BCs}. For the first loading case we rotate the right edge in plane through a given displacement in $x$-direction as a linear function of $y$-coordinates while for the second loading case a moment is applied on the right edge by means of a traction in $x$-direction as a linear function of $y$-coordinates.  The left boundary for both loading cases is fixed in $x$-direction and free to move in $y$-direction. Furthermore, we fix the middle point on the right edge in $y$-direction. We intend by introducing these two loading cases to prove that they deliver identical results for the microstructured metamaterial beams. This equivalence should then be demonstrated as well by the relaxed micromorphic model when appropriate boundary conditions are set. Furthermore, we assume $\kappa=1$ and $\overline{t}=10^9 \, \textrm{N}/\textrm{m} $\,.

 \begin{figure}[ht]
\center
	\unitlength=1mm
	\begin{picture}(140,35)
	\put(0,5){\def\svgwidth{14 cm}{\small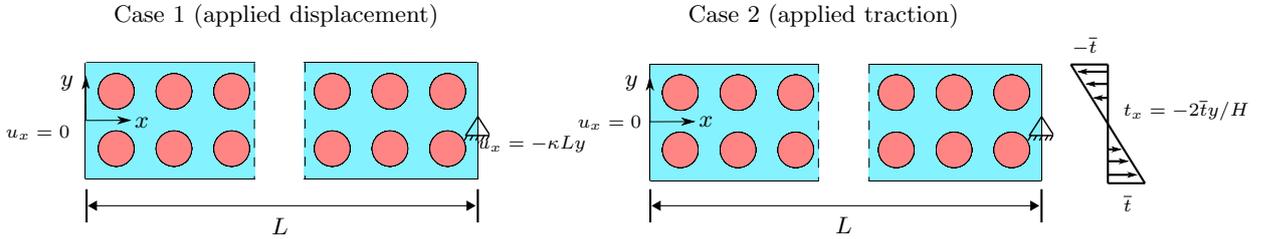}}
	\end{picture}
	\caption{The boundary conditions of the fully resolved metamaterial shown exemplarily for $n=2$ ($H \times L = 2 \, l \times 24 \, l$).}
	\label{Figure:beam_metamaterials_BCs}
\end{figure}

\FloatBarrier
After solving the fully resolved microstructure, the effective curvature $\overline\kappa$ is obtained by the following least square minimization   
\begin{equation}
\label{eq:ls}
\sum_I^{n_\textrm{node}} ( \bd^u_I -  \overline w(\bX_I))^2 \textrm{d} x \rightarrow \textrm{min}\,,
\end{equation}
which leads, considering Eq. \ref{Eq:deflection}, to 
\begin{equation}
\overline \kappa = \frac{\sum\limits_I^{n_\textrm{node}} (\bd^u_I)_2 \frac{(\bX_I)_1^2 - L^2 }{2}}{\sum\limits_I^{n_\textrm{node}} (\frac{(\bX_I)_1^2 - L^2)}{2})^2},
\end{equation} 
where $\bX_I$ and $\bd^u_I$ are the coordinates and the displacement degrees of freedom at  node $I$. The bending stiffness can be calculated following  Eq. \ref{Eq:bending_stiffness} where the moment $\overline  M$ can be calculated using the nodes reactions on the left or right edges. Alternatively, the bending stiffness can be calculated by means of the maximum deflection at the left edge of the beam. We obtain from Eqs. \ref{Eq:bending_stiffness} and \ref{Eq:deflection} substituting $x=0$ and considering $\overline w(0) = w^\textrm{FEM}(0)$ since the deflection's fluctuation of the heterogeneous solution is small compared to the maximum deflection
\begin{equation}
\label{eq:bs2}
\overline D = - \frac{\overline  M L^2}{2 \,  w^\textrm{FEM}(0) } \, , 
\end{equation}   
where $w^\textrm{FEM}(0)$ is the deflection of the FEM solution averaged over the left edge ($x=0$). Calculating the bending stiffness using Eqs. \ref{Eq:bending_stiffness} or \ref{eq:bs2} delivers the same result which we tested numerically.

The effective material properties of the large specimens can be obtained by the standard computational periodic first-order homogenization produced by a unit cell with periodic boundary condition which is identified as $\Cmacro$ in Section \ref{sec:cmacro}. As we will show later the macro elasticity tensor $\Cmacro$ is not isotropic and shows {cubic} symmetry.  The size-effects are shown via the so-called normalized bending stiffness $\overline{D}/D_\textrm{macro}$ plotted in Fig. \ref{Figure:metamateials_beam_stiffness} which relates the actual stiffness of the fully discretized metamaterial to the one obtained from  homogenized linear elasticity with $\Cmacro$ which reads analytically
\begin{equation}
\label{eq:macro_ana}
D_\textrm{macro} = \frac{E_\textrm{macro} \, H^3}{12\, (1-\nu_\text{macro}^2)}\,.
\end{equation}   

The normalized bending stiffness approaches the value one when we increase the specimen size. Applying a rotation (loading case 1) or a moment (loading case 2) leads to similar results as expected.

 \begin{figure}[ht]
	\unitlength=1mm
	\center
	\begin{picture}(100,60)
		 \put(27,52){\color[rgb]{0,0,0}\makebox(0,0)[lb]{\small\smash{$n=1$}}}%	 
		 \put(52,40){\color[rgb]{0,0,0}\makebox(0,0)[lb]{\small\smash{$n=2$}}}%	  	 	
		 \put(65,38){\color[rgb]{0,0,0}\makebox(0,0)[lb]{\small\smash{$n=3$}}}%
		 \put(85,36){\color[rgb]{0,0,0}\makebox(0,0)[lb]{\small\smash{$n=5$}}}%	  	   	 	
		 \put(47,15){\color[rgb]{0,0,0}\makebox(0,0)[lb]{\small\smash{case 1 (applied displacement)}}}%
		 \put(47,19.5){\color[rgb]{0,0,0}\makebox(0,0)[lb]{\small\smash{case 2  (applied traction)}}}%	  	   
         \includegraphics[width=0.6 \textwidth]{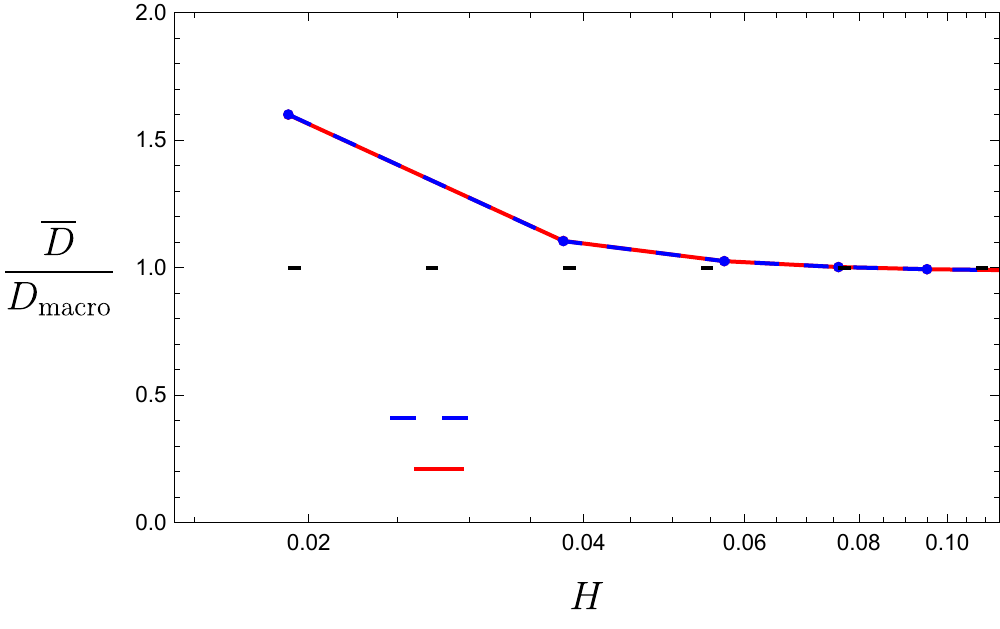}
		\end{picture}	  		\caption{The normalized bending stiffness varying the beam size $H \times L = n\,l \times 12 \, n \,  l$.  }
	\label{Figure:metamateials_beam_stiffness}
\end{figure}
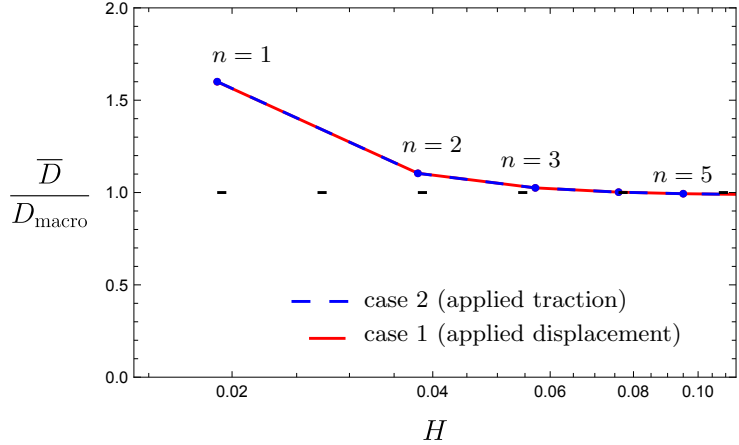

\FloatBarrier

\section{Size-effects of the relaxed micromorphic continuum subjected to pure bending}
\label{sec:iden}
The previous size-effects exhibited by the fully resolved heterogeneous material should be reproduced by the  relaxed micromorphic model. However, the boundary conditions and material parameters identification are not obvious.

\FloatBarrier
\subsection{Identification of $\Cmacro$}
\label{sec:cmacro}
 The macro elasticity tensor $\Cmacro$ corresponds to the case $\Lc \rightarrow 0$ equivalent to large values of $n$ where the macro homogeneous response is expected.  A unit cell with periodic boundary conditions should be used, see for example \cite{ZohWri:2005:ait}.  
 The geometry of the unit cell has no role for this standard analysis. Our analysis shows that $\Cmacro$ has the cubic symmetry property for our assumed metamaterial and it reads in Voigt notation
\begin{equation}
\Cmacro = \left( \begin{array}{c c c}
2 \mu_\textrm{macro} + \lambda_\textrm{macro} & \lambda_\textrm{macro} & 0 \\ 
\lambda_\textrm{macro} & 2 \mu_\textrm{macro} + \lambda_\textrm{macro} & 0 \\
0 & 0 & \mu^*_\textrm{macro}
\end{array}
\right)  \, , 
\end{equation}
where three parameters need to be defined. We obtain by our standard numerical analysis  
\begin{equation}
\Cmacro = \left( \begin{array}{c c c}
47.86 &  17.61 & 0 \\ 
17.61 & 47.86 &  0 \\ 
0 & 0 & 9.98 \\
\end{array}
\right)  [\textrm{GPa}] \quad \Rightarrow \quad  \begin{array}{c c c c}
\lambda_\textrm{macro} & = &  \, 17.61 \, & \textrm{GPa}\\ 
\mu_\textrm{macro} & = & 15.13 \, & \textrm{GPa}\\ 
\mu^*_\textrm{macro} & =  & \, 9.98 \, & \textrm{GPa}\\
\end{array}\,.
\end{equation}
\FloatBarrier
 \subsection{Identification of $\Cmicro$ (first approach)} 
 \label{sec:cmicro1}
  The micro elasticity tensor $\Cmicro$ in the relaxed micromorphic model is identified as the maximum stiffness on the micro-scale which must { exhibit} the cubic symmetry similar to $\Cmacro$ according to the extended Neumanns's principle \cite{NefEidMad:2019:ios}. In order to achieve stiff estimates for $\Cmicro$ we apply first affine Dirichlet boundary conditions. Furthermore, we have to choose unit cells which preserve the cubic symmetry under the applied Dirichlet boundary conditions. However, different variants of unit cells must be investigated for the affine Dirichlet boundary conditions.  For each choice of a unit cell $i=1,..,r$ with the affine Dirichlet boundary conditions, we obtain the corresponding apparent stiffness tensor denoted as $\mathbb{C}^\textrm{D}_i$. The positive definite micro elasticity tensor will be set as the least upper bound of the apparent stiffness of the microstructure measured in the energy norm following the L\"owner matrix supremum problem, see for details \cite{NefEidMad:2019:ios}. 

 For the assumed metamaterial, four different variants of the unit cell are suitable, see Fig. \ref{Figure:unit_cells}, which lead to the elasticity tensors $\mathbb{C}^\textrm{D}_i, i=1,..,4$ with the cubic symmetry property as intended. The results are summarized in Table \ref{Table:Cmicro}. 
 \begin{figure}[ht]
	\unitlength=1mm
	\center
 		 \begin{subfigure}[b]{1.0 \textwidth}
		 \center
	\begin{picture}(130,45)
	\put(0,5){\def\svgwidth{12 cm}{\small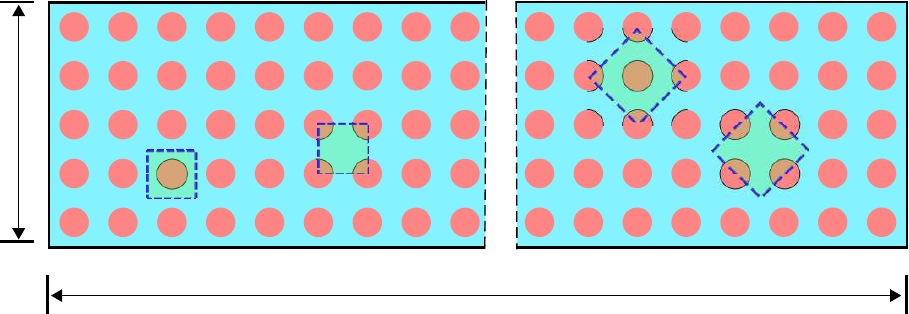}}
			 \put(65,4){\color[rgb]{0,0,0}\makebox(0,0)[lb]{\small\smash{$L$}}}%	  	
			 \put(-2,30){\color[rgb]{0,0,0}\makebox(0,0)[lb]{\small\smash{$H$}}}%	  	
	\end{picture}
   	\caption*{}
    \end{subfigure}
 		 \begin{subfigure}[b]{0.19 \textwidth}
		 \center
	\begin{picture}(15,25)
	\put(0,0){\def\svgwidth{1.5 cm}{\small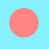}}
	\end{picture}
\caption*{(1)}
    \end{subfigure}
  		 \begin{subfigure}[b]{0.19 \textwidth}
		 \center
	\begin{picture}(15,25)
	\put(0,0){\def\svgwidth{1.5 cm}{\small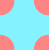}}
	\end{picture}
\caption*{(2)}
    \end{subfigure}
  		 \begin{subfigure}[b]{0.29 \textwidth}
		 \center
  	\begin{picture}(30,25)
	\put(0,0){\def\svgwidth{3 cm}{\small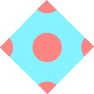}}
	\end{picture}
\caption*{(3)}
    \end{subfigure}
    		 \begin{subfigure}[b]{0.29 \textwidth}
		 \center
	\begin{picture}(30,25)
	\put(0,0){\def\svgwidth{3  cm}{\small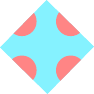}}
	\end{picture}
   		  	\caption*{(4)}
    \end{subfigure}
\caption{The possible choices of the unit cells with cubic symmetry. The edge length of the unit cell equals to $l$ for (1) and (2) and $\sqrt{2} \, l$ for (3) and (4).}
	\label{Figure:unit_cells}
\end{figure} 

\begin{table}[ht]
\caption{The elasticity parameters of the unit cells in Fig. \ref{Figure:unit_cells} under affine Dirichlet boundary conditions. The elasticity parameters define elasticity tensors which exhibit { cubic} symmetry  similar to $\Cmacro$.}
\label{Table:Cmicro}
\center
\begin{tabular}{|c| c| c| c|} \hline   
unit cell & $\lambda_i \, \, [\textrm{GPa}]$ & $\mu_i \, \, [\textrm{GPa}]$ & $\mu^*_i \, \, [\textrm{GPa}] $ \\  \hline   
1 & $18.26 \, \, $ & $15.34  \, $  & $\bf 14.61\, $ \\  \hline   
2 & $\bf 20.15 \, \, $ & $\bf 15.83 \, \, $  & $ 14.44 \,  $ \\  \hline   
3 & $19.25 \,  \, $ & $15.54 \, \, $  & $ 13.19 \,  $ \\  \hline   
4 & $19.56 \,  \, $ & $15.66 \,  \, $  & $ 12.68 \,  $ \\  \hline   
\end{tabular}
\end{table}

 The micro elasticity tensor $\Cmicro^\textrm{L\"owner}$ is defined then by the L\"owner matrix supremum problem as
\begin{equation}
\tilde{\Bvarepsilon}:  \Cmicro^\textrm{L\"owner} :\tilde{\Bvarepsilon} \ge \tilde{\Bvarepsilon}:  \mathbb{C}^\textrm{D}_i :\tilde{\Bvarepsilon}  \quad \textrm{where} \quad \textrm{i=1,...,4} \quad , \quad \forall \tilde{\Bvarepsilon} \in \textrm{Sym}(3)\,. 
\end{equation}
which turns for the cubic symmetry case to the following one written in Voigt notation  
 \begin{equation}
 \begin{aligned}
& \left( 
\begin{array}{c}
\tilde{\varepsilon}_{11} \\
\tilde{\varepsilon}_{22} \\
2\tilde{\varepsilon}_{12} \\
\end{array}  \right) : 
\left( \begin{array}{c c c}
2 \mu_\textrm{micro}^\textrm{L\"owner}  + \lambda_\textrm{micro}^\textrm{L\"owner}  & \lambda_\textrm{micro}^\textrm{L\"owner}  & 0 \\ 
\lambda_\textrm{micro}^\textrm{L\"owner}  & 2 \mu_\textrm{micro}^\textrm{L\"owner}  + \lambda_\textrm{micro}^\textrm{L\"owner}  & 0 \\
0 & 0 & (\mu^*)_\textrm{micro}^\textrm{L\"owner} 
\end{array}
\right)  :  \left( 
\begin{array}{c}
\tilde{\varepsilon}_{11} \\
\tilde{\varepsilon}_{22} \\
2\tilde{\varepsilon}_{12} \\
\end{array}  \right)  \ge \\ &  \left( 
\begin{array}{c}
\tilde{\varepsilon}_{11} \\
\tilde{\varepsilon}_{22} \\
2\tilde{\varepsilon}_{12} \\
\end{array}  \right) : \left( \begin{array}{c c c}
2 \mu_i + \lambda_i & \lambda_i & 0 \\ 
\lambda_i & 2 \mu_i + \lambda_i & 0 \\
0 & 0 & \mu^*_i
\end{array}
\right) : \left( 
\begin{array}{c}
\tilde{\varepsilon}_{11} \\
\tilde{\varepsilon}_{22} \\
2\tilde{\varepsilon}_{12} \\
\end{array}  \right)   \textrm{for} \quad \textrm{i=1,....,4} \quad ,  \forall  \left( 
\begin{array}{c}
\tilde{\varepsilon}_{11} \\
\tilde{\varepsilon}_{22} \\
\tilde{\varepsilon}_{12} \\
\end{array}  \right) \in \R^3\,. 
 \end{aligned}
\end{equation}

The solution of the previous problem reads 
\begin{equation}
(\mu^*)_\textrm{micro}^\textrm{L\"owner}   \ge \underset{i}{\textrm{max}} \{ \mu^*_i  \} \,,\qquad \mu^\textrm{L\"owner}_\textrm{micro}  \ge \underset{i}{\textrm{max}} \{ \mu_i  \} \,,\qquad \lambda^\textrm{L\"owner} _\textrm{micro} + \mu^\textrm{L\"owner} _\textrm{micro}   \ge \underset{i}{\textrm{max}} \{ \mu_i  + \lambda_i\}\,, 
\end{equation}
for $i=1,..,4$\,. We take therefore (see Table \ref{Table:Cmicro})
\begin{equation}
 \begin{aligned}
(\mu^*)^\textrm{L\"owner} _\textrm{micro} := \mu^*_1  = 14.61 \, \textrm{GPa}\,, \quad \mu^\textrm{L\"owner} _\textrm{micro} :=   \mu_2  = 15.83  \, \textrm{GPa} \,, \\ \lambda^\textrm{L\"owner} _\textrm{micro}  :=  \mu_2  + \lambda_2 -  \mu_\textrm{micro} = 20.15 \, \textrm{GPa}  \, ,  
 \end{aligned}
\end{equation}
and thus
\begin{equation}
{\mathbb{C}}^\textrm{L\"owner}_\textrm{micro} := \left( \begin{array}{c c c}
51.81 & 20.15 & 0 \\ 
20.15 & 51.81 & 0 \\
0 & 0 & 14.61
\end{array} 
\right)   [\textrm{GPa}]\,.
\end{equation}
However, the previous estimate will serve as a lower bound for $\Cmicro$.
In Fig. \ref{Figure:Compare_tensors_D}, we show the size-effect of the fully resolved metamaterial beams and the linear elasticity solutions with different elasticity tensors: I) $\Cmacro$, II)  $\Cmicro^\textrm{L\"owner}$, III) $\Cmatrix$ of the homogeneous isotropic matrix, and IV) $\CVoigt$ which is isotropic and obtained by the equal strain assumption ${\CVoigt = \phi_\textrm{matrix} \Cmatrix + \phi_\textrm{inclusion} \Cinclusion}$  where $\phi_\textrm{matrix}$ and  $\phi_\textrm{inclusion}$ are the volume fractions of the matrix and inclusion, respectively, {which leads to $\lambda_\textrm{Voigt} = 36.77 \, \textrm{GPa}$ and $\mu_\textrm{Voigt} = 18.44 \, \textrm{GPa}$}.  
The calculated value for ${\mathbb{C}}^\textrm{L\"owner}_\textrm{micro}$ is too soft compared to the microstructured beams and even linear elasticity with $\CVoigt$ is softer than the solution of the fully resolved metamaterial beam for $n=1$. {This can be explained by the fact that the typical bending mode, e.g. due to a pure bending moment as in the paper, cannot be mapped with affine Dirichlet Boundary conditions. Here, a "Voigt bound" for higher modes (not for affine deformations) would be required, which, to our knowledge, does not exist. Note that the tensor $\Cmicro$, although appearing in the relaxed  micromorphic model and in the classical micromorphic model, does not have the same meaning in the latter, which is related to the bounded stiffness property of the former.}  {Since $ \Cmatrix$ represents the largest stiffness},  we {may} relate $\Cmicro$ to the matrix stiffness $ \Cmatrix$ and introduce a scalar $\alpha \ge 1$ so that we have $ \Cmicro := \alpha \, {\mathbb{C}}^\textrm{L\"owner}_\textrm{micro} $. We define an upper limit for  $ \Cmicro $ as 
\begin{equation}
\label{eq:int_matrix} 
\tilde{\Bvarepsilon}:  \Cmatrix :\tilde{\Bvarepsilon} \ge \tilde{\Bvarepsilon}: \Cmicro :\tilde{\Bvarepsilon} = \tilde{\Bvarepsilon}: \alpha \, {\mathbb{C}}^\textrm{L\"owner}_\textrm{micro}  :\tilde{\Bvarepsilon}   , \quad \forall \tilde{\Bvarepsilon} \in \textrm{Sym}(3)\,. 
\end{equation}
By introducing Eq. \ref{eq:int_matrix} we keep the anisotropic symmetry property of $\Cmicro$ {while} the elasticity tensor $\Cmatrix$  is isotropic.  We obtain then 
\begin{equation}
\mu^*_\textrm{matrix} = \mu_\textrm{matrix}   \ge \alpha (\mu^*)_\textrm{micro}^\textrm{L\"owner}  \,,\quad \mu_\textrm{matrix}  \ge \alpha \mu_\textrm{micro}^\textrm{L\"owner} 
\,,\quad \lambda_\textrm{matrix} + \mu_\textrm{matrix}   \ge \alpha ( \lambda_\textrm{micro}^\textrm{L\"owner} + \mu_\textrm{micro}^\textrm{L\"owner})\,, 
\end{equation}
which leads to 
\begin{equation}
\alpha \in [ 1, \textrm{min}( \frac{\mu^*_\textrm{matrix}}{(\mu^*)_\textrm{micro}^\textrm{L\"owner}},  \frac{\mu_\textrm{matrix}}{\mu_\textrm{micro}^\textrm{L\"owner}},   \frac{\mu_\textrm{matrix} + \lambda_\textrm{matrix} }{\mu_\textrm{micro}^\textrm{L\"owner} +  \lambda_\textrm{micro}^\textrm{L\"owner} }) ] = [1,1.66] \,.
\end{equation}
Fig. \ref{Figure:Compare_tensors_D} shows that linear elasticity with $\Cmicro = 1.66 \, {\mathbb{C}}_\textrm{micro}^\textrm{L\"owner}$ is stiffer than the fully resolved metamateriel for $n=1$ and therefore it is a valid choice. 
However, assuming $\Cmicro = \Cmatrix$ does not break the extended Neumann's principle. We will investigate later numerically the consequences of the different choices for $\Cmicro$.

 \begin{figure}[ht]
	\unitlength=1mm
	\center
	\begin{picture}(100,60)
		 \put(22,45){\color[rgb]{1,0,0}\makebox(0,0)[lb]{\small\smash{\tiny $n=1$}}}%	 
		 \put(52,40){\color[rgb]{1,0,0}\makebox(0,0)[lb]{\small\smash{\tiny $n=2$}}}%	  	 	
		 \put(86,40){\color[rgb]{1,0,0}\makebox(0,0)[lb]{\small\smash{\tiny$n=5$}}}%	
         \includegraphics[width=0.6 \textwidth]{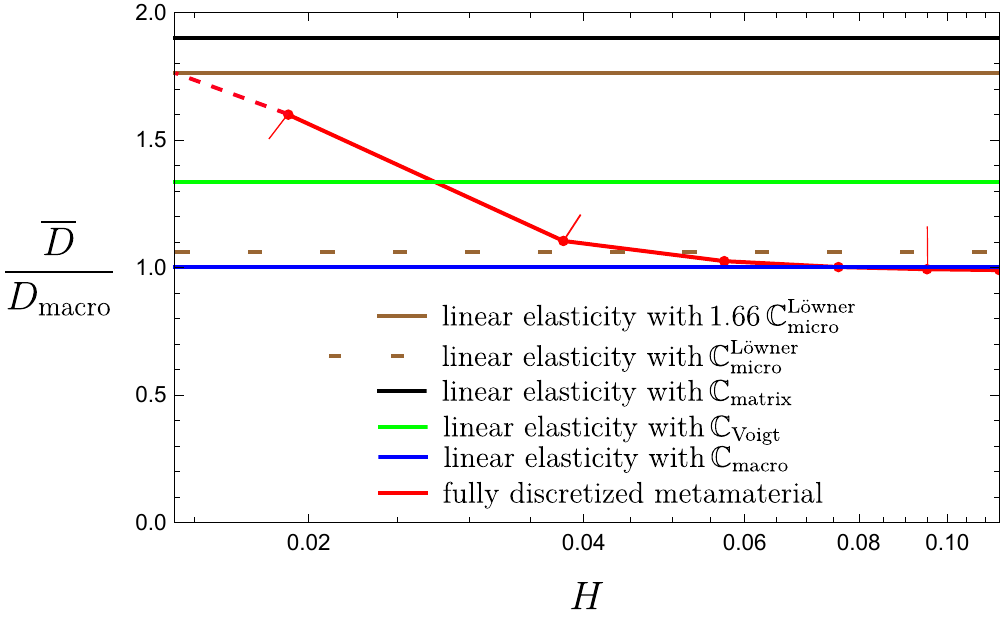}
   	\end{picture}
\caption{The normalized bending stiffness varying the beam size $H \times L = n\,l \times 12 \, n \, l$  compared to the ones obtained by linear elasticity with different elasticity tensors shown in Section  \ref{sec:cmicro1}.}
	\label{Figure:Compare_tensors_D}
\end{figure}

\FloatBarrier
 \subsection{Identification of $\Cmicro$ (second approach)} 
  \label{sec:cmicro2}

  The affine Dirichlet boundary conditions are unable to capture the size-effects as we have shown in Section \ref{sec:cmicro1}. In order to estimate the stiffness $\Cmicro$ for the relaxed micromorphic model we choose in the following approach the most simple ansatz
\begin{equation}
\Cmicro = \beta  \, \Cmacro\, \qquad \textrm{with} \qquad \beta > 1\,.
\end{equation}
In general, the size dependency can not be modeled by a single scalar $\beta$ alone, of course. We introduce this numerical study to get a first estimate. The parameter $\beta$ is determined via the energy equivalence of a heterogeneous microstructured domain and an equivalent homogeneous domain with the same dimensions governed by linear elasticity with elasticity tensor $\Cmicro = \beta \, \Cmacro$, see Fig. \ref{Figure:unit_cells_calculations}. We consider here a higher-order deformation mode which is the bending mode.  The bending mode is enforced by applying non-affine Dirichlet boundary conditions on the whole boundary. They are derived from the analytical solution of the pure bending problem of the homogeneous problem in \cite{RizHueMadNef:2021:aso4} 
\begin{equation}
\bu = \overline\bu = \frac{\kappa}{2} \left( \begin{array}{c}
- 2 \, x y \\ 
\dfrac{\lambda_\textrm{macro}}{2 \mu_\textrm{macro} + \lambda_\textrm{macro}}  y^2 + x^2
\end{array} 
\right) \quad \textrm{on} \quad \partial \B \,,
\end{equation}
which leads to a constant curvature $\kappa$ for the homogeneous case with no shear strain and one active stress component $\sigma_{11}$ 
\begin{equation}
\overline\Bvarepsilon = \sym \nabla \overline \bu =  \left( \begin{array}{c c}
- \kappa \, y & 0 \\
 0 & \dfrac{\lambda_\textrm{macro} \kappa \, y}{2 \mu_\textrm{macro} + \lambda_\textrm{macro}} 
\end{array} \right) \,,  \qquad  \overline\Bsigma = \left(  \begin{array}{c c}
\dfrac{-4  \mu_\textrm{macro} ( \mu_\textrm{macro} + \lambda_\textrm{macro} ) \beta \kappa \, y}{ 2 \mu_\textrm{macro} + \lambda_\textrm{macro}} & 0 \\
0 & 0 \\ 
\end{array} \right) \,. 
\end{equation}

 \begin{figure}[htb]
\center
	\unitlength=1mm
	\begin{picture}(120,40)
	\put(0,5){\def\svgwidth{12 cm}{\small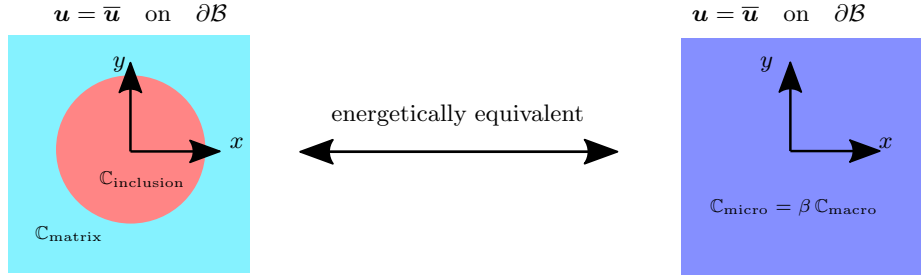}}
	\end{picture}
	\caption{ Illustration shows the procedure used to calculate $\beta$. } 
	\label{Figure:unit_cells_calculations}
\end{figure}

We search for the stiffest possible component on the microstructure under flexural deformation mode (highest values of $\beta$) by investigating different arrangements of unit cells. Six different arrangements were considered, see Fig. \ref{Figure:unit_cells_bending}. The largest obtained value is $\beta = 1.64$. Increasing the size of the  arrangements of the unit cells, considered in Fig. \ref{Figure:unit_cells_bending}, we retrieve the macro property where $\beta$ converges to the value one as it should. This behavior is shown examplarily for unit cell (a) in Fig. \ref{Figure:unit_cell_a_beta_nxn}.

 \begin{figure}[ht]
	\unitlength=1mm
	\center
 \begin{subfigure}[b]{0.32 \textwidth}
		 \center
	\begin{picture}(20,20)
	\put(0,0){\def\svgwidth{2 cm}{\small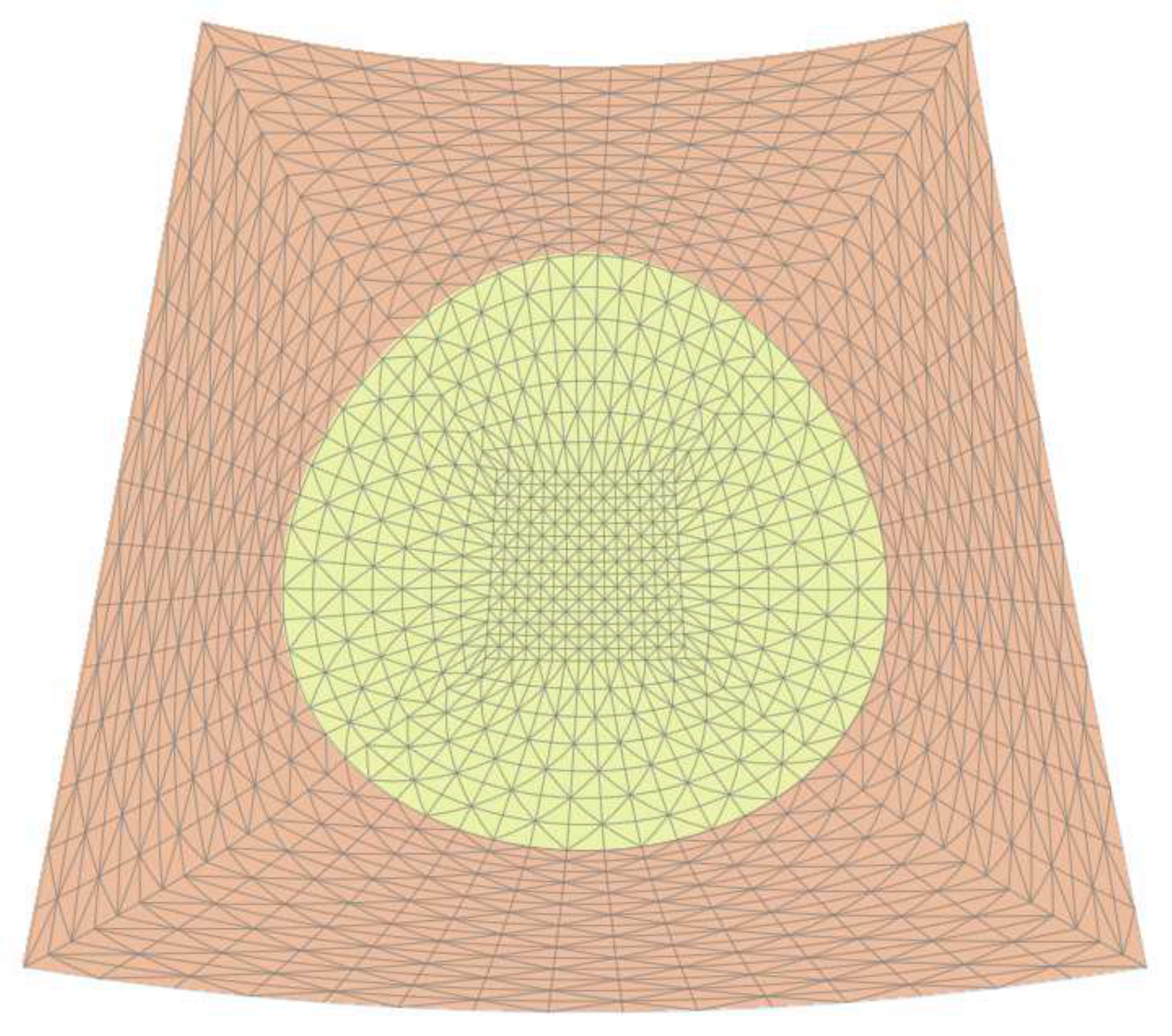}}
	\end{picture}
\caption*{(a)}
    \end{subfigure}
\begin{subfigure}[b]{0.32 \textwidth}
		 \center
	\begin{picture}(20,20)
	\put(0,0){\def\svgwidth{2 cm}{\small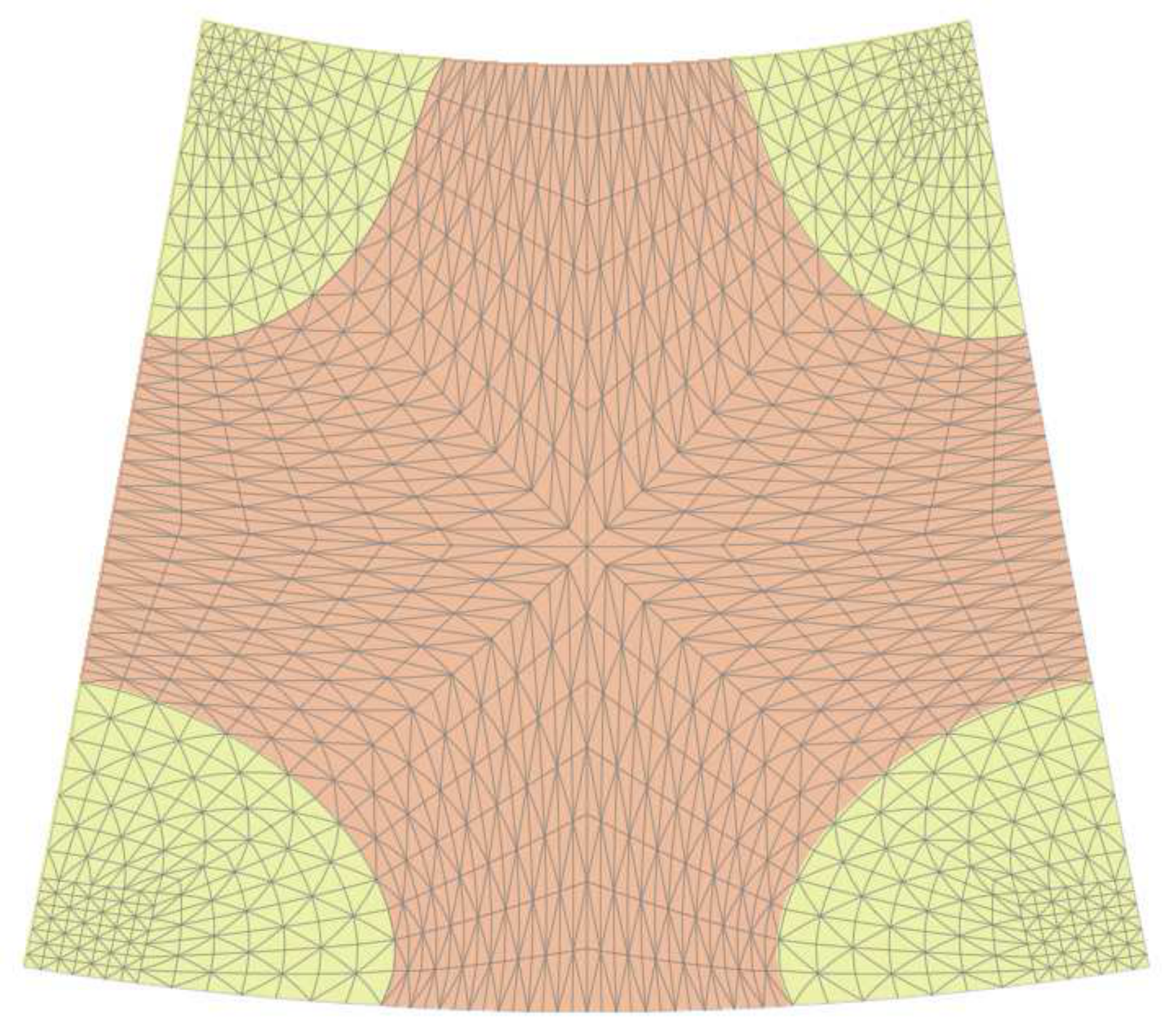}}
	\end{picture}
\caption*{(b)}
\end{subfigure}
\begin{subfigure}[b]{0.32 \textwidth}
		 \center
    	\begin{picture}(20,20)
	\put(0,0){\def\svgwidth{2 cm}{\small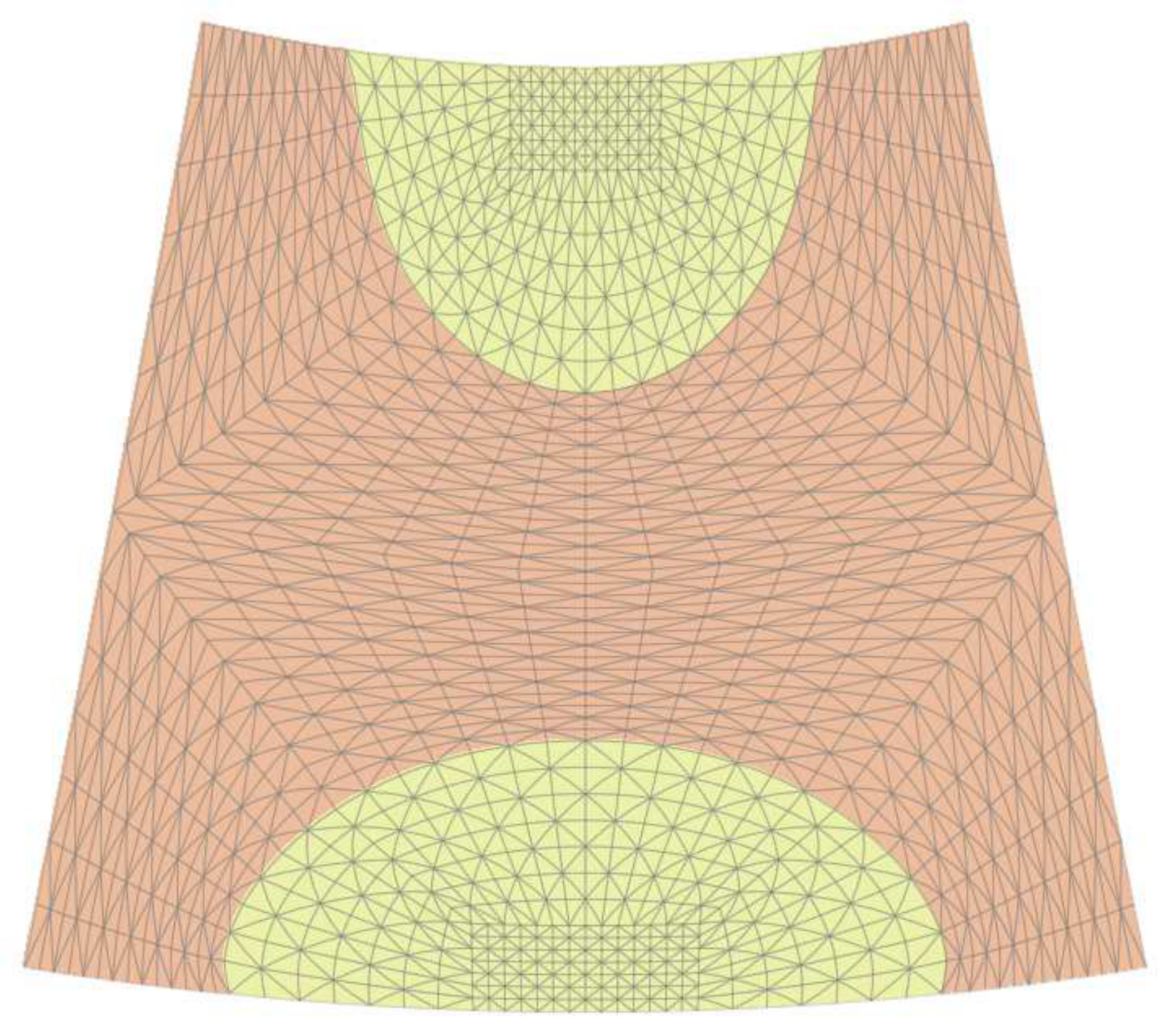}}
	\end{picture}
\caption*{(c)}
    \end{subfigure}
\begin{subfigure}[b]{0.32 \textwidth}
		 \center
    	\begin{picture}(20,35)
	\put(0,0){\def\svgwidth{2 cm}{\small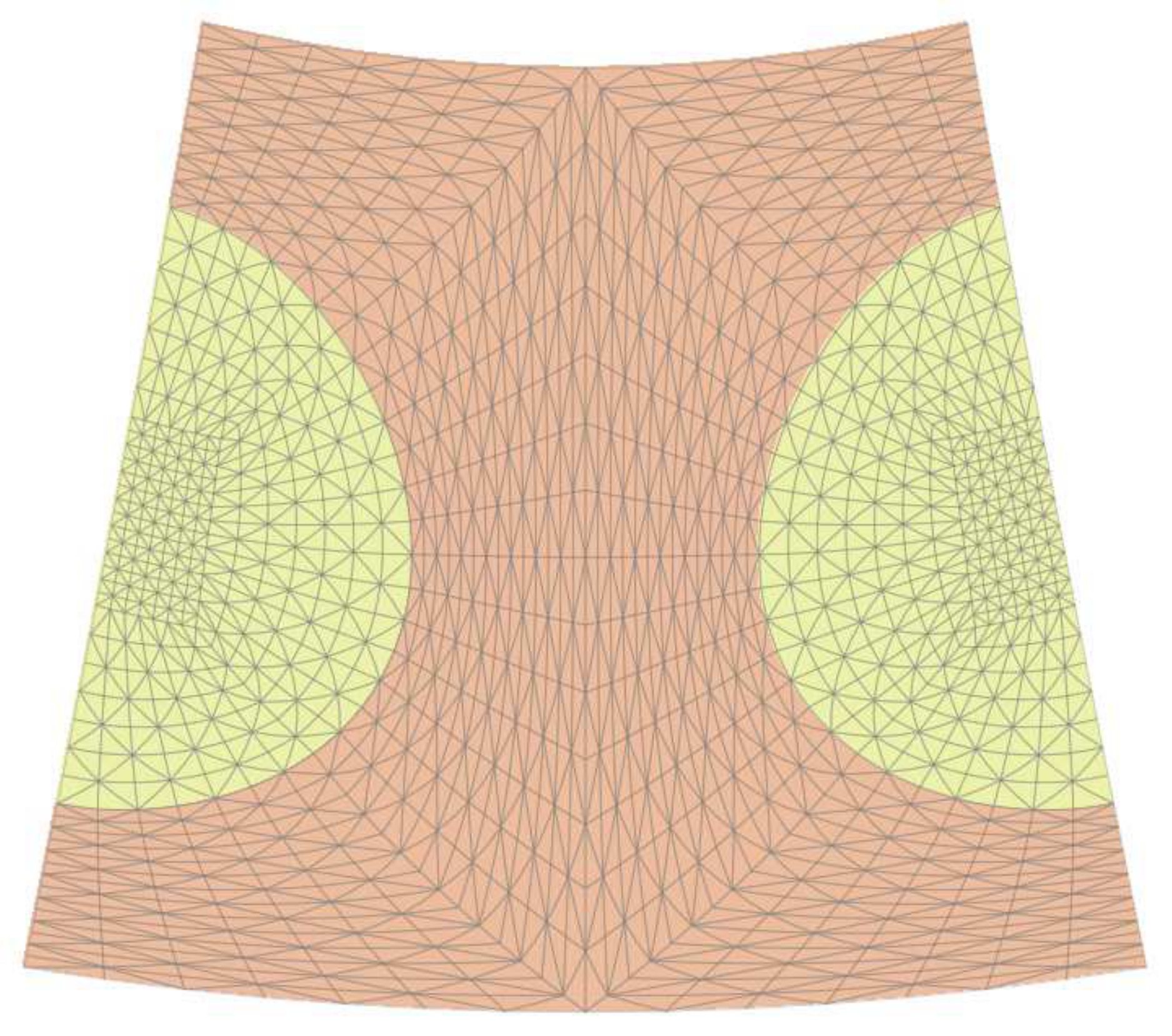}}
	\end{picture}
\caption*{(d)}
    \end{subfigure}
 \begin{subfigure}[b]{0.32 \textwidth}
		 \center
  	\begin{picture}(35,35)
	\put(0,0){\def\svgwidth{3.5 cm}{\small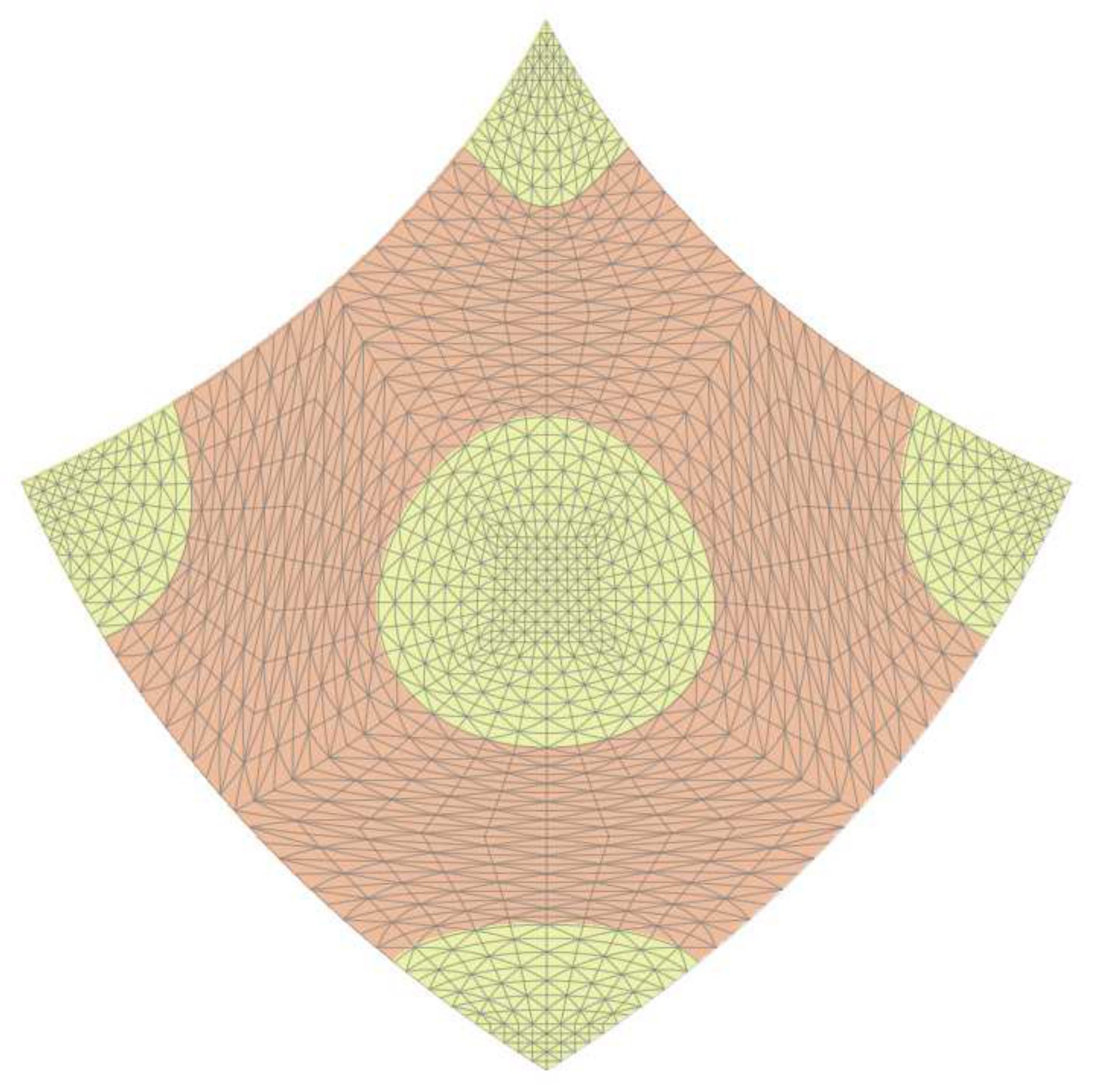}}
	\end{picture}
\caption*{(e)}
    \end{subfigure}
    		 \begin{subfigure}[b]{0.32 \textwidth}
		 \center
	\begin{picture}(35,35)
	\put(0,0){\def\svgwidth{3.5  cm}{\small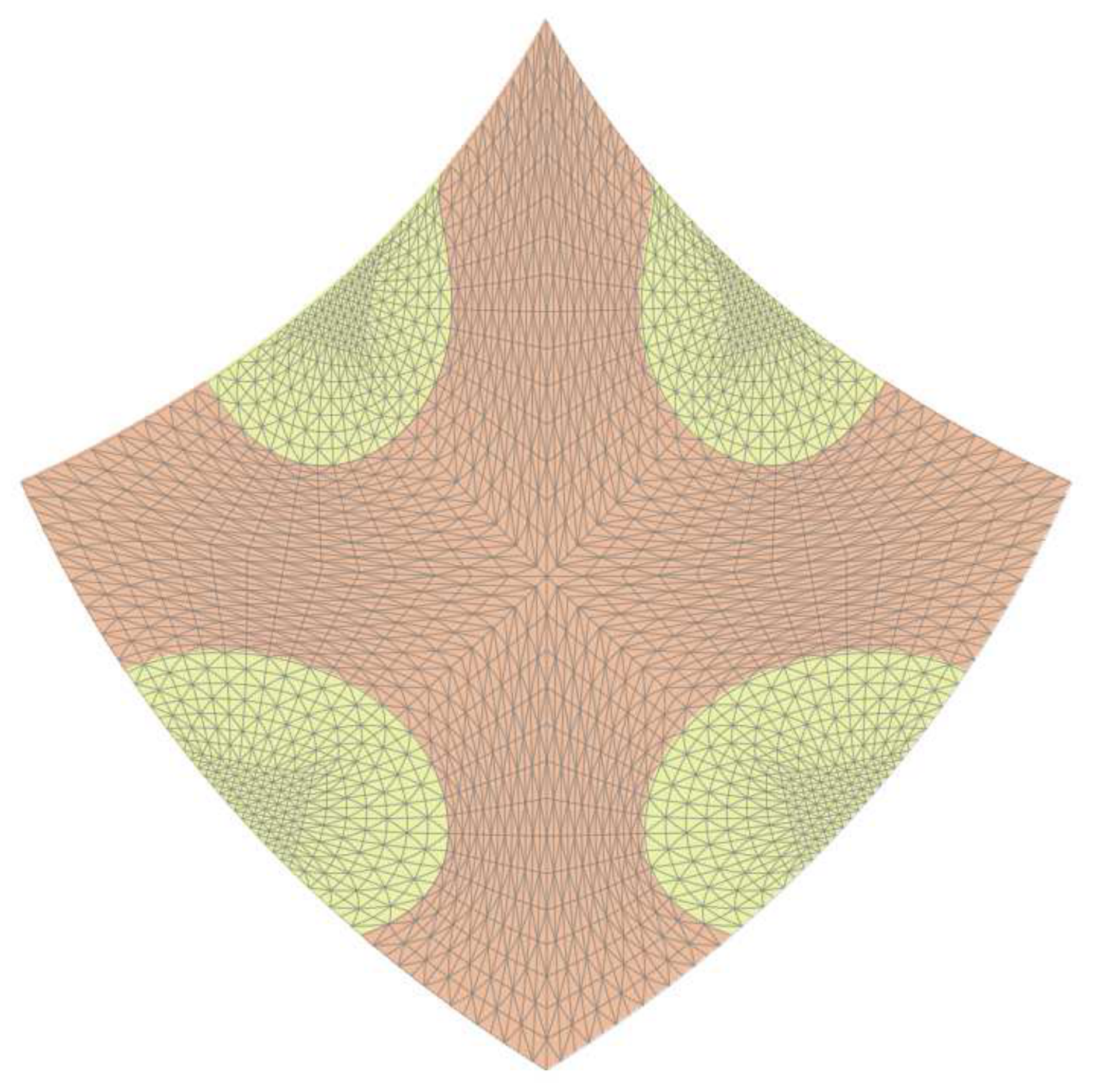}}
	\end{picture}
   		  	\caption*{(f)}
    \end{subfigure}
\caption{ The values of the parameter $\beta$ calculated for different unit cells.  Unit cell (a) provides the stiffest flexural stiffness with $\beta = 1.64$.   }
	\label{Figure:unit_cells_bending}
\end{figure}

 \begin{figure}[ht]
	\unitlength=1mm
	\center
	\begin{picture}(90,65)	
  \put(-10,0){\includegraphics[width=0.6 \textwidth]{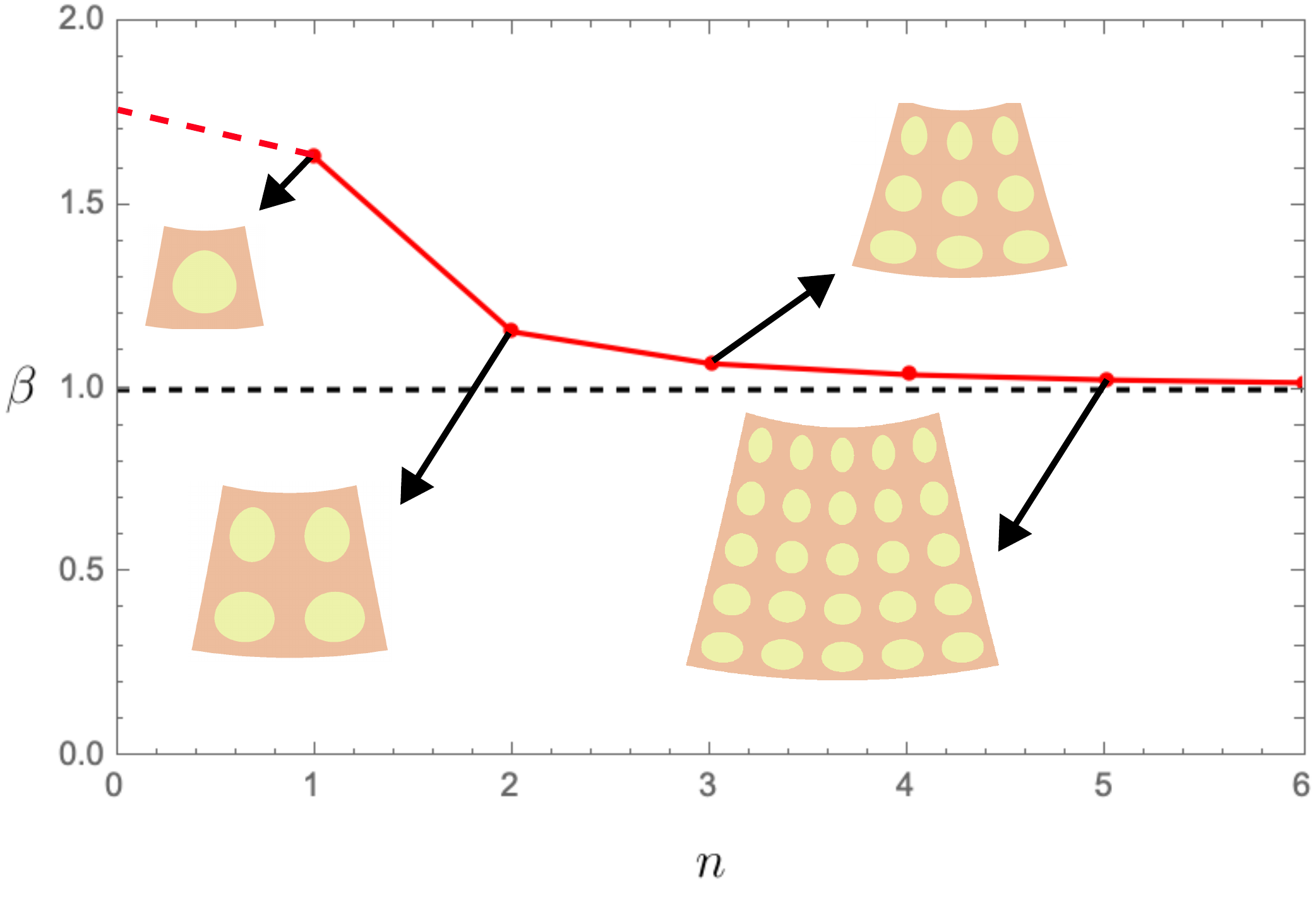}}
  		 \put(0,39){\color[rgb]{0,0,0}\makebox(0,0)[lb]{\small\smash{\small $\beta=1.64$}}}%	 
  		 \put(5,15){\color[rgb]{0,0,0}\makebox(0,0)[lb]{\small\smash{\small $\beta=1.16$}}}%	
		 \put(54,60){\color[rgb]{0,0,0}\makebox(0,0)[lb]{\small\smash{\small $\beta=1.07$}}}%	  	 	
		 \put(46,13){\color[rgb]{0,0,0}\makebox(0,0)[lb]{\small\smash{\small $\beta=1.03$}}}%	
		 \put(-8,56.5){\color[rgb]{1,0,0}\makebox(0,0)[lb]{\small\smash{\small $1.75$}}}%	
	\end{picture}         
\caption{ The parameter $\beta$ converges to the value one when increasing the size of a cluster of unit cells ($n \times n$) shown exemplarily for type (a) in Fig. \ref{Figure:unit_cells_bending}. We also show the extrapolated value $\beta = 1.75$ .}
	\label{Figure:unit_cell_a_beta_nxn}
\end{figure}
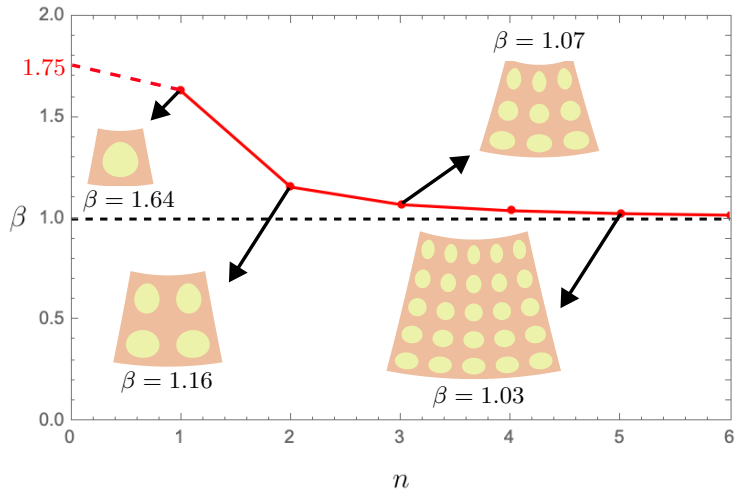

 The choice $\Cmicro = 1.64 \, \Cmacro$ guarantees that a homogeneous continuum with elasticity tensor ${\Cmicro = 1.64 \, \Cmacro}$ is stiffer than the fully discretized metamaterial. In Fig. \ref{Figure:Compare_tensors_D_1}, we show the size-effect of the fully resolved metamaterial beams and the linear elasticity
solutions with elasticity tensors ${\Cmicro = 1.64 \, \Cmacro}$ and $\Cmacro$. The upper limit ${\Cmicro = 1.64 \, \Cmacro}$
is slightly stiffer than the relatively stiffest metamaterial beam ($n = 1$) which confirms its validity. However, to provide a better fitting, we extrapolate $\Cmicro=1.75 \, \Cmacro$ as an improved upper bound. Unique identification of the micro elasticity tensor $\Cmicro$ remains an open question for future works.

 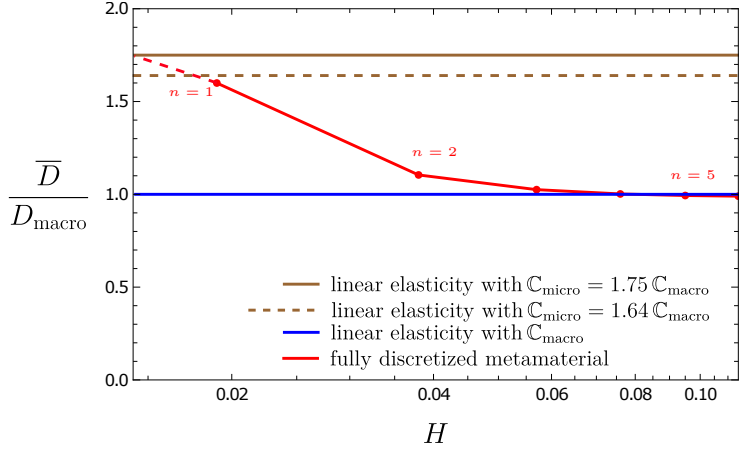
\begin{figure}[ht]
	\unitlength=1mm
	\center
	\begin{picture}(90,65)	
		 \put(11,46){\color[rgb]{1,0,0}\makebox(0,0)[lb]{\small\smash{\tiny $n =1$}}}%	 
		 \put(43,38){\color[rgb]{1,0,0}\makebox(0,0)[lb]{\small\smash{\tiny $n=2$}}}%	  	 	
		 \put(77,35){\color[rgb]{1,0,0}\makebox(0,0)[lb]{\small\smash{\tiny$n=5$}}}%	
  \put(-10,0){\includegraphics[width=0.6 \textwidth]{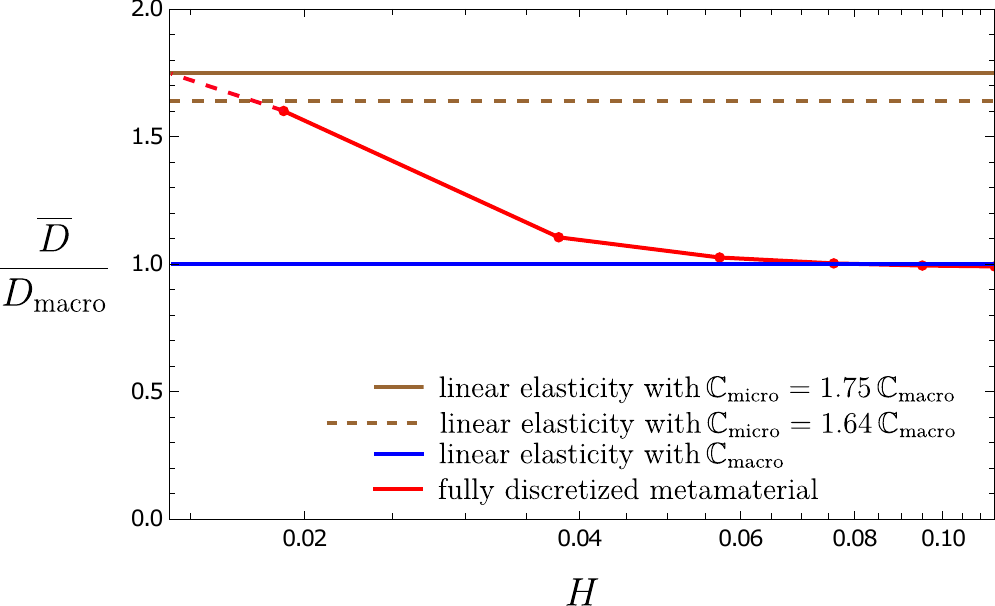}}
	\end{picture}         
\caption{ The normalized bending stiffness  varying the beam size $H \times L = n\,l \times 12 \, n \, l$  compared to the ones obtained by linear elasticity with different elasticity tensors shown in Section \ref{sec:cmicro2}.}
	\label{Figure:Compare_tensors_D_1}
\end{figure}

 \subsection{Identification of $\Ce$} 
The elasticity tensor $\Ce$ is calculated via the micro-macro Reuss-like homogenization formula 
\begin{equation}
\Cmacro^{-1} =  \Cmicro^{-1} +  \Ce^{-1} \quad \Rightarrow \quad \Ce = \Cmicro (\Cmicro-\Cmacro)^{-1} \Cmacro\,. 
\end{equation} 
The  obtained elasticity tensor $\Ce$ is automatically positive definite since $\Cmicro > \Cmacro$ and has cubic symmetry property. However, no obvious physical interpretation can be assigned to the tensor $\Ce$.

\FloatBarrier
\subsection{The boundary conditions of the micro-distortion field}
 The boundary conditions (BCs) of the micro-distortion field are key components for the relaxed micromorphic model. The boundary conditions should be chosen in a way that induces a curvature in the model, i.e. $\Curl \Bdis \ne \bzero$. Otherwise, insufficient boundary condition of the micro-distortion field can cause {unwanted} behavior of the relaxed micromorphic model. This behavior is represented by showing no size-effects or not reaching the intended upper limit (linear elasticity with $\Cmicro$) for $\Lc \rightarrow \infty$. 
 
\subsubsection{Symmetric  force stress case:} 
  \label{sec:symm_for}
 We assume here $\Cc=0$ which causes symmetric force stress $\Bsigma$ and symmetric $\Curl \bbm$ because Eq. \ref{eq:sfs1}d becomes symmetric.  We test the sufficiency of the boundary condition by comparing the solution of the relaxed micromorphic model for varied values of the characteristic  length with the solutions obtained by the standard linear elasticity with elasticity tensors $\Cmicro$ and $\Cmacro$. More specifically, the relaxed micromorphic model should reproduce linear elasticity with elasticity tensors $\Cmicro$ and $\Cmacro$ for $\Lc \rightarrow \infty$ and $\Lc \rightarrow 0$, respectively, see \cite{DagRizKhaLewMadNef:2021:tcc,RizHueKhaGhiMadNef:2021:aso1,RizHueMadNef:2021:aso3,RizHueMadNef:2021:aso4,RizKhaGhiMadNef:2021:aso2,SchSarSchNef:2022:lhb}. 
 
 We design a test by fixing the geometry $H \times L= 2 \, l \times 24 \, l$ with assuming ${\Cmicro = 1.75 \,  {\mathbb{C}}_\textrm{macro}}$ and setting $\mu=\mu_\textrm{macro}$. The boundary conditions of the displacement field are taken similar to the ones applied on the fully resolved metamaterials in Fig. \ref{Figure:beam_metamaterials_BCs}.  For the first case with applied rotation, the consistent coupling condition is applied on the right and left edges via a penalty approach, see Fig.  \ref{Fig:rmm:bcs:case}. Indeed, applying the consistent coupling condition on the Dirichlet boundary of the displacement field is adequate to fulfill the theoretical limits of the relaxed micromorphic model.  Removing the consistent coupling condition on left or right edges leads to vanishing curvature and turns the relaxed micromorphic model into standard linear elasticity with $\Cmacro$.  The previous behavior is demonstrated in Fig.  \ref{Fig:rmm:bcs:case2:Pi}.  The exact same behavior is observed for the second loading case with applied traction if we apply the consistent coupling condition on the boundary corresponding to the first loading case, see Fig. \ref{Fig:rmm:bcs:case}. Consequently, the relaxed micromorphic model results in consistent results for both loading cases, see Fig. \ref{Fig:rmm:bcs:case2:Pi}. 
 
 \begin{figure}[ht]
\center
	\unitlength=1mm
	\begin{picture}(150,20)
	\put(15,0){\def\svgwidth{4 cm}{\small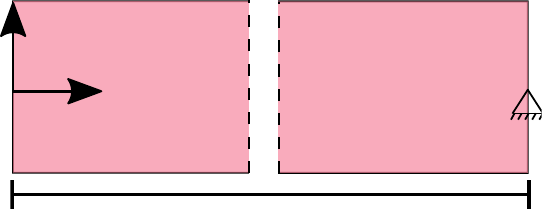}}
	\put(100,0){\def\svgwidth{4 cm}{\small\input{figures/rmm_case2_geo.eps_tex}}}
	\end{picture}
\caption{{The boundary value problems of the homogeneous relaxed micromorphic model for both loading cases. These boundary value problems are equivalent to the two cases of the heterogeneous metamaterial shown in Fig. \ref{Figure:beam_metamaterials_BCs}.  The upper and lower edges are traction-free.} }
	\label{Fig:rmm:bcs:case}
\end{figure} 

  \begin{figure}[ht]
  \center
         \includegraphics[width=0.5 \textwidth]{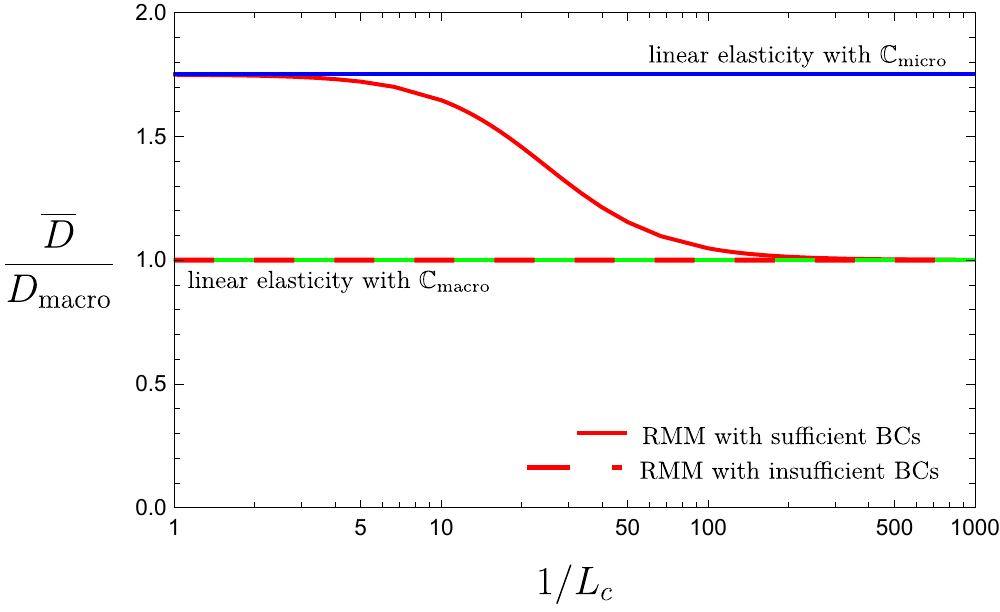}
   		  	\caption{The normalized bending stiffness obtained by the relaxed micromorphic model for both loading cases assuming $\Cc =\bzero$ ($\mu_c=0$) while varying the characteristic length $\Lc$. Sufficient BCs indicate to apply the consistent coupling condition on the left and right edges, see Fig. \ref{Fig:rmm:bcs:case}. Removing the consistent coupling condition on left or right edge is considered as insufficient and leads to no size-effect. }
\label{Fig:rmm:bcs:case2:Pi}
\end{figure} 
\FloatBarrier
\subsubsection{Non-symmetric  force stress case:} 
Here, we assume $\Cc = 2 \mu_c \II$ where $\II$ is the fourth order identity tensor and $\mu_c$ is the  Cosserat couple modulus acting as a spring constant between the skew-symmetric parts of $\nabla \bu$ and $\Bdis$. We study the influence of varying the Cosserat couple modulus $\mu_c \in [0,0.01,0.1,1] * \mu_\textrm{macro} $ considering different scenarios of the boundary condition of $\Bdis$. The geometry and the remaining material parameters are taken as for the symmetric case, see Section~ \ref{sec:symm_for}. 

In Fig. \ref{Fig:rmm:bcs1:case2:Pi}, we show the normalized bending stiffness for the cases (a) the consistent coupling condition is applied either on the left or right edge, (b) no consistent boundary condition is considered and (c) the consistent coupling condition is applied on both left and right edges.  Size-effects are noticed even if the consistent coupling condition is not placed on the right and left edges simultaneously which is not the case for the symmetric force case ($\mu_c = 0$). Increasing the Cosserat couple modulus  $\mu_c$ raises the stiffness of the relaxed micromorphic model closer to linear elasticity with $\Cmicro$ for $\Lc \rightarrow \infty$ and even reach it in Fig. \ref{Fig:rmm:bcs1:case2:Pi}(a). However, it is not guaranteed that the relaxed micromorphic model achieves linear elasticity with $\Cmicro$ for $\Lc \rightarrow \infty$, see Fig. \ref{Fig:rmm:bcs1:case2:Pi}(b). The results of enforcing the consistent coupling condition on both left and right edges are equivalent for the symmetric and non-symmetric cases in Figs. \ref{Fig:rmm:bcs:case2:Pi} and \ref{Fig:rmm:bcs1:case2:Pi}(c), respectively, and the Cosserat couple modulus has no influence. 

 \begin{figure}[ht]
    	\unitlength=1mm
	\center
 		 \begin{subfigure}[b]{0.48 \textwidth}
		 \center
         \includegraphics[width=1 \textwidth]{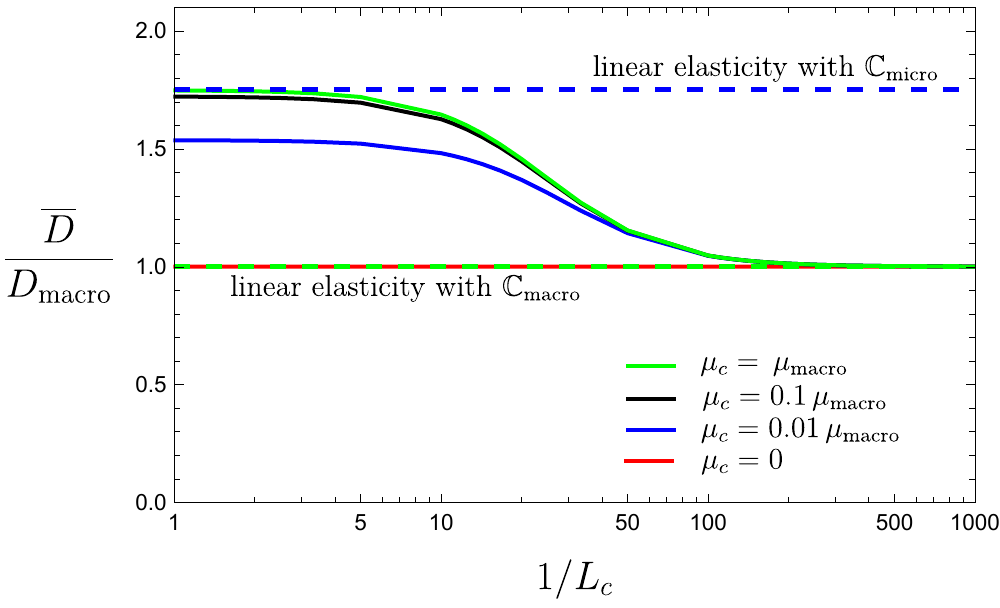}
   		  	\caption{\scriptsize  consistent coupling condition either on the left or right edge}
    \end{subfigure}
  \vspace{0.5 cm}
 \hspace{0.2 cm}
 		 \begin{subfigure}[b]{0.48 \textwidth}
		 \center
         \includegraphics[width= \textwidth]{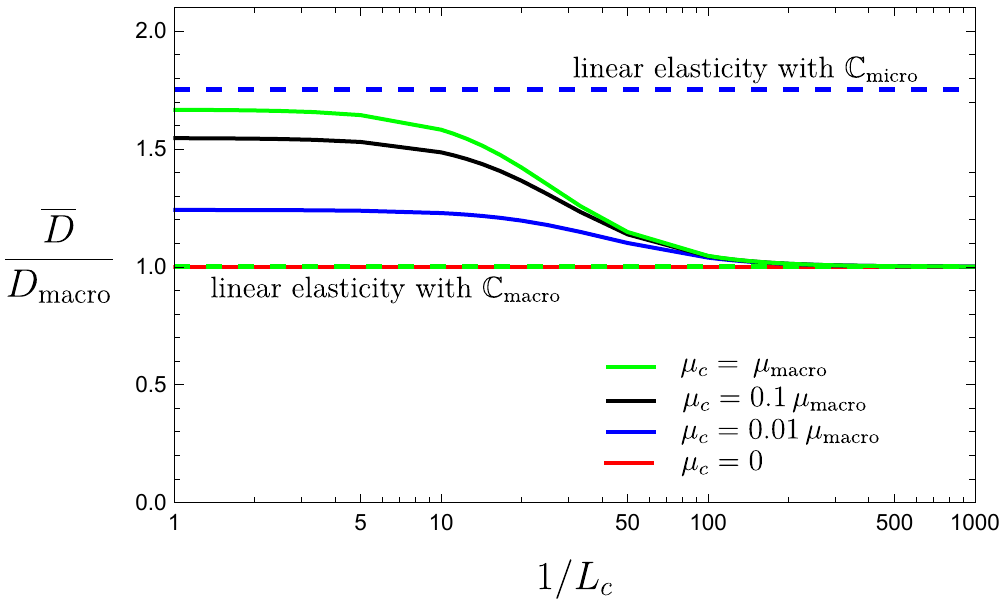}
   		  	\caption{\scriptsize no consistent coupling condition}
    \end{subfigure}
 		 \begin{subfigure}[b]{0.48 \textwidth}
		 \center
         \includegraphics[width= \textwidth]{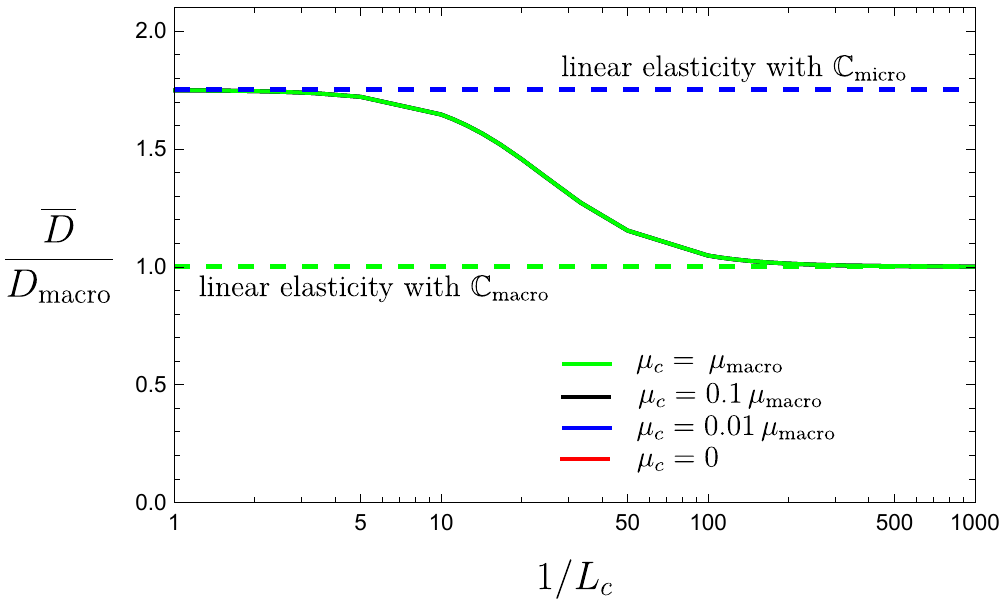}
   		  	\caption{\scriptsize consistent coupling condition on the left and right edges}
    \end{subfigure}

		  	\caption{The normalized bending stiffness obtained by the relaxed micromorphic model for both loading cases with non-symmetric force stress and with varying the characteristic length $\Lc$. Different scenarios are investigated for the boundary conditions of the micro-distortion field.}
\label{Fig:rmm:bcs1:case2:Pi}
\end{figure} 

\subsubsection{Cosserat limit, special case of a skew-symmetric micro-distortion field:}
\label{sec:micro3}
For the case of ${\Cmicro \rightarrow \infty}$, the micro-distortion field $\Bdis$ must be skew-symmetric and the Cosserat model is recovered, c.f.  \cite{AlaGanSadNasAkb:2022:cmo,BleNef:2022:ssi,GhiRizMadNef:2022:cme,KhaChiMadNef:2022:eau}. We investigate here the influence of different scenarios of the boundary conditions for the micro-distortion field $\Bdis$ similar to Section \ref{sec:symm_for}: (a) the consistent coupling condition is applied on either the left or right edge, (b) without enforcing the consistent boundary condition and (c) the consistent coupling condition is applied on both left and right edges for $\Cmicro = 1000 \, \Cmacro$. Different values of the Cosserat couple modulus  $\mu_c$ are assumed for varied values of the characteristic length $\Lc$ in Fig. \ref{Fig:rmm:bcs1:case2:Pi_C}. Our analysis shows that when the consistent coupling condition is not applied at both right and left ends, large values of $\Lc$ result in a beam that does not bend, causing a nonphysical bending stiffness. This highlights the crucial role of the consistent coupling condition, not just in the relaxed micromorphic model, but also in the Cosserat case. Hence due to the bending stiffness issue, we have opted to show the relative total energy $\Pi/\Pi_\textrm{macro}$ for this analysis alternatively.

 \begin{figure}[ht]
    	\unitlength=1mm
	\center
 		 \begin{subfigure}[b]{0.48 \textwidth}
		 \center
         \includegraphics[width=0.7 \textwidth]{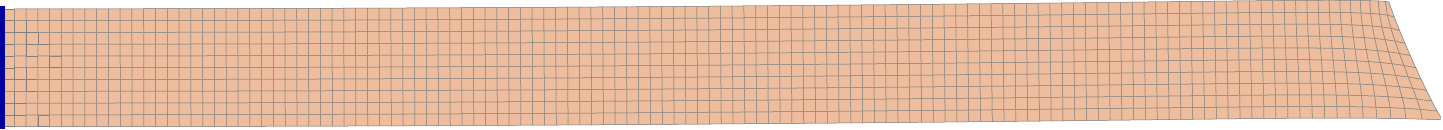}
   		  	\caption{\scriptsize  consistent coupling condition on the left edge}
    \end{subfigure}
 		 \begin{subfigure}[b]{0.48 \textwidth}
		 \center
         \includegraphics[width=0.65 \textwidth]{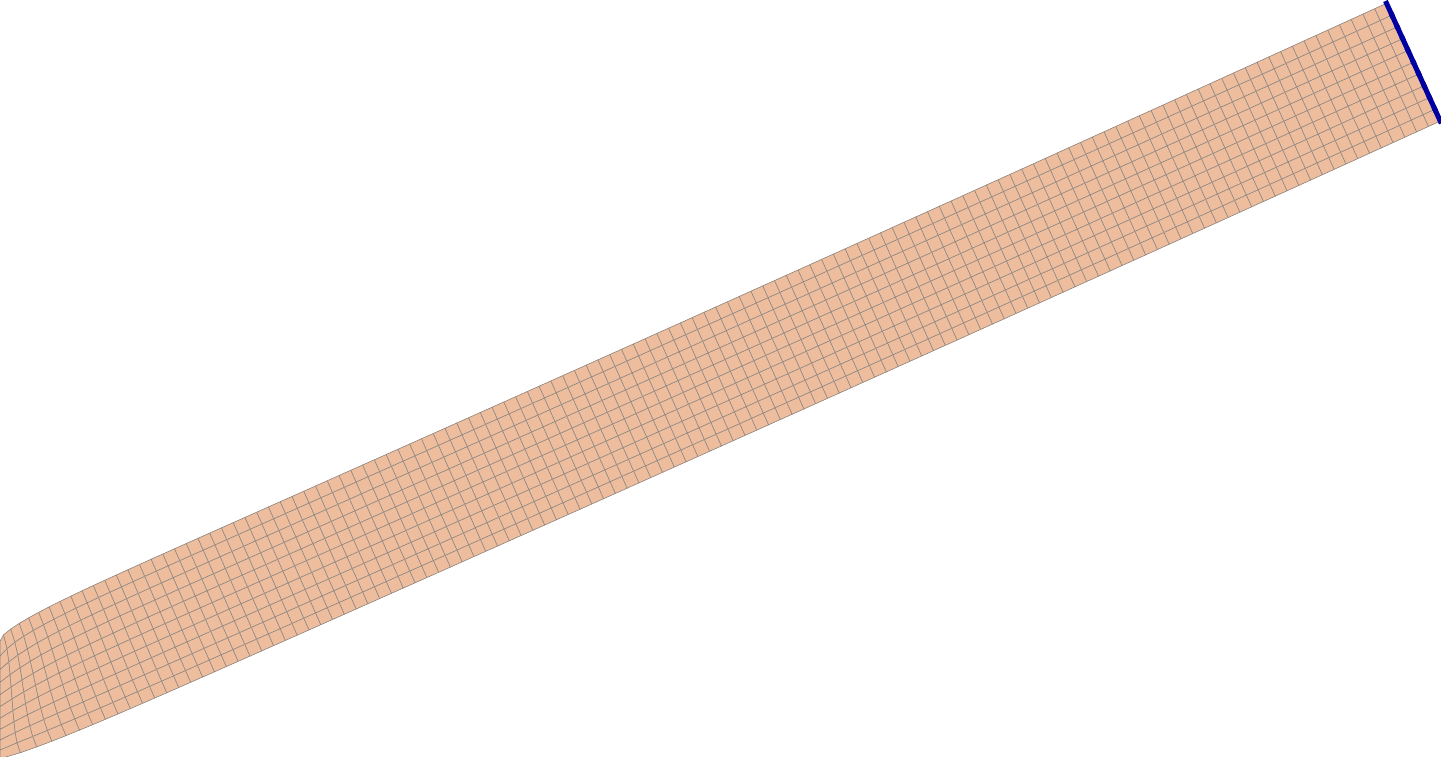}
   		  	\caption{\scriptsize  consistent coupling condition on the right edge}
    \end{subfigure}    
    \vspace{0.5  cm}
 		 \begin{subfigure}[b]{0.48  \textwidth}
		 \center
         \includegraphics[width=0.7 \textwidth]{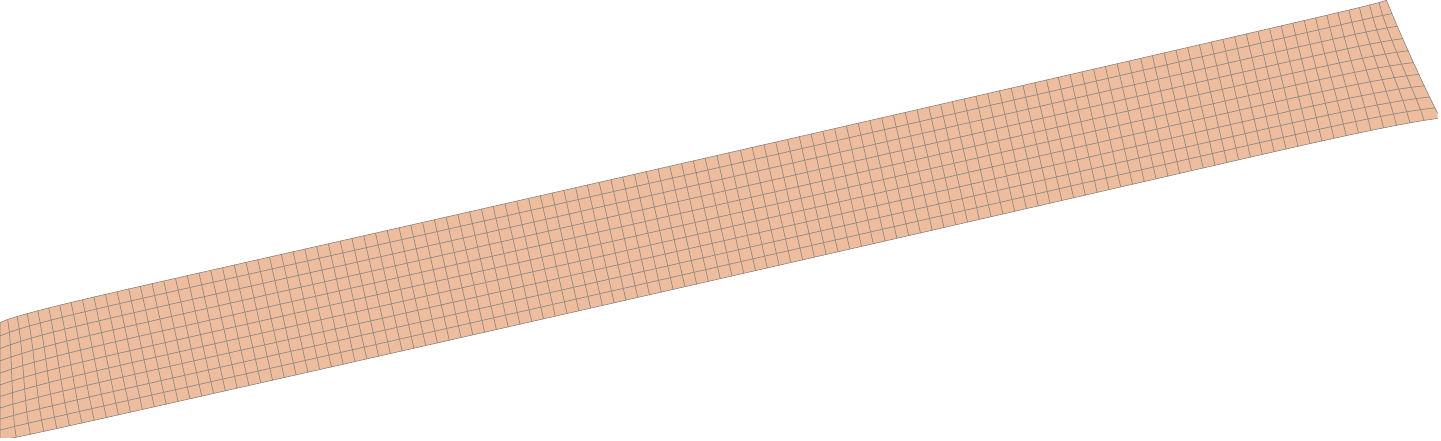}
   		  	\caption{\scriptsize no consistent coupling condition}
    \end{subfigure}
     		 \begin{subfigure}[b]{0.48 \textwidth}
		 \center
         \includegraphics[width=0.7 \textwidth]{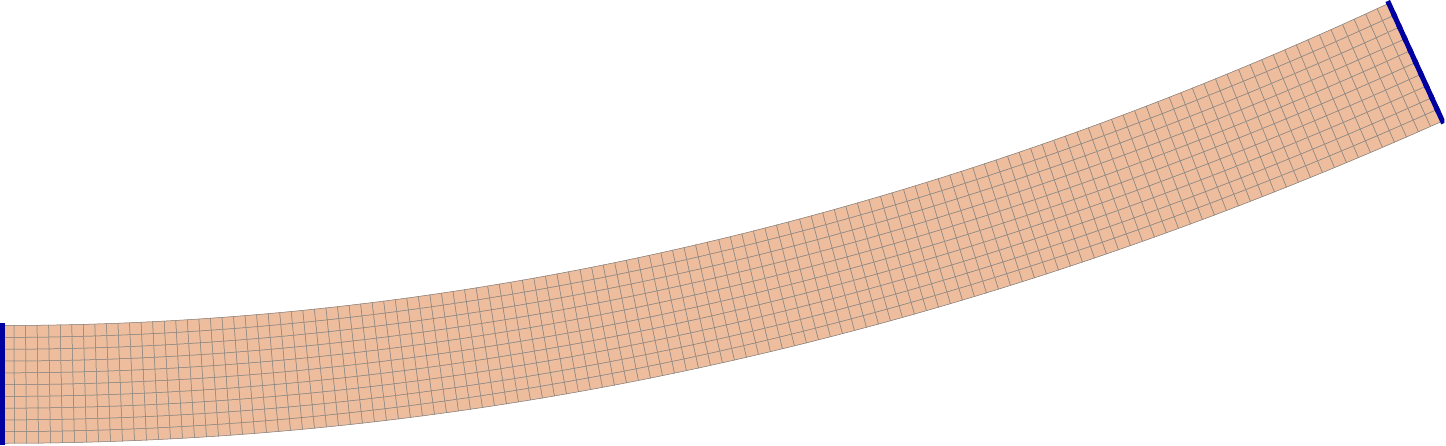}
   		  	\caption{\scriptsize consistent coupling condition on the left and right edges}
    \end{subfigure}

		  	\caption{The deformed beams for the case $\Cmicro = 1000 \, \Cmacro$ "Cosserat type" with $\Lc= \, 1000 \, \textrm{m}$ and $\mu_c = 2 \, \mu_\textrm{macro}$. Bending of the beam can only be induced when the consistent coupling condition is applied on both its left and right ends. }
\end{figure} 

 We notice that linear elasticity with elasticity tensor $\Cmicro$ is recognized as an upper limit only when the consistent coupling condition is enforced on both left and right edges in Fig. \ref{Fig:rmm:bcs1:case2:Pi_C}(c). Weak size-effects are noticed when the consistent coupling condition is not enforced, Fig. \ref{Fig:rmm:bcs1:case2:Pi_C}(a-b). While size-effects are prompted only for non-vanishing  Cosserat couple modulus $\mu_c \ne 0$ for cases (a) and (b), enforcing the consistent coupling condition on both left and right edges simultaneously allows the model to act on the intended theoretical range with no influence of the Cosserat couple modulus  $\mu_c$ which is well known for the Cosserat model. This can be explained by the fact that the skew-symmetric part of the micro-distortion field is the same as the skew-symmetric part of the gradient of the displacement, see \cite{RizHueMadNef:2021:aso4}, which is the case for both the relaxed micromorphic continuum in Fig. \ref{Fig:rmm:bcs1:case2:Pi}(c) and the Cosserat model in Fig. \ref{Fig:rmm:bcs1:case2:Pi_C}(c).

 \begin{figure}[ht]
    	\unitlength=1mm
	\center
 		 \begin{subfigure}[b]{0.48 \textwidth}
		 \center
         \includegraphics[width=1 \textwidth]{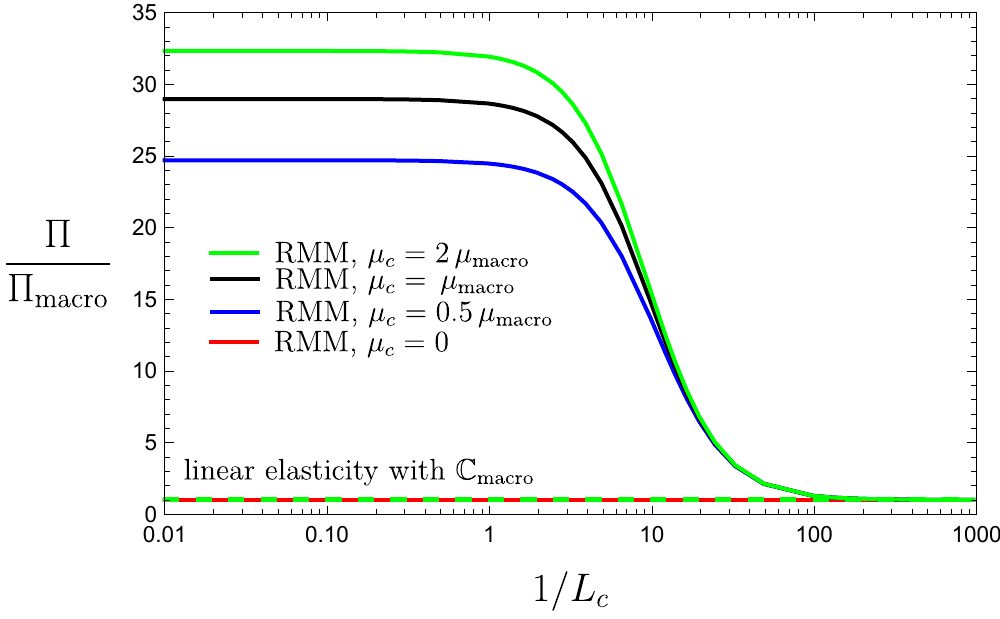}
   		  	\caption{\scriptsize  consistent coupling condition either on the left or right edge}
    \end{subfigure}
    \vspace{0.5  cm}
 		 \begin{subfigure}[b]{0.48  \textwidth}
		 \center
         \includegraphics[width= \textwidth]{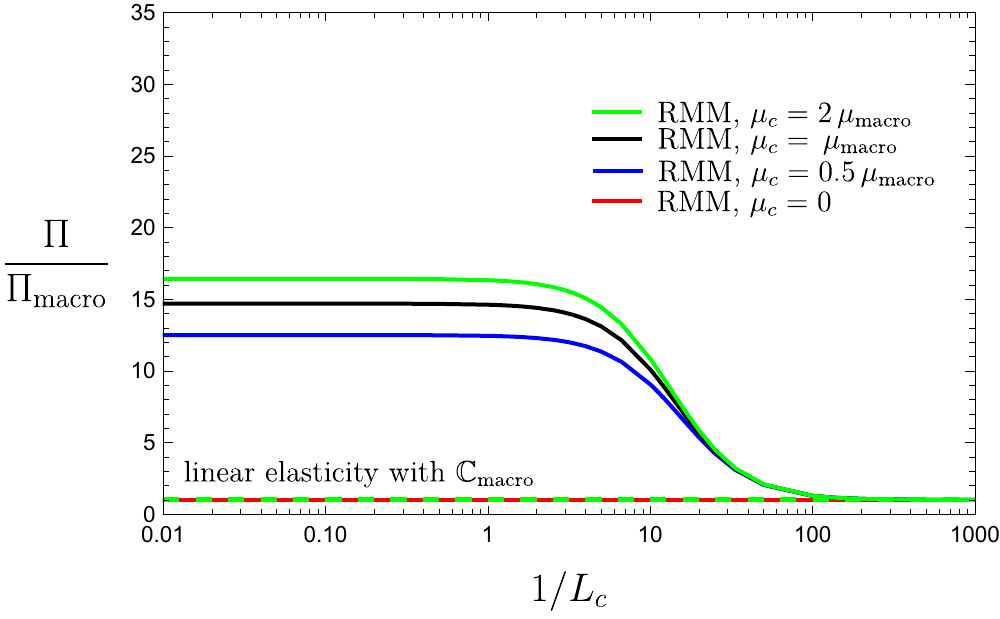}
   		  	\caption{\scriptsize no consistent coupling condition}
    \end{subfigure}
     		 \begin{subfigure}[b]{0.48 \textwidth}
		 \center
         \includegraphics[width= \textwidth]{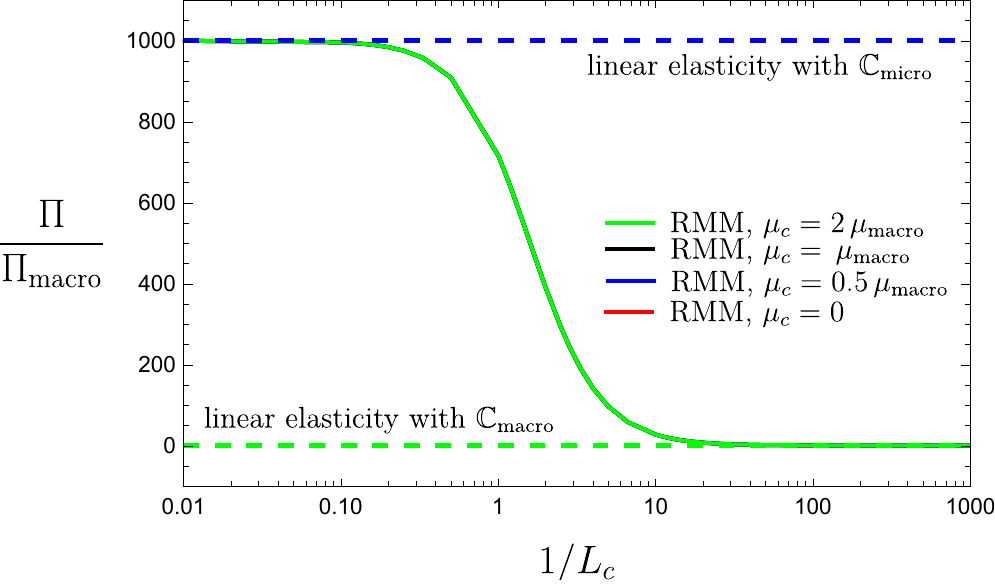}
   		  	\caption{\scriptsize consistent coupling condition on the left and right edges}
    \end{subfigure}

		  	\caption{The relative total energy obtained by the relaxed micromorphic model for both loading cases with non-symmetric force stress and with varying the characteristic length $\Lc$. Here,  we assume ${\Cmicro = 1000 \, \Cmacro}$ leading to a skew-symmetric micro-distortion field which retrieves the Cosserat model since the curvature expression is then equivalent with the Cosserat framework, see \cite{GhiRizMadNef:2022:cme}.  Different scenarios are investigated for the boundary conditions of the micro-distortion field.}
\label{Fig:rmm:bcs1:case2:Pi_C}
\end{figure}

 \FloatBarrier
 \subsection{Scaling of the curvature} 
 
The curvature for the 2D case is isotropic because $\Curl \Bdis$ is reduced to a vector. Therefore, the curvature will be controlled by only one parameter with assuming that $\IL = \II$ is the fourth order identity tensor. Since the parameters $\mu$ and $\Lc$ should be set constant independent of the specimen size, the curvature is modified by incorporating the size of the beams through the number $n$.
Fig. \ref{Figure:metamateials_beam_stiffness} shows that stiffer response is observed for smaller values of the number $n$ ($n=1$ is the stiffest). The relaxed micromorphic model exhibits stiffer response for bigger values of the characteristic length $\Lc$ (inversely proportional to $n$), see for example Fig. \ref{Fig:rmm:bcs:case2:Pi}. Therefore, we replace the last term in Eq. \ref{eq:W} by
 \begin{equation}
 \label{eq:cal}
 \frac{1}{2}   \mu  \, {\left(\frac{\Lc}{n}\right)}^2  \textrm{Curl} \Bdis :\textrm{Curl} \Bdis \,,
 \end{equation}
 where $n$ denotes the number of unit cells in the second-direction.
Hence, for a constant $\Lc$ smaller values are obtained for the term $\Lc/n$ by increasing the beam size (increasing $n$) which reproduces the intended size-effects (smaller is stiffer). This modification  is not ad hoc, but follows from a rigorous scaling argument,  c.f. \cite{NefEidMad:2019:ios} { and applies as such to higher-gradient models or the classical micromorphic model as well. Note that the shear modulus $\mu$ appears for dimensional reasons and is a priori not related to the shear moduli appearing in $\Cmacro$ or $\Cmicro$. }

 \FloatBarrier

 \section{Final calibration} 
 \label{sec:fit}
Now, we provide an identification scheme for the scale-independent material parameters of the relaxed micromorphic model. The boundary conditions of the micro-distortion field are determined in order to guarantee the intended behavior of the relaxed micromorphic model and the influence of the characteristic length $\Lc$ for both loading cases. For this calibration we assume symmetric force stress, i.e. $\Cc = \bzero$. As we discussed in Sections \ref{sec:cmicro1},  \ref{sec:cmicro2} and \ref{sec:micro3}, different choices can be made for $\Cmicro$, e.g. $\Cmicro = 1.66 \, \Cmicro^\textrm{L\"owner} $,  $\Cmicro =  \Cmatrix$, $\Cmicro = 1.75 \, \Cmacro$ and  $\Cmicro = 1000 \, \Cmacro$. Considering $\Cmicro=10 000 \, \Cmacro$ yield similar results to $\Cmicro = 1000 \, \Cmacro$, as expected, which can be explained by the fact that we are operating in a range  close to the lower bound $\Cmacro$. For each choice of $\Cmicro$, the curvature should be calibrated by means of $\Lc$ and $\mu$. Without loss of generality,  we can always assume the shear modulus $\mu=\mu_\textrm{macro}$  and then the characteristic length $\Lc$ should be selected in order to capture the size-effects of the fully discretized metamaterial, Fig. \ref{Figure:stiffness_rmm_case1-lc}.  Alternatively, the characteristic length $\Lc$ can be set in advance, e.g. $\Lc=l$, and then the shear modulus $\mu$ should be calibrated, see Fig. \ref{Figure:stiffness_rmm_case1-mu} and Eq. \ref{eq:cal}. The decisive quantity is the product $\mu \Lc^2$. Since the Cosserat curvature coincides with the curvature expression of the relaxed micromorphic model \cite{GhiRizMadNef:2022:cme}, one would expect that using similar values for $\mu \Lc^2$ is a sensible choice. As Figures \ref{Figure:stiffness_rmm_case1-lc}(d) and  \ref{Figure:stiffness_rmm_case1-mu}(d)  show, this is not the case. For a rough Cosserat fit different orders of magnitude for $\mu \Lc^2$ have to be taken which are getting arbitrary. Furthermore, the data points can be fitted also with a Cosserat type model but it should be remarked that the unbounded stiffness (beyond $n=1$) leads to a sensitive identification of the parameters. The same problem would appear by using second gradient or the classical micromorphic theories.

 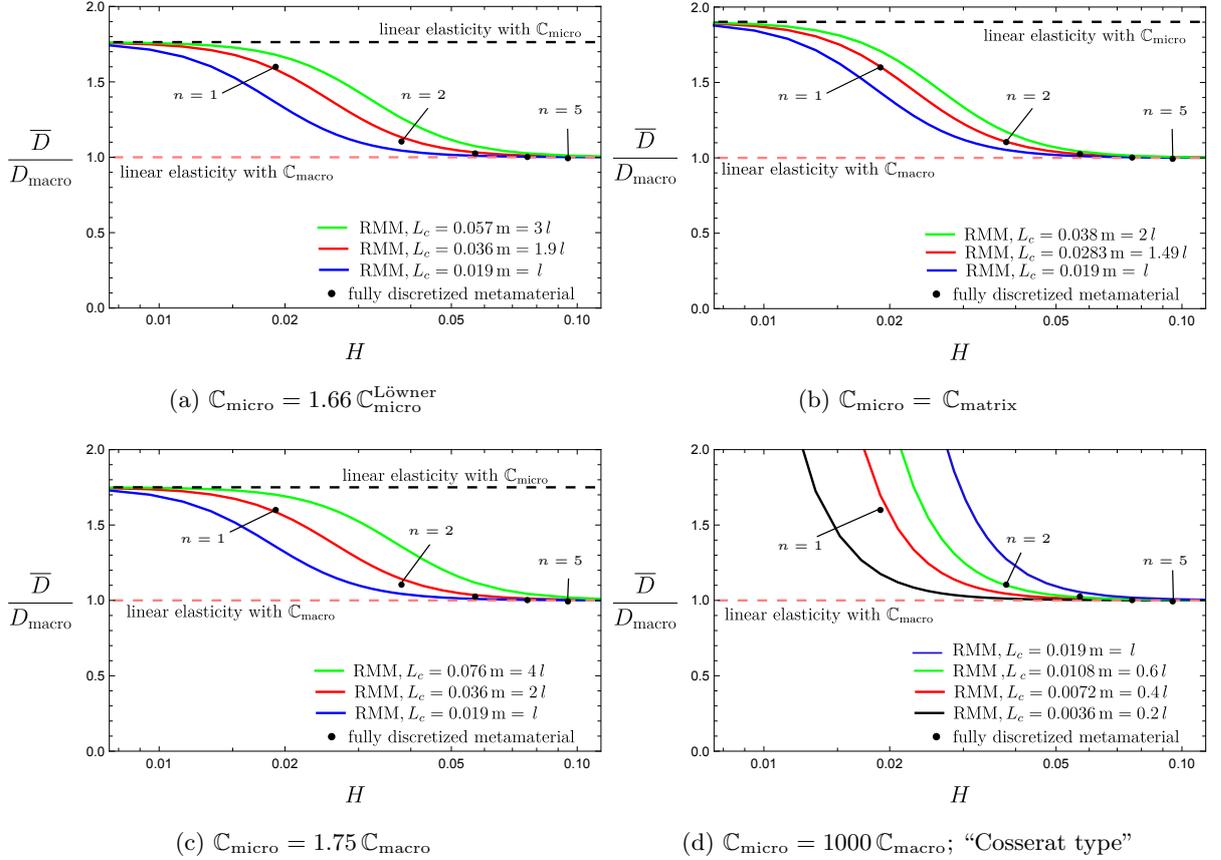
\begin{figure}[ht]
	\unitlength=1mm
	\center
 		 \begin{subfigure}[b]{0.49 \textwidth}
		 \center
		 	\begin{picture}(100,47)
		 \put(22,36){\color[rgb]{0,0,0}\makebox(0,0)[lb]{\small\smash{\tiny $n=1$}}}%	 
		 \put(52,36){\color[rgb]{0,0,0}\makebox(0,0)[lb]{\small\smash{\tiny $n=2$}}}%	  	 	
		 \put(70,34){\color[rgb]{0,0,0}\makebox(0,0)[lb]{\small\smash{\tiny$n=5$}}}%	
         \includegraphics[width=1 \textwidth]{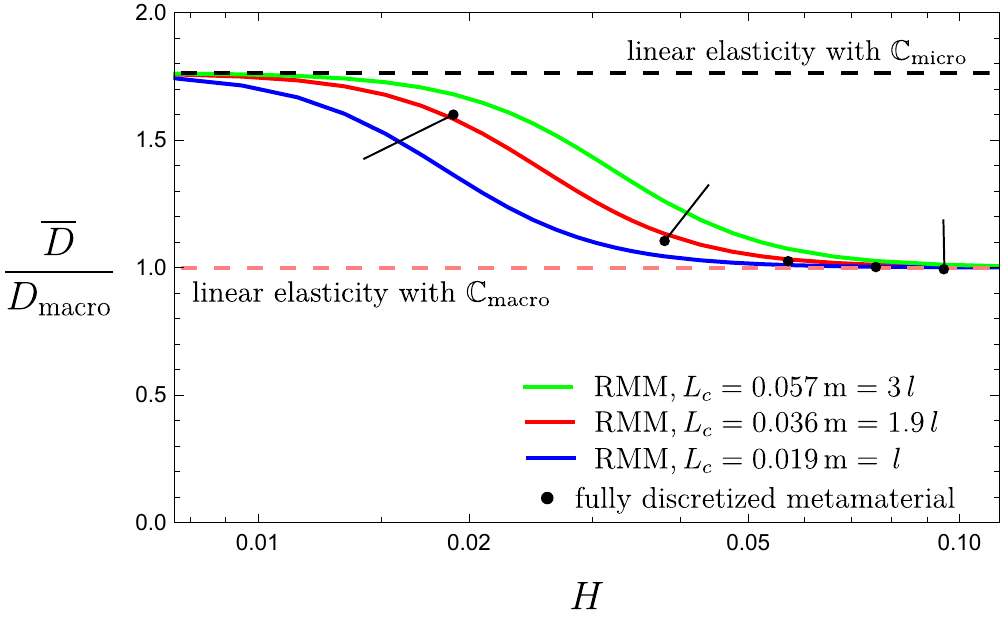}
	\end{picture}
   		  	\caption{$\Cmicro=  1.66 \, {\mathbb{C}}^\textrm{L\"owner}_\textrm{micro} $}
    \end{subfigure}
 		 \begin{subfigure}[b]{0.49 \textwidth}
		 \center
		 	\begin{picture}(100,47)
		 \put(22,36){\color[rgb]{0,0,0}\makebox(0,0)[lb]{\small\smash{\tiny $n=1$}}}%	 
		 \put(52,36){\color[rgb]{0,0,0}\makebox(0,0)[lb]{\small\smash{\tiny $n=2$}}}%	  	 	
		 \put(70,34){\color[rgb]{0,0,0}\makebox(0,0)[lb]{\small\smash{\tiny$n=5$}}}%	
         \includegraphics[width=1 \textwidth]{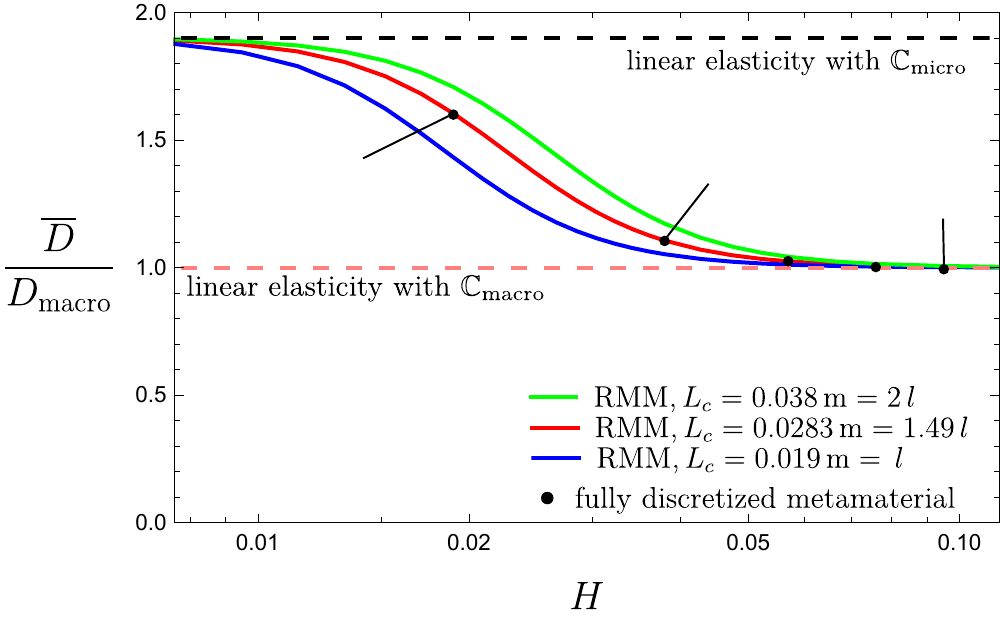}
   	\end{picture}
   		  	\caption{$\Cmicro=  \, \Cmatrix $}
    \end{subfigure}
      		 \begin{subfigure}[b]{0.49 \textwidth}
		 \center
		 	\begin{picture}(100,52)
		 \put(23,36){\color[rgb]{0,0,0}\makebox(0,0)[lb]{\small\smash{\tiny $n=1$}}}%	 
		 \put(53,37){\color[rgb]{0,0,0}\makebox(0,0)[lb]{\small\smash{\tiny $n=2$}}}%	  	 	
		 \put(70,33){\color[rgb]{0,0,0}\makebox(0,0)[lb]{\small\smash{\tiny$n=5$}}}%	
  \includegraphics[width=1 \textwidth]{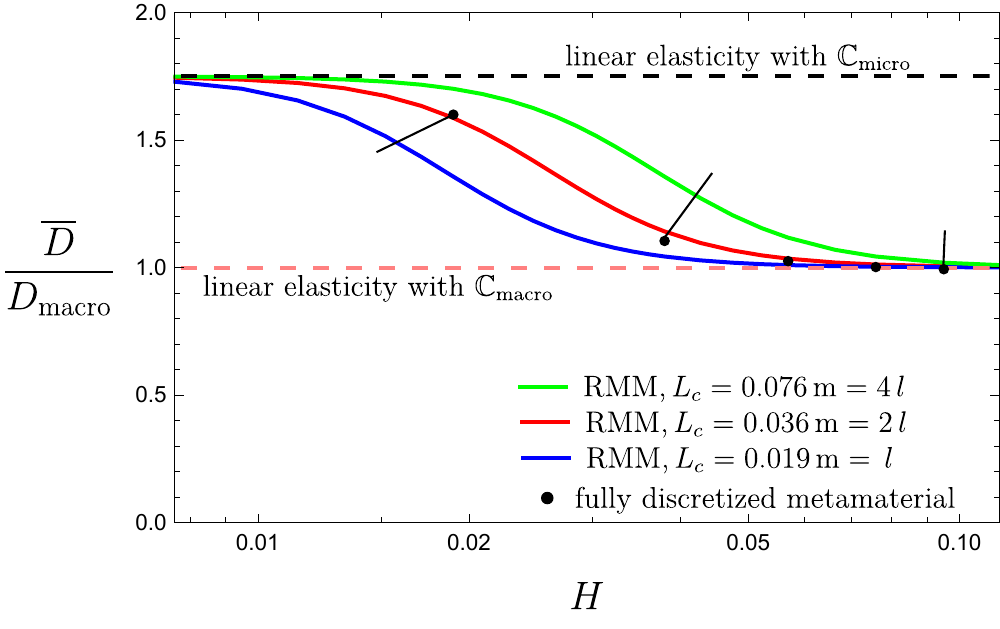}
     	\end{picture}
   		  	\caption{  $\Cmicro= 1.75  \, \Cmacro $}
    \end{subfigure}
      		 \begin{subfigure}[b]{0.49 \textwidth}
		 \center
		 	\begin{picture}(100,52)
		 \put(22,35){\color[rgb]{0,0,0}\makebox(0,0)[lb]{\small\smash{\tiny $n=1$}}}%	 
		 \put(52,36){\color[rgb]{0,0,0}\makebox(0,0)[lb]{\small\smash{\tiny $n=2$}}}%	  	 	
		 \put(70,33){\color[rgb]{0,0,0}\makebox(0,0)[lb]{\small\smash{\tiny$n=5$}}}%	
  \includegraphics[width=1 \textwidth]{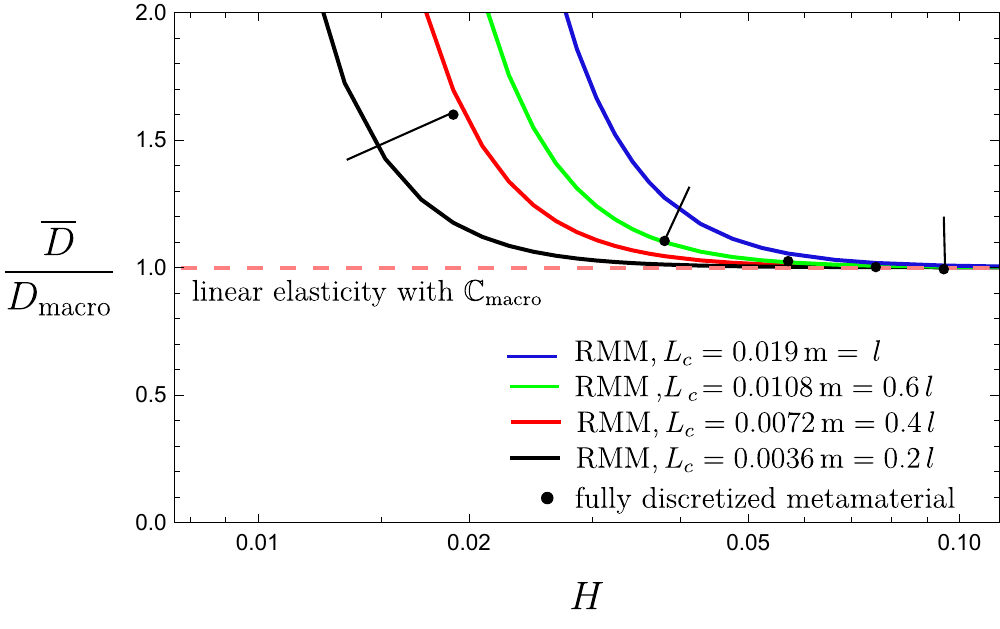}
     	\end{picture}
   		  	\caption{$\Cmicro= 1000  \, \Cmacro $; ``Cosserat type" }
    \end{subfigure}
		  	\caption{The normalized bending stiffness varying the beam size $H \times L = n \, l \times 12 \, n \, l$ obtained by the fully discretized metamaterial and the relaxed micromorphic model. We analyze here different choices for $\Cmicro$ with varying $\Lc$ and fixing $\mu=\mu_\textrm{macro}$. {Assuming~$\Cmicro = 10000 \, \Cmacro$~yields~the~same~results~as~in~(d).}}
	\label{Figure:stiffness_rmm_case1-lc}
\end{figure}

 \begin{figure}[ht]
	\unitlength=1mm
	\center
 		 \begin{subfigure}[b]{0.49 \textwidth}
		 \center
		 	\begin{picture}(100,50)
		 \put(22,36){\color[rgb]{0,0,0}\makebox(0,0)[lb]{\small\smash{\tiny $n=1$}}}%	 
		 \put(52,36){\color[rgb]{0,0,0}\makebox(0,0)[lb]{\small\smash{\tiny $n=2$}}}%	  	 	
		 \put(70,34){\color[rgb]{0,0,0}\makebox(0,0)[lb]{\small\smash{\tiny$n=5$}}}%	
         \includegraphics[width=1 \textwidth]{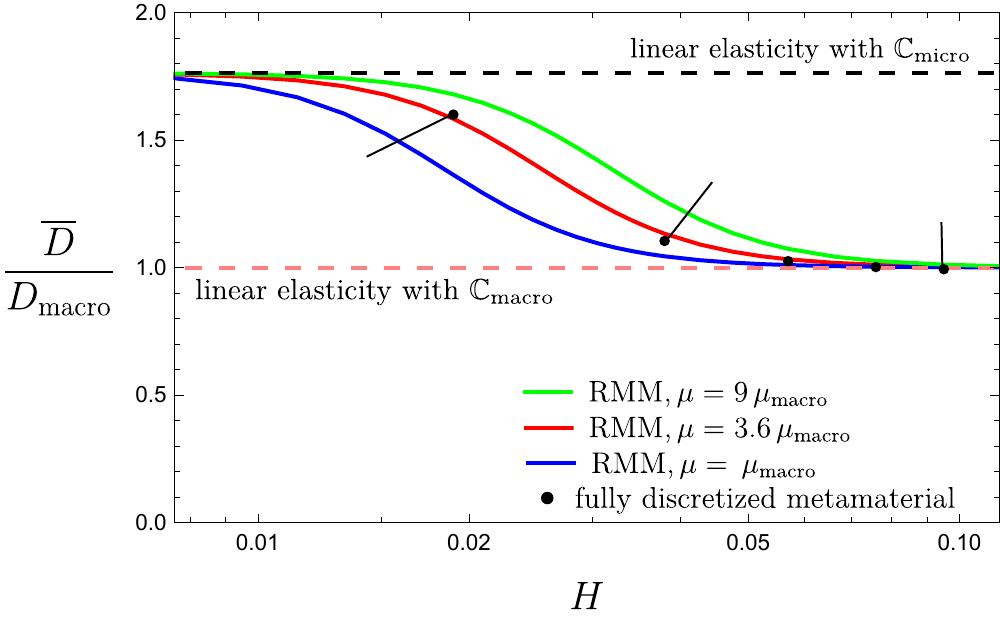}
	\end{picture}
   		  	\caption{$\Cmicro=  1.66 \, {\mathbb{C}}^\textrm{L\"owner}_\textrm{micro} $}
    \end{subfigure}
 		 \begin{subfigure}[b]{0.49 \textwidth}
		 \center
		 	\begin{picture}(100,50)
		 \put(22,36){\color[rgb]{0,0,0}\makebox(0,0)[lb]{\small\smash{\tiny $n=1$}}}%	 
		 \put(52,36){\color[rgb]{0,0,0}\makebox(0,0)[lb]{\small\smash{\tiny $n=2$}}}%	  	 	
		 \put(70,34){\color[rgb]{0,0,0}\makebox(0,0)[lb]{\small\smash{\tiny$n=5$}}}%	
         \includegraphics[width=1 \textwidth]{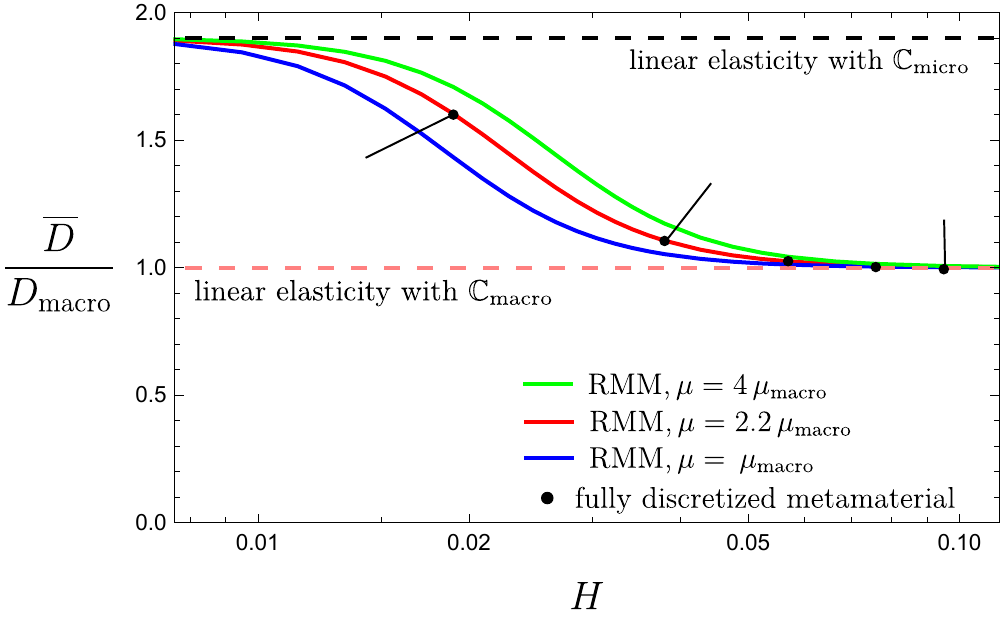}
   	\end{picture}
   		  	\caption{$\Cmicro=  \, \Cmatrix $}
    \end{subfigure}
      		 \begin{subfigure}[b]{0.49 \textwidth}
		 \center
		 	\begin{picture}(100,52)
		 \put(23,36){\color[rgb]{0,0,0}\makebox(0,0)[lb]{\small\smash{\tiny $n=1$}}}%	 
		 \put(53,37){\color[rgb]{0,0,0}\makebox(0,0)[lb]{\small\smash{\tiny $n=2$}}}%	  	 	
		 \put(70,33){\color[rgb]{0,0,0}\makebox(0,0)[lb]{\small\smash{\tiny$n=5$}}}%	
  \includegraphics[width=1 \textwidth]{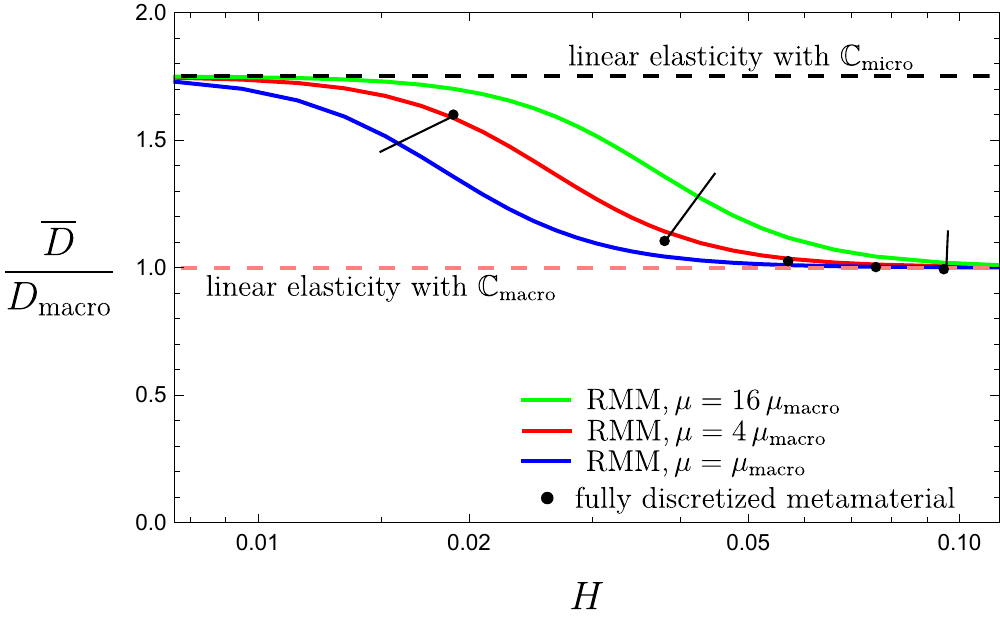}
     	\end{picture}
   		  	\caption{ $\Cmicro= 1.75  \, \Cmacro $}
    \end{subfigure}
      		 \begin{subfigure}[b]{0.49 \textwidth}
		 \center
		 	\begin{picture}(100,52)
		 \put(22,35){\color[rgb]{0,0,0}\makebox(0,0)[lb]{\small\smash{\tiny $n=1$}}}%	 
		 \put(52,35){\color[rgb]{0,0,0}\makebox(0,0)[lb]{\small\smash{\tiny $n=2$}}}%	  	 	
		 \put(70,33){\color[rgb]{0,0,0}\makebox(0,0)[lb]{\small\smash{\tiny$n=5$}}}%	
  \includegraphics[width=1 \textwidth]{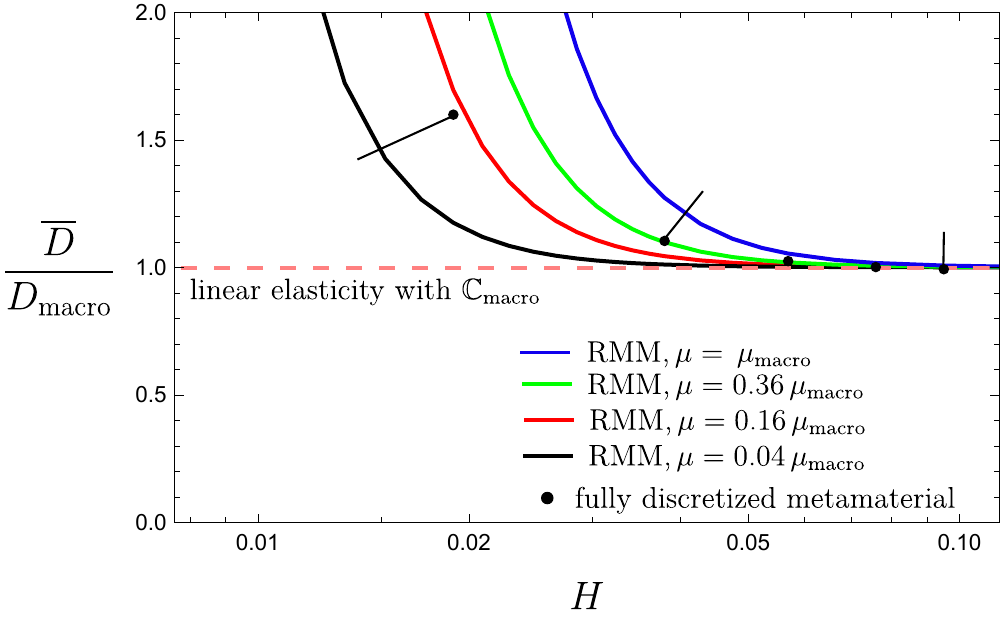}
     	\end{picture}
   		  	\caption{$\Cmicro= 1000  \, \Cmacro $; ``Cosserat type" }
    \end{subfigure}
		  	\caption{The normalized bending stiffness varying the beam size $H \times L = n \, l \times 12 \, n \, l$  obtained by the fully discretized metamaterial and the relaxed micromorphic model. We analyze here different choices for $\Cmicro$ with varying $\mu$ and fixing $\Lc=l$. The results are equivalent for both loading cases. The relaxed micromorphic model shows bounded stiffness given by $\Cmicro$ in contrast to the Cosserat model.}
	\label{Figure:stiffness_rmm_case1-mu}
\end{figure}
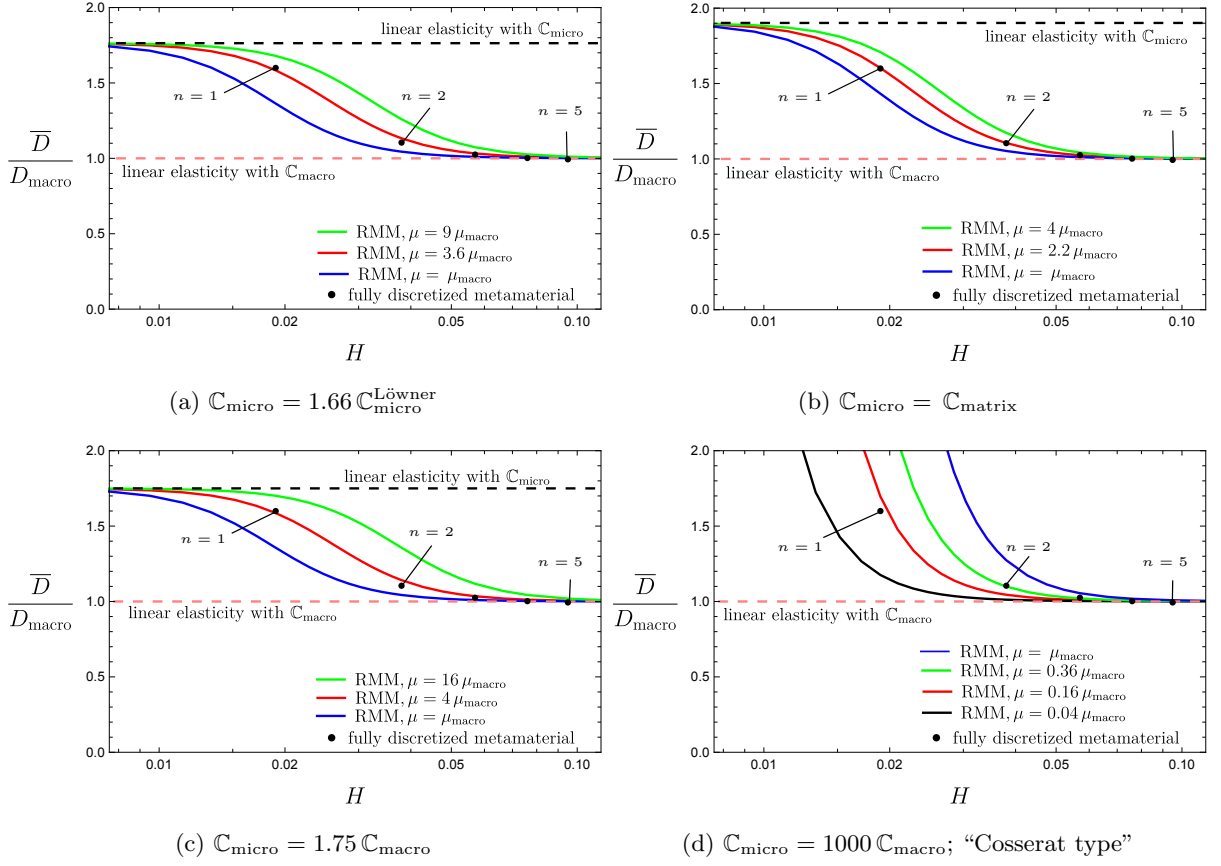 

\FloatBarrier

\section{Validation: further numerical examples}
\label{sec:validation}
This study assesses the obtained material parameters of the relaxed micromorphic model for two additional loading scenarios apart from the pure bending. The fully discretized metamaterial samples considered have the dimensions and material parameters as outlined before in Section \ref{sec:reference}. In the relaxed micromorphic model, we consider the symmetric case where $\mu_c=0$. The macro-scale elasticity tensor, $\Cmacro$, is defined in Section \ref{sec:cmacro} and the curvature is scaled to the specimen's size using Eq. \ref{eq:cal} under the assumption of $\mu = \mu_\textrm{macro}$. The micro-scale elasticity tensor will be determined using the same four different assumptions outlined in Section \ref{sec:fit}.

\subsection{Simple shearing} 

The boundary conditions are derived from the solution of an infinite stripe under simple shear in \cite{RizHueMadNef:2021:aso3}

 \begin{equation}
 \bu =  \overline \bu  = \left( \begin{array}{c}
a y \\
0  
 \end{array} \right)  \qquad \textrm{on} \qquad \partial \B \,,
 \end{equation}
 
which leads to the following strain and stress tensors for the homogeneous macro-elasticity case
 
 \begin{equation}
\overline \Bvarepsilon = \left(
 \begin{array}{c c}
0 & a /2 \\ 
a/2 & 0 
 \end{array}
 \right) , \qquad  \overline \Bsigma = \left(
 \begin{array}{c c}
0 & a \, \mu^*_\textrm{macro} \\ 
a \, \mu^*_\textrm{macro} & 0 
 \end{array}
 \right) \,.
 \end{equation}
The boundary value problems for the relaxed micromorphic model and the reference full detailed metamaterial are depicted in Fig. \ref{Figure:Geo3}. 
Dirichlet boundary condition for the displacement field and the consistent coupling condition must be satisfied over the entire boundary. 

  \begin{figure}[ht] 
\center
	\unitlength=1mm
	\begin{picture}(150,45)
	\put(0,0){\def\svgwidth{15cm}{\small\input{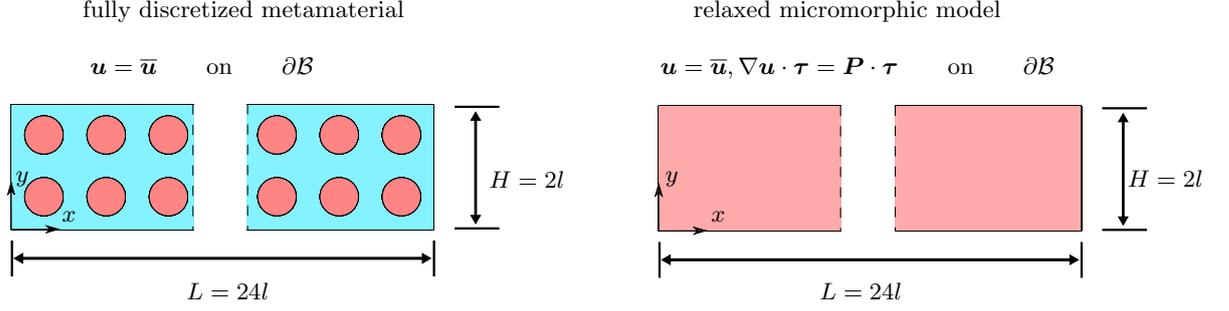}}}
	\end{picture}
	\caption{The geometry of the boundary value problem (shear) shown for $n=2$ for the fully discretized metamaterial and the homogeneous relaxed micrmorphic continuum.}
	\label{Figure:Geo3}
\end{figure} 

The size-effect is analyzed through the relative shear force $\frac{T}{T_\textrm{macro}}$, which is shown in Fig. \ref{Fig:shear_results}. The macro-scale shear force is given by $T_\textrm{macro} = a \, \mu^*_\textrm{macro} \, L$. The shear response of the assumed metamaterial is less influenced by its size compared to its response to bending. We notice that the choices $\Cmicro=  1.66 \, {\mathbb{C}}^\textrm{L\"owner}_\textrm{micro} $ and $\Cmicro= 1.75  \, \Cmacro $ deliver close results for the bending in Fig. \ref{Figure:stiffness_rmm_case1-lc} but different results for the simple shear in Fig \ref{Fig:shear_results} which can be explained by their different anisotropy properties. 
 
 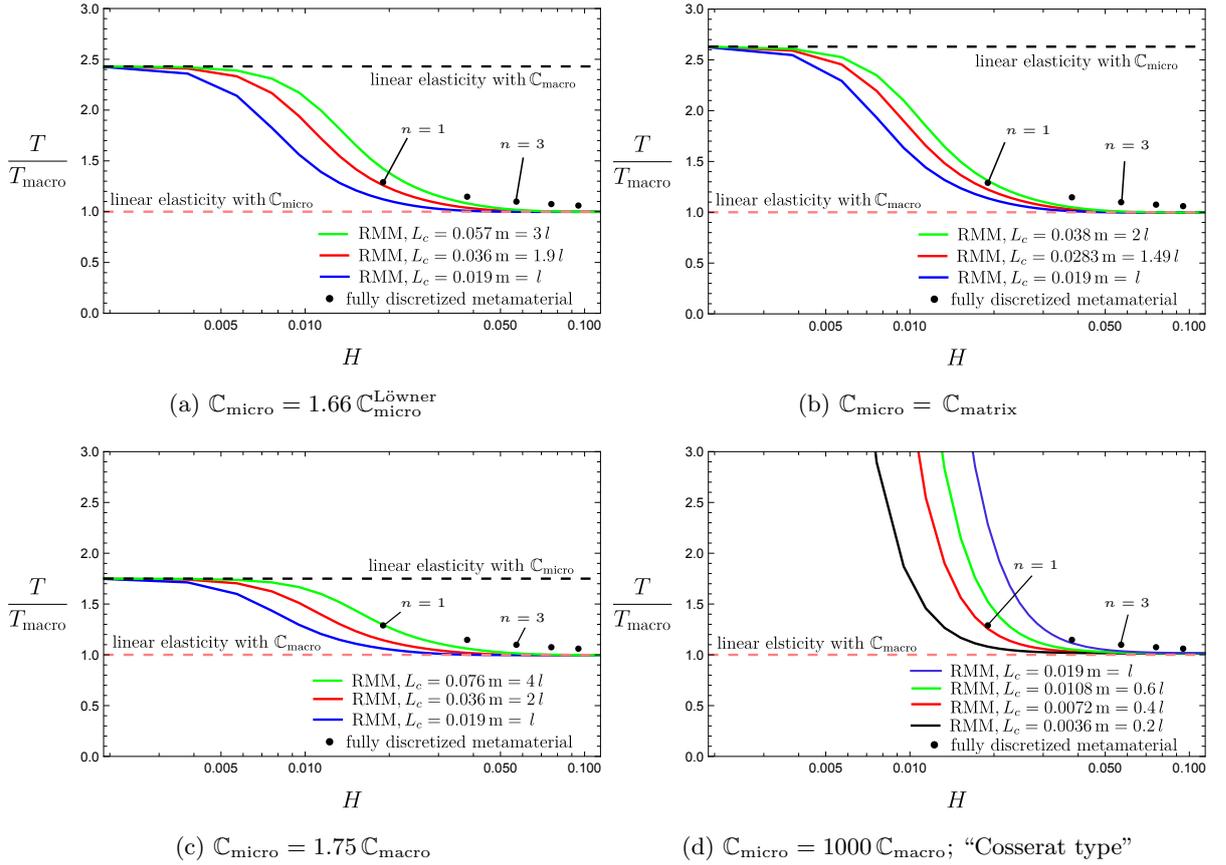
\begin{figure}[ht]
	\unitlength=1mm
	\center
 		 \begin{subfigure}[b]{0.49 \textwidth}
		 \center
		 	\begin{picture}(100,50)
		 \put(52,32){\color[rgb]{0,0,0}\makebox(0,0)[lb]{\small\smash{\tiny $n=1$}}}%	  	 	
		 \put(65,30){\color[rgb]{0,0,0}\makebox(0,0)[lb]{\small\smash{\tiny$n=3$}}}%	
         \includegraphics[width=1 \textwidth]{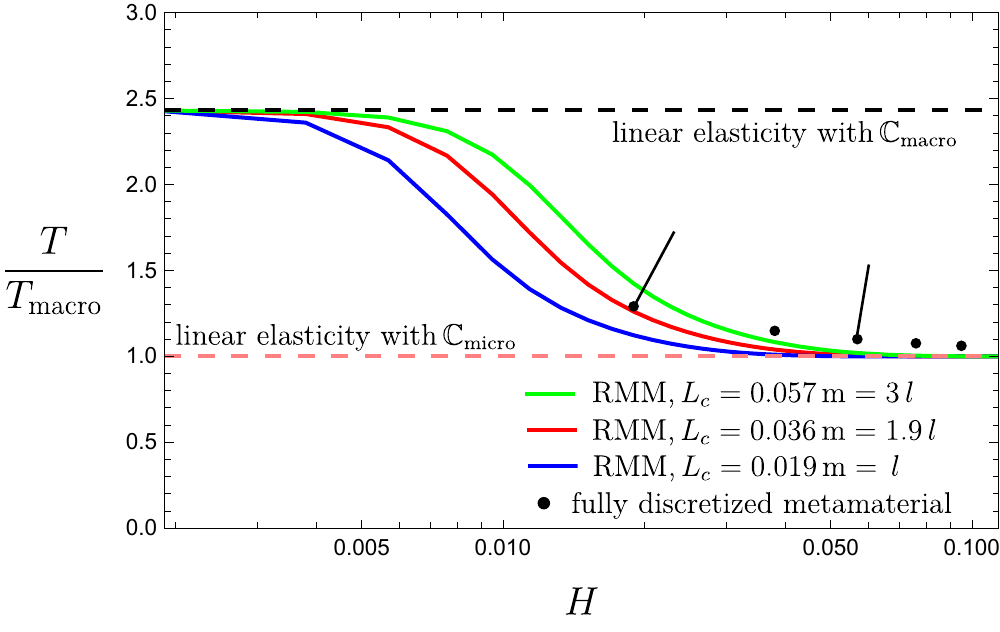}
	\end{picture}
   		  	\caption{$\Cmicro=  1.66 \, {\mathbb{C}}^\textrm{L\"owner}_\textrm{micro} $}
    \end{subfigure}
 		 \begin{subfigure}[b]{0.49 \textwidth}
		 \center
		 	\begin{picture}(100,50)
		 \put(52,32){\color[rgb]{0,0,0}\makebox(0,0)[lb]{\small\smash{\tiny $n=1$}}}%	  	 	
		 \put(65,30){\color[rgb]{0,0,0}\makebox(0,0)[lb]{\small\smash{\tiny$n=3$}}}%	
         \includegraphics[width=1 \textwidth]{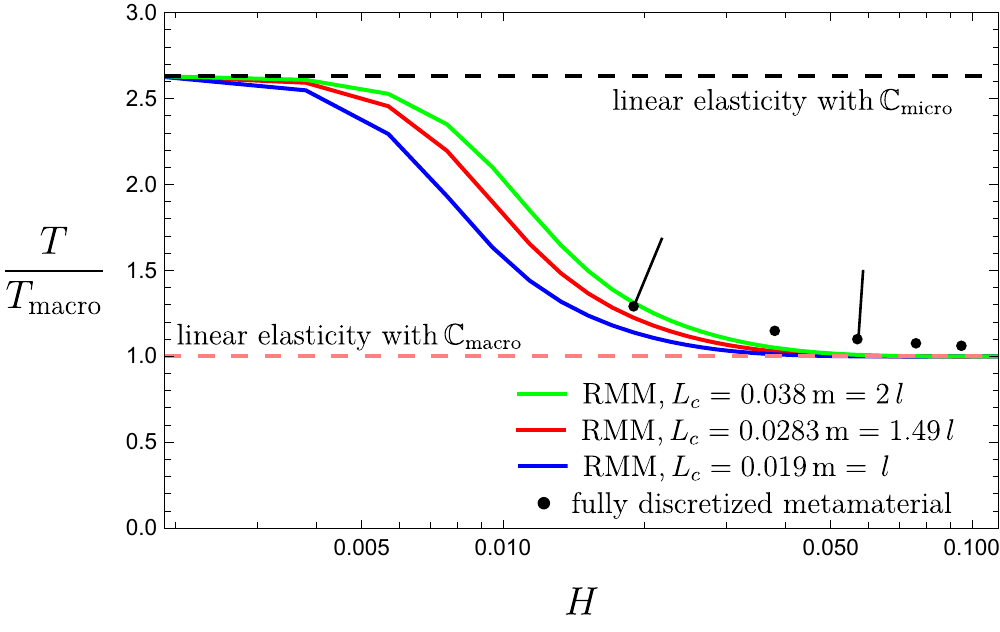}
   	\end{picture}
   		  	\caption{$\Cmicro=  \, \Cmatrix $}
    \end{subfigure}
      		 \begin{subfigure}[b]{0.49 \textwidth}
		 \center
		 	\begin{picture}(100,52)
		 \put(52,28){\color[rgb]{0,0,0}\makebox(0,0)[lb]{\small\smash{\tiny $n=1$}}}%	  	 	
		 \put(65,26){\color[rgb]{0,0,0}\makebox(0,0)[lb]{\small\smash{\tiny$n=3$}}}%	
  \includegraphics[width=1 \textwidth]{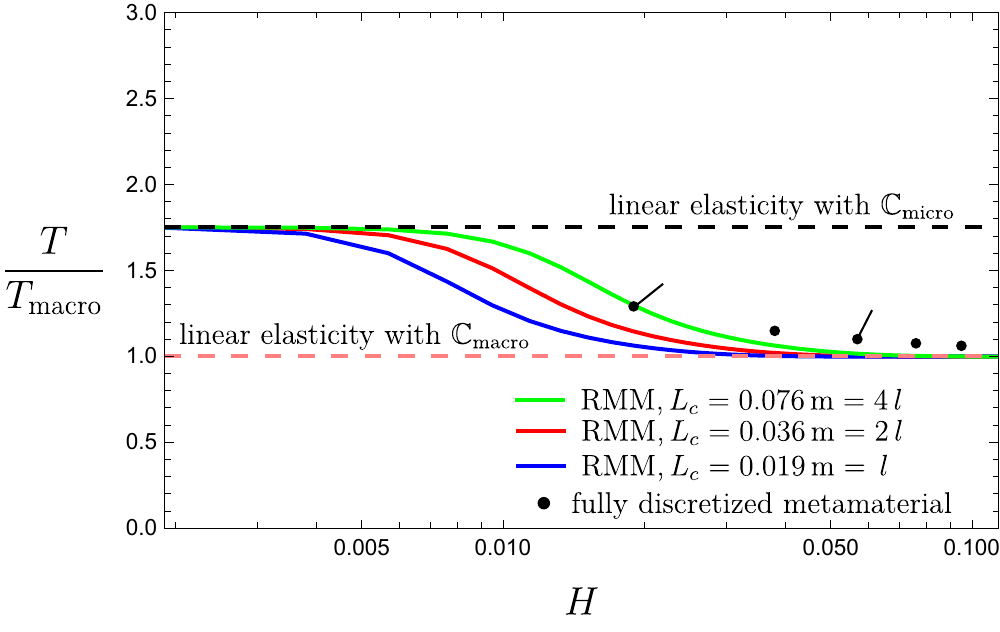}
     	\end{picture}
   		  	\caption{ $\Cmicro= 1.75  \, \Cmacro $}
    \end{subfigure}
      		 \begin{subfigure}[b]{0.49 \textwidth}
		 \center
		 	\begin{picture}(100,52)
		 \put(53,33){\color[rgb]{0,0,0}\makebox(0,0)[lb]{\small\smash{\tiny $n=1$}}}%	  	 	
		 \put(65,28){\color[rgb]{0,0,0}\makebox(0,0)[lb]{\small\smash{\tiny$n=3$}}}%	
  \includegraphics[width=1 \textwidth]{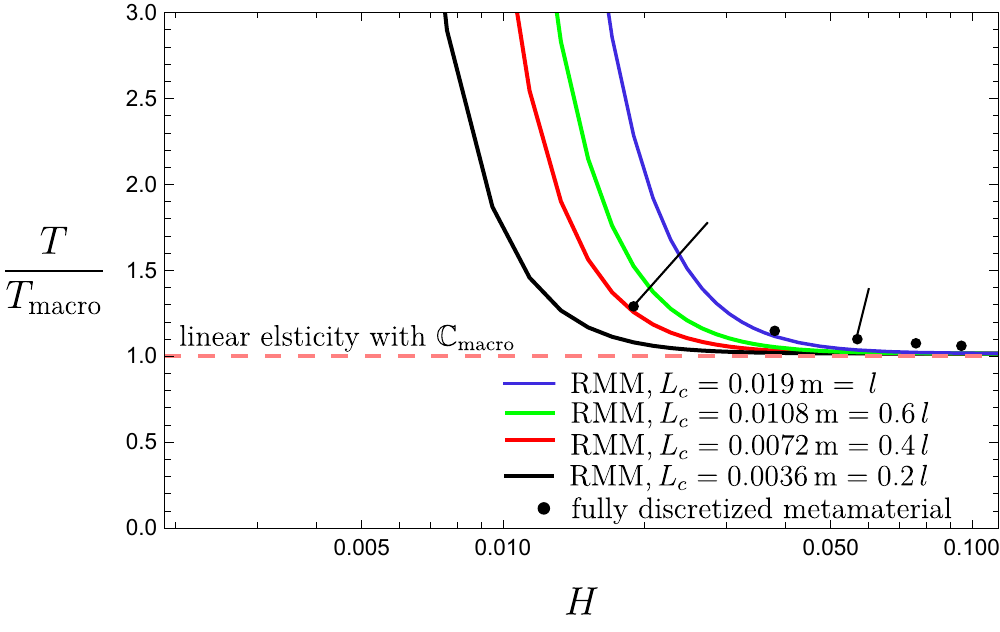}
     	\end{picture}
   		  	\caption{$\Cmicro= 1000  \, \Cmacro $; ``Cosserat type" }
    \end{subfigure}
		  	\caption{The relative shear force varying the specimen's size $H \times L$ for different choices of $\Cmicro$. }
\label{Fig:shear_results}
\end{figure} 

\FloatBarrier

\subsection{Cantilever under traction load}

In this setup, the right edge of the metamaterial is fixed in both directions while a constant traction of $t_y = \overline{t}$ is applied in the $y$-direction on the left side. The boundary value problems for both the fully discretized metamaterial and the relaxed micromorphic model are depicted in Fig. \ref{Figure:Geo4}. The micro elasticity can be recovered for large values of $\Lc$ when a consistent coupling condition is applied to the entire boundary. However, for small values of $\Lc$, a boundary layer is created at the upper and lower edges, requiring a fine mesh. This issue can be resolved by partially applying the consistent boundary condition, $(\nabla \bu \cdot \Btau)_y = (\Bdis \cdot \Btau)_y$.

  \begin{figure}[ht] 
\center
	\unitlength=1mm
	\begin{picture}(150,45)
	\put(0,0){\def\svgwidth{15cm}{\small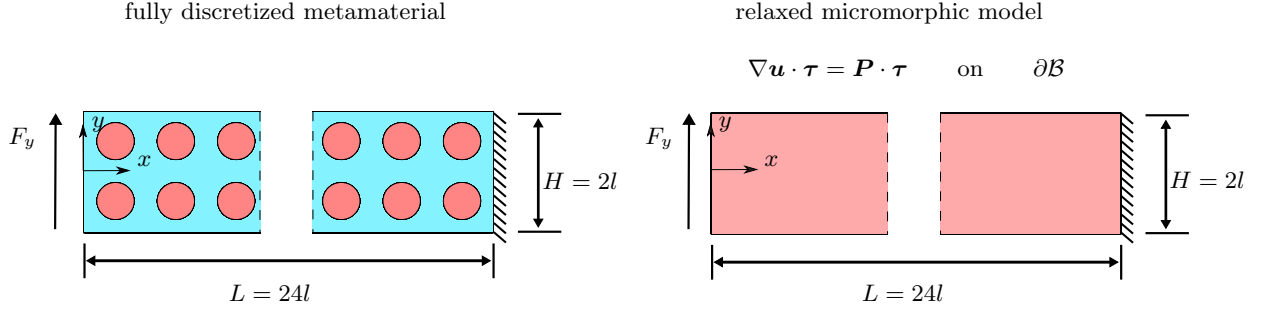}}
	\end{picture}

	\caption{ The geometry of the boundary value problem shown for $n=2$ for the fully discretized metamaterial and the homogeneous relaxed micrmorphic continuum.}
	\label{Figure:Geo4}
\end{figure} 

The equivalent beam model of the assumed cantilever, with the deformed shape illustrated for $n=2$, is displayed in Fig. \ref{Figure:canti_model}. The cantilever is subjected to a constant shear force $F_y = \overline{t} \, H$ and a linear moment that is zero on the left end and maximum on the right end $M = F_y \, x$. 

  \begin{figure}[ht] 
\center
	\unitlength=1mm
	\begin{picture}(150,40)
	\put(0,0){\def\svgwidth{15cm}{\small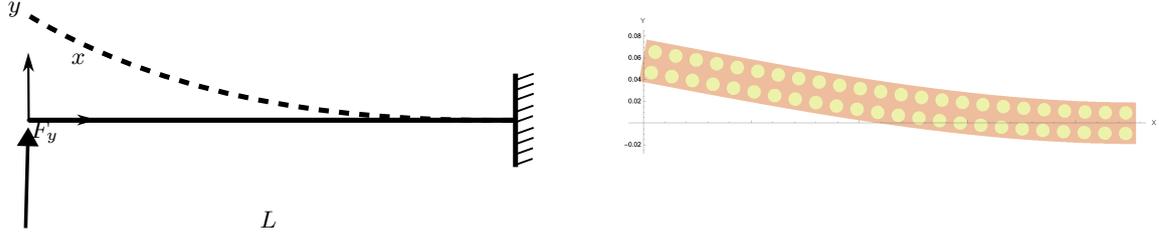}}
	\end{picture}
	\caption{The beam model of the cantilever and the deformed shape for $H \times L = 2 \, l \times 24 \, l$.}
	\label{Figure:canti_model}
\end{figure} 

The size-effect is analyzed by determining the inverse of the relative maximum displacement, expressed as $\frac{w_\textrm{macro}(0)}{w(0)}$. This calculation is illustrated in Fig. \ref{Fig:canti_results}. The macro-scale displacement is calculated using the formula $w_\textrm{macro}(0) = \frac{4 (1-\nu_\textrm{macro}^2) F_y L^3}{E_\textrm{macro} H^3}$. The results of both the fully discretized metamaterial and the relaxed micromorphic model show good agreement, as the dominant size-effect is bending. However, if consistent boundary conditions are not applied across the entire boundary, agreement is not achieved.

 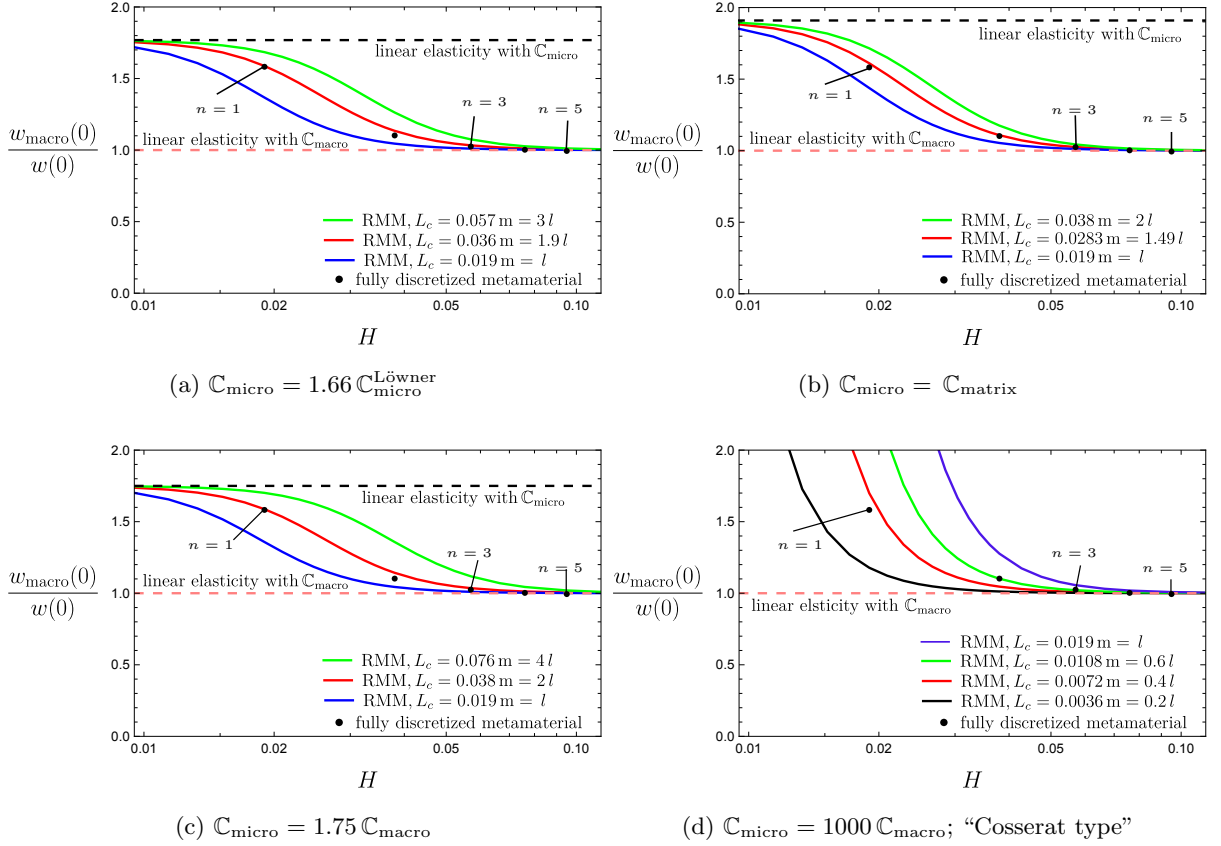
\begin{figure}[ht]
	\unitlength=1mm
	\center
 		 \begin{subfigure}[b]{0.49 \textwidth}
		 \center
		 	\begin{picture}(100,50)
		 \put(25,32){\color[rgb]{0,0,0}\makebox(0,0)[lb]{\small\smash{\tiny $n=1$}}}%	  	 	
		 \put(60,33){\color[rgb]{0,0,0}\makebox(0,0)[lb]{\small\smash{\tiny$n=3$}}}%	
		 \put(70,32){\color[rgb]{0,0,0}\makebox(0,0)[lb]{\small\smash{\tiny$n=5$}}}%	
         \includegraphics[width=1 \textwidth]{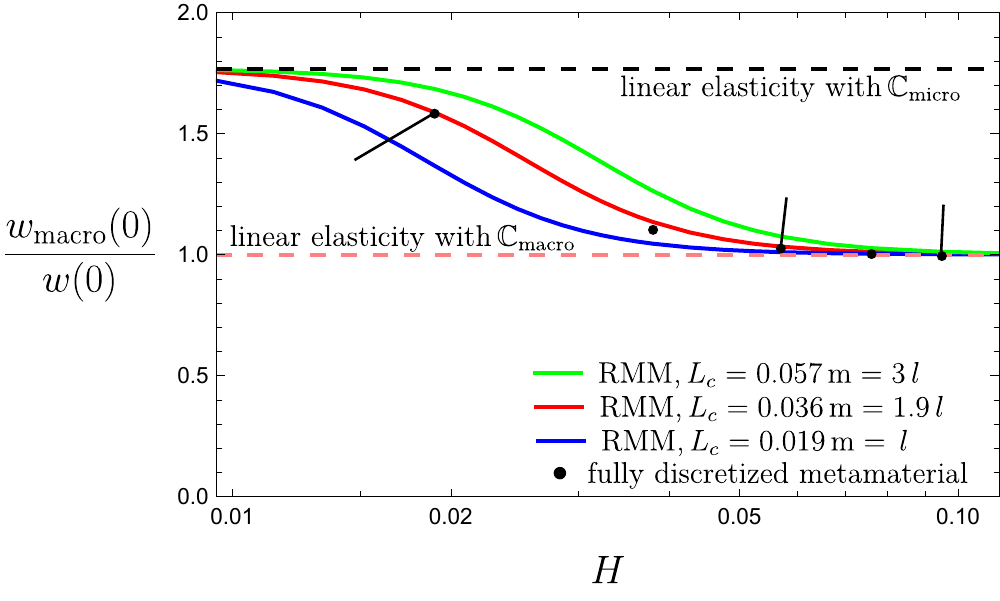}
	\end{picture}
   		  	\caption{$\Cmicro=  1.66 \, {\mathbb{C}}^\textrm{L\"owner}_\textrm{micro} $}
    \end{subfigure}
 		 \begin{subfigure}[b]{0.49 \textwidth}
		 \center
		 	\begin{picture}(100,50)
		 \put(26,34){\color[rgb]{0,0,0}\makebox(0,0)[lb]{\small\smash{\tiny $n=1$}}}%	  	 	
		 \put(58,32){\color[rgb]{0,0,0}\makebox(0,0)[lb]{\small\smash{\tiny$n=3$}}}%	
		 \put(70,31){\color[rgb]{0,0,0}\makebox(0,0)[lb]{\small\smash{\tiny$n=5$}}}%	
         \includegraphics[width=1 \textwidth]{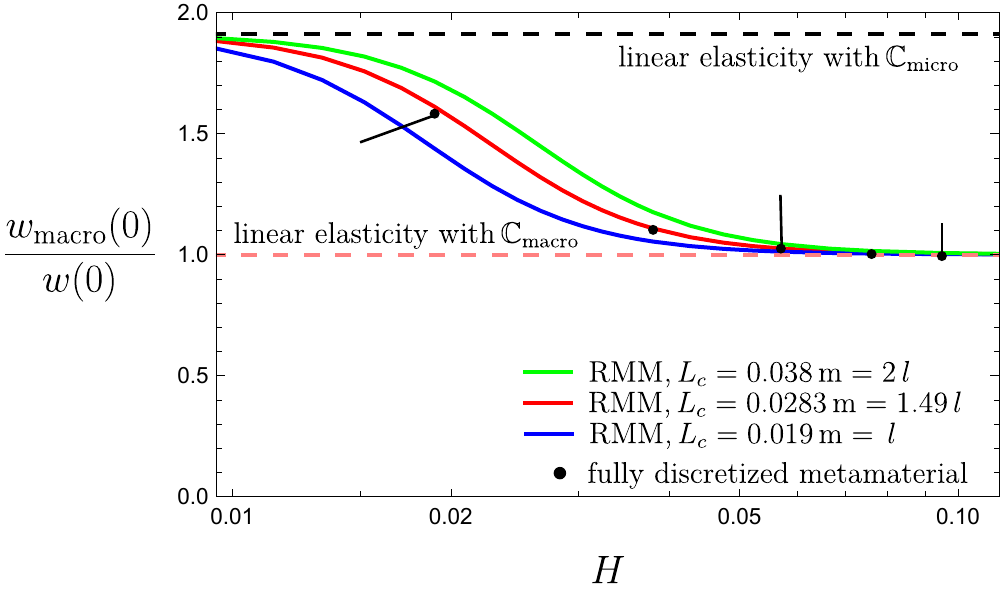}
   	\end{picture}
   		  	\caption{$\Cmicro=  \, \Cmatrix $}
    \end{subfigure}
      		 \begin{subfigure}[b]{0.49 \textwidth}
		 \center
		 	\begin{picture}(100,52)
		 \put(24,33){\color[rgb]{0,0,0}\makebox(0,0)[lb]{\small\smash{\tiny $n=1$}}}%	  	 	
		 \put(58,32){\color[rgb]{0,0,0}\makebox(0,0)[lb]{\small\smash{\tiny$n=3$}}}%	
		 \put(70,30){\color[rgb]{0,0,0}\makebox(0,0)[lb]{\small\smash{\tiny$n=5$}}}%	
  \includegraphics[width=1 \textwidth]{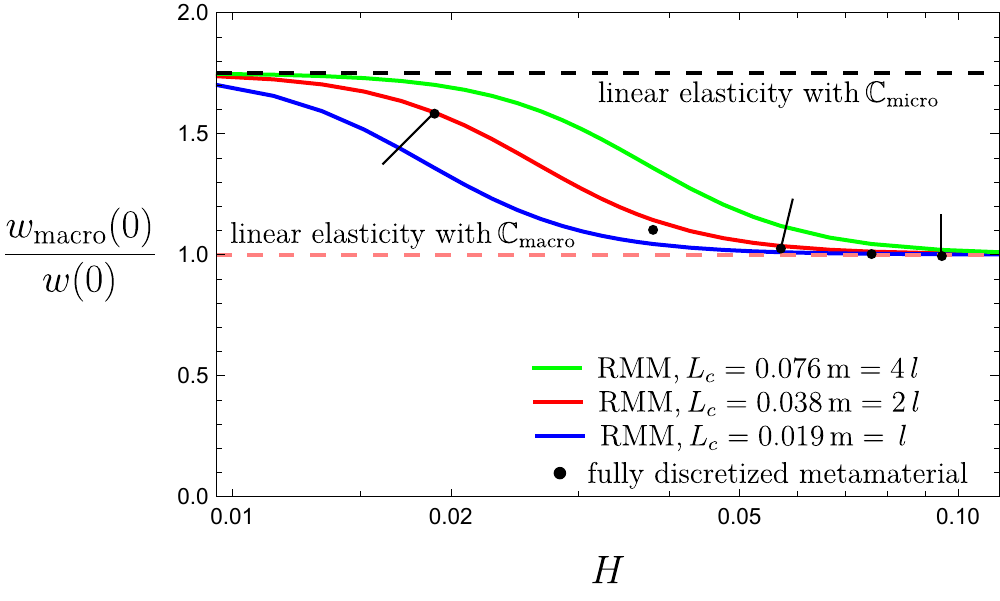}
     	\end{picture}
   		  	\caption{ $\Cmicro= 1.75  \, \Cmacro $}
    \end{subfigure}
      		 \begin{subfigure}[b]{0.49 \textwidth}
		 \center
		 	\begin{picture}(100,52)
		 \put(22,33){\color[rgb]{0,0,0}\makebox(0,0)[lb]{\small\smash{\tiny $n=1$}}}%	  	 	
		 \put(58,32){\color[rgb]{0,0,0}\makebox(0,0)[lb]{\small\smash{\tiny$n=3$}}}%	
		 \put(70,30){\color[rgb]{0,0,0}\makebox(0,0)[lb]{\small\smash{\tiny$n=5$}}}%
  \includegraphics[width=1 \textwidth]{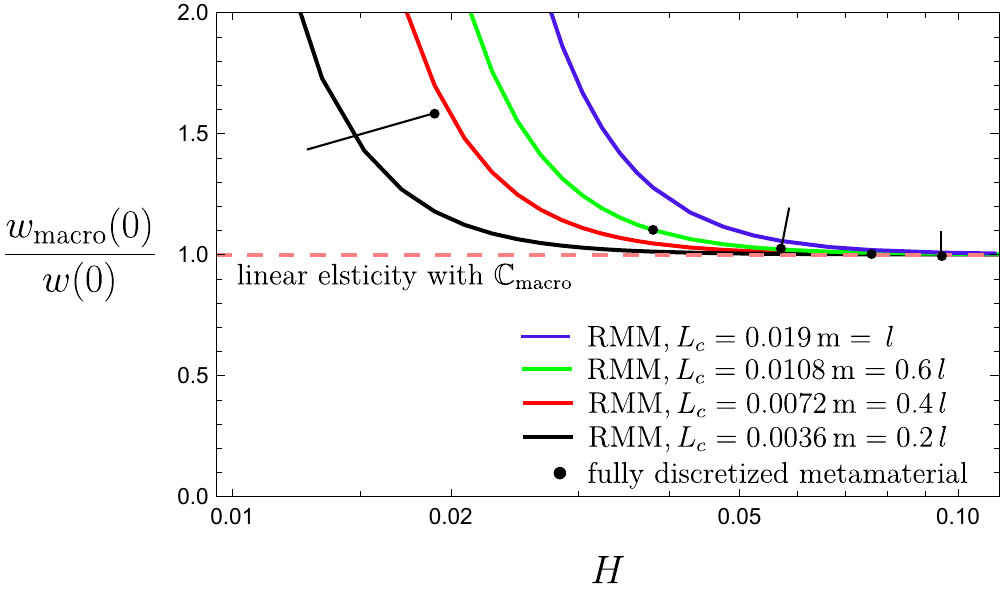}
     	\end{picture}
   		  	\caption{$\Cmicro= 1000  \, \Cmacro $; ``Cosserat type" }
    \end{subfigure}
		  	\caption{The inverse of the relative deflection at the free end of the cantilever ($x=0$) for varying the specimen's size $H \times L$ for different choices of $\Cmicro$. }
\label{Fig:canti_results}
\end{figure}

\FloatBarrier
 
 \section{Conclusions} 
  \label{sec:con}
 We introduced the relaxed micromorphic model with a brief description of the suitable tangential-conforming finite element formulation. We studied the size-effect phenomena of fully resolved beams under bending.  We have shown that applying a rotation (via a given displacement) or moment (applied traction) on the fully discretized metamaterial leads to similar results which we should get as well when we use the relaxed micromorphic model. We defined the macro elasticity tensor $\Cmacro$ by means of the standard periodic homogenization corresponding to large specimens. The micro elasticity tensor is connected to the stiffest possible response of the assumed metamaterial.  We introduced an approach to defining $\Cmicro$ which is based on the least upper bound of the apparent stiffness of the microstructure measured in the energy norm following the L\"owner matrix supremum problem where different variants of unit cells are considered under the affine Dirichlet boundary conditions. However, the flexural deformation mode is not captured by affine Dirichlet boundary conditions and the resulting elasticity tensor is much softer than the bent fully resolved metamaterial beams. Therefore, we scaled up the resulting elasticity tensor keeping its anisotropic cubic symmetry. Another procedure is tested to identify the micro elasticity tensor by non-affine boundary conditions (bending) on the unit cell or cluster of unit cells with the possible largest flexural rigidity. The boundary conditions were investigated for both loading cases (rotation or moment) for the symmetric and non-symmetric force stress. The consistent coupling boundary condition permits the model to work on the whole intended range bounded by linear elasticity with micro and macro elasticity tensors from above and below, respectively. We scaled the curvature measurement, which is isotropic in 2D, to account for the beam's size where a final fitting is conducted to decide the values of characteristic length and the shear modulus associated with the $\Curl$ of the micro-distortion field. The relaxed micromorphic delivers successfully the size-effects in a consistent manner for both loading cases. Finally, the relaxed micromorphic model was tested for two loading scenarios apart from pure bending with the consistent boundary condition applied on the entire boundary, highlighting its importance. Good agreement was obtained, however, the unique identification of the micro-elasticity tensor remains an open topic for future improvement. We established that the micro-elasticity tensor $\Cmicro$ must be stiffer than the apparent stiffness under the affine Dirichlet boundary conditions, but not stiffer than the homogeneous matrix.

{\bf Acknowledgment} \\
Funded by the Deutsche Forschungsgemeinschaft (DFG, German research Foundation) -  Project number 440935806 (SCHR 570/39-1, SCHE 2134/1-1, NE 902/10-1) within the DFG priority program 2256. 

%\input{relaxed_bending_sec}
%\input{examples_sec}
%\input{conclusions_sec}
% === list of references
% ------- layout-datei --------------
%\bibliographystyle{abbrv}

\bibliographystyle{plainnat}
%\bibliographystyle{plaindin}
%\bibliographystyle{plaindin_shortname2}
%\bibliographystyle{elsarticle-num-names}
% ------- bib-datei --------------
\FloatBarrier
{\footnotesize
\bibliography{micmag_01}

\begin{thebibliography}{126}
\providecommand{\natexlab}[1]{#1}
\providecommand{\url}[1]{\texttt{#1}}
\expandafter\ifx\csname urlstyle\endcsname\relax
  \providecommand{\doi}[1]{doi: #1}\else
  \providecommand{\doi}{doi: \begingroup \urlstyle{rm}\Url}\fi

\bibitem[Abali(2019)]{Abi:2019:rtp}
B.~E. Abali.
\newblock Revealing the physical insight of a length-scale parameter in
  metamaterials by exploiting the variational formulation.
\newblock \emph{Continuum Mechanics and Thermodynamics}, 31:\penalty0 885--894,
  2019.

\bibitem[Abali and Barchiesi(2021)]{AbaBar:2021:ami}
B.~E. Abali and E.~Barchiesi.
\newblock Additive manufacturing introduced substructure and computational
  determination of metamaterials parameters by means of the asymptotic
  homogenization.
\newblock \emph{Continuum Mechanics and Thermodynamics}, 33:\penalty0
  993--1009, 2021.

\bibitem[Abali et~al.(2019)Abali, Yang, and Papadopoulos]{AbaYanPap:2019:aca}
B.~E. Abali, H.~Yang, and P.~Papadopoulos.
\newblock A computational approach for determination of parameters in
  generalized mechanics.
\newblock In Holm Altenbach, Wolfgang~H. M{\"u}ller, and Bilen~Emek Abali,
  editors, \emph{Higher Gradient Materials and Related Generalized Continua},
  pages 1--18. Springer International Publishing, Cham, 2019.

\bibitem[Abali et~al.(2022)Abali, Vazic, and Newell]{AbaVazNew:2022:iom}
B.~E. Abali, B.~Vazic, and P.~Newell.
\newblock Influence of microstructure on size effect for metamaterials applied
  in composite structures.
\newblock \emph{Mechanics Research Communications}, 122:\penalty0 103877, 2022.

\bibitem[Aifantis(2011)]{Aif:2011:otg}
E.~C. Aifantis.
\newblock On the gradient approach - {R}elation to {E}ringen's nonlocal theory.
\newblock \emph{International Journal of Engineering Science}, 49\penalty0
  (12):\penalty0 1367--1377, 2011.

\bibitem[Aivaliotis et~al.(2020)Aivaliotis, Tallarico, d`Agostino, Daouadji,
  Neff, and Madeo]{AivTalAgoDaoNefMad:2020:fan}
A.~Aivaliotis, D.~Tallarico, M.~V. d`Agostino, A.~Daouadji, P.~Neff, and
  A.~Madeo.
\newblock Frequency- and angle-dependent scattering of a finite-sized
  meta-structure via the relaxed micromorphic model.
\newblock \emph{Archive of Applied Mechanics}, 90:\penalty0 1073--1096, 2020.

\bibitem[Al-Basyouni et~al.(2015)Al-Basyouni, Tounsi, and
  Mahmoud]{AlAbdMah:2015:sdb}
K.~S. Al-Basyouni, A.~Tounsi, and S.~R. Mahmoud.
\newblock Size dependent bending and vibration analysis of functionally graded
  micro beams based on modified couple stress theory and neutral surface
  position.
\newblock \emph{Composite Structures}, 125:\penalty0 621--630, 2015.

\bibitem[Alavi et~al.(2021)Alavi, Ganghoffer, Reda, and
  Sadighi]{AlaGanRedSad:2021:com}
S.~E. Alavi, J.-F. Ganghoffer, H.~Reda, and M.~Sadighi.
\newblock Construction of micromorphic continua by homogenization based on
  variational principles.
\newblock \emph{Journal of the Mechanics and Physics of Solids}, 153:\penalty0
  104278, 2021.

\bibitem[Alavi et~al.(2022{\natexlab{a}})Alavi, Ganghoffer, and
  Sadighi]{AlaGanSad:2022:cch}
S.~E. Alavi, J.-F. Ganghoffer, and M.~Sadighi.
\newblock Chiral {C}osserat homogenized constitutive models of architected
  media based on micromorphic homogenization.
\newblock \emph{Mathematics and Mechanics of Solids}, 27\penalty0
  (10):\penalty0 2287--2313, 2022{\natexlab{a}}.

\bibitem[Alavi et~al.(2022{\natexlab{b}})Alavi, Ganghoffer, Sadighi,
  Nasimsobhan, and Akbarzadeh]{AlaGanSadNasAkb:2022:cmo}
S.E. Alavi, J.F. Ganghoffer, M.~Sadighi, M.~Nasimsobhan, and A.H. Akbarzadeh.
\newblock Continualization method of lattice materials and analysis of size
  effects revisited based on cosserat models.
\newblock \emph{International Journal of Solids and Structures},
  254-255:\penalty0 111894, 2022{\natexlab{b}}.

\bibitem[Altan and Aifantis(1997)]{AltAif:1997:osa}
B.~S. Altan and E.~C. Aifantis.
\newblock On some aspects in the special theory of gradient elasticity.
\newblock \emph{Journal of the Mechanical Behavior of Materials}, 8\penalty0
  (3):\penalty0 231--282, 1997.

\bibitem[Askes and Aifantis(2011)]{AskAif:2011:gei}
H.~Askes and E.~C. Aifantis.
\newblock Gradient elasticity in statics and dynamics: An overview of
  formulations, length scale identification procedures, finite element
  implementations and new results.
\newblock \emph{International Journal of Solids and Structures}, 48\penalty0
  (13):\penalty0 1962--1990, 2011.

\bibitem[Askes et~al.(2008)Askes, Metrikine, Pichugin, and
  Bennett]{AskMetPicBen:2008:fsg}
H.~Askes, A.~V. Metrikine, A.~V. Pichugin, and T.~Bennett.
\newblock Four simplified gradient elasticity models for the simulation of
  dispersive wave propagation.
\newblock \emph{Philosophical Magazine}, 88\penalty0 (28-29):\penalty0
  3415--3443, 2008.

\bibitem[Auffray et~al.(2010)Auffray, Bouchet, and
  Br{\'e}chet]{AufBouBre:2010:sge}
N.~Auffray, R.~Bouchet, and Y.~Br{\'e}chet.
\newblock Strain gradient elastic homogenization of bidimensional cellular
  media.
\newblock \emph{International Journal of Solids and Structures}, 47\penalty0
  (13):\penalty0 1698--1710, 2010.

\bibitem[Bacigalupo and Gambarotta(2010)]{BacGam:2010:soc}
A.~Bacigalupo and L.~Gambarotta.
\newblock Second-order computational homogenization of heterogeneous materials
  with periodic microstructure.
\newblock \emph{Journal of Applied Mathematics and Mechanics}, 90\penalty0
  (10-11):\penalty0 796--811, 2010.

\bibitem[Bacigalupo et~al.(2018)Bacigalupo, Paggi, {Dal Corso}, and
  Bigoni]{BacPagDalBig:2018:iof}
A.~Bacigalupo, M.~Paggi, F.~{Dal Corso}, and D.~Bigoni.
\newblock Identification of higher-order continua equivalent to a {C}auchy
  elastic composite.
\newblock \emph{Mechanics Research Communications}, 93:\penalty0 11--22, 2018.
\newblock Mechanics from the 20th to the 21st Century: The Legacy of G\'erard
  A. Maugin.

\bibitem[Barbagallo et~al.(2017)Barbagallo, Madeo, d'Agostino, Abreu, Ghiba,
  and Neff]{BabMadDagAbrGhiNeff:2017:taf}
G.~Barbagallo, A.~Madeo, M.~V. d'Agostino, R.~Abreu, I.-D. Ghiba, and P.~Neff.
\newblock Transparent anisotropy for the relaxed micromorphic model:
  Macroscopic consistency conditions and long wave length asymptotics.
\newblock \emph{International Journal of Solids and Structures}, 120:\penalty0
  7--30, 2017.

\bibitem[Barbagallo et~al.(2019)Barbagallo, Tallarico, d'Agostino, Aivaliotis,
  Neff, and Madeo]{BarTalDagAivNefMad:2019:rmm}
G.~Barbagallo, D.~Tallarico, M.~V. d'Agostino, A.~Aivaliotis, P.~Neff, and
  A.~Madeo.
\newblock Relaxed micromorphic model of transient wave propagation in
  anisotropic band-gap metastructures.
\newblock \emph{International Journal of Solids and Structures}, 162:\penalty0
  148--163, 2019.

\bibitem[Berkache et~al.(2017)Berkache, Deogekar, Goda, Picu, and
  Ganghoffer]{BerDeoGodPicGan:2017:cof}
K.~Berkache, S.~Deogekar, I.~Goda, R.C. Picu, and J.-F. Ganghoffer.
\newblock Construction of second gradient continuum models for random fibrous
  networks and analysis of size effects.
\newblock \emph{Composite Structures}, 181:\penalty0 347--357, 2017.

\bibitem[Biswas and Poh(2017)]{BisPoh:2017:amc}
R.~Biswas and L.~H. Poh.
\newblock A micromorphic computational homogenization framework for
  heterogeneous materials.
\newblock \emph{Journal of the Mechanics and Physics of Solids}, 102:\penalty0
  187--208, 2017.

\bibitem[Blesgen and Neff(2022)]{BleNef:2022:ssi}
T.~Blesgen and P.~Neff.
\newblock Simple shear in nonlinear cosserat micropolar elasticity: Existence
  of minimizers, numerical simulations, and occurrence of microstructure.
\newblock \emph{Mathematics and Mechanics of Solids}, page 10812865221122191,
  2022.

\bibitem[Boutin(1996)]{Bou:1996:mei}
C.~Boutin.
\newblock Microstructural effects in elastic composites.
\newblock \emph{International Journal of Solids and Structures},
  33(7):\penalty0 1023--1053, 1996.

\bibitem[Carcaterra et~al.(2015)Carcaterra, dell`Isola, Esposito, and
  Pulvirenti]{CarDelEspPul:2015:mdo}
A.~Carcaterra, F.~dell`Isola, R.~Esposito, and M.~Pulvirenti.
\newblock Macroscopic description of microscopically strongly inhomogenous
  systems: A mathematical basis for the synthesis of higher gradients
  metamaterials.
\newblock \emph{Archive for Rational Mechanics and Analysis}, 218:\penalty0
  1239--1262, 2015.

\bibitem[Cosserat and Cosserat()]{CosCos:1909:tof}
E.~Cosserat and F.~Cosserat.
\newblock Th\`eorie des corps d\`eformable.
\newblock Librairie Scientifique A. Hermann et Fils, engl. translation by {D.
  H. Delphenich}, pdf available at
  (http://www.mathematik.tu-darmstadt.de/fbereiche/analysis/pde/staff/neff/patrizio/Cosserat.html),
  reprint 2009 by Hermann Librairie Scientifique, ISBN 978 27056 6920 1, Paris,
  1909.

\bibitem[d`Agostino et~al.(2020)d`Agostino, Barbagallo, Ghiba, Eidel, Neff, and
  Madeo]{DagBarGhiEidNefMad:2020:edo}
M.~V. d`Agostino, G.~Barbagallo, I.-D. Ghiba, B.~Eidel, P.~Neff, and A.~Madeo.
\newblock Effective description of anisotropic wave dispersion in mechanical
  band-gap metamaterials via the relaxed micromorphic model.
\newblock \emph{Journal of Elasticity}, 139:\penalty0 299--329, 2020.

\bibitem[d'Agostino et~al.(2022)d'Agostino, Rizzi, Khan, Lewintan, Madeo, and
  Neff]{DagRizKhaLewMadNef:2021:tcc}
M.~V. d'Agostino, G.~Rizzi, H.~Khan, P.~Lewintan, A.~Madeo, and P.~Neff.
\newblock The consistent coupling boundary condition for the classical
  micromorphic model: existence, uniqueness and interpretation of the
  parameters.
\newblock \emph{Continuum Mechanics and Thermodynamics}, 2022.
\newblock \doi{10.1007/s00161-022-01126-3}.

\bibitem[{Del Vescovo} and Giorgio(2014)]{DelGio:2014:dpo}
D.~{Del Vescovo} and I.~Giorgio.
\newblock Dynamic problems for metamaterials: Review of existing models and
  ideas for further research.
\newblock \emph{International Journal of Engineering Science}, 80:\penalty0
  153--172, 2014.

\bibitem[Demore et~al.(2022)Demore, Rizzi, Collet, Neff, and
  Madeo]{DemRizColNefMad:2022:uem}
F.~Demore, G.~Rizzi, M.~Collet, P.~Neff, and A.~Madeo.
\newblock Unfolding engineering metamaterials design: Relaxed micromorphic
  modeling of large-scale acoustic meta-structures.
\newblock \emph{Journal of the Mechanics and Physics of Solids}, 168:\penalty0
  104995, 2022.

\bibitem[El~Dhaba(2020)]{Dhaba:2020:rmm}
A.~R. El~Dhaba.
\newblock Reduced micromorphic model in orthogonal curvilinear coordinates and
  its application to a metamaterial hemisphere.
\newblock \emph{Scientific Reports}, 10, 2020.

\bibitem[Eremeyev et~al.(2021)Eremeyev, Cazzani, and
  dell'Isola]{EreCazDel:2021:ond}
V.~A. Eremeyev, A.~Cazzani, and F.~dell'Isola.
\newblock On nonlinear dilatational strain gradient elasticity.
\newblock \emph{Continuum Mechanics and Thermodynamics}, 33:\penalty0
  1429--1463, 2021.

\bibitem[Eringen(1968)]{Eri:1968:mom}
A.~C. Eringen.
\newblock Mechanics of micromorphic continua.
\newblock In \emph{Mechanics of Generalized Continua}, pages 18--35. Springer,
  Berlin, Heidelberg., 1968.

\bibitem[Eringen and Suhubi(1964)]{EriSub:1964:nto}
A.~C. Eringen and E.~S. Suhubi.
\newblock Nonlinear theory of simple micro-elastic solids-{I}.
\newblock \emph{International Journal of Engineering Science}, 2\penalty0
  (2):\penalty0 189--203, 1964.

\bibitem[Fischer et~al.(2011)Fischer, Klassen, Mergheim, Steinmann, and
  M\"uller]{FisKlaMerSteMue:2011:iao}
P.~Fischer, M.~Klassen, J.~Mergheim, P.~Steinmann, and R.~M\"uller.
\newblock Isogeometric analysis of {2D} gradient elasticity.
\newblock \emph{Computational Mechanics}, 47:\penalty0 1432--0924, 2011.

\bibitem[Fischer et~al.(2020)Fischer, Hillen, and Eberl]{FisHilEbe:2020:mmo}
S.~C.~L. Fischer, L.~Hillen, and C.~Eberl.
\newblock Mechanical metamaterials on the way from laboratory scale to
  industrial applications: Challenges for characterization and scalability.
\newblock \emph{Materials}, 13\penalty0 (16), 2020.

\bibitem[Forest(2002)]{For:2002:hma}
S.~Forest.
\newblock Homogenization methods and mechanics of generalized continua - part
  2.
\newblock \emph{Theoretical and Applied Mechanics}, 28-29:\penalty0 113--144,
  2002.

\bibitem[Forest(2016)]{For:2016:nro}
S.~Forest.
\newblock Nonlinear regularization operators as derived from the micromorphic
  approach to gradient elasticity, viscoplasticity and damage.
\newblock \emph{Proceedings of the Royal Society}, 472\penalty0
  (2188):\penalty0 20150755, 2016.

\bibitem[Forest and Sab(1998)]{ForSab:1998:com}
S.~Forest and K.~Sab.
\newblock Cosserat overall modeling of heterogeneous materials.
\newblock \emph{Mechanics Research Communications}, 25\penalty0 (4):\penalty0
  449--454, 1998.

\bibitem[Forest and Trinh(2011)]{ForTri:2011:gca}
S.~Forest and D.K. Trinh.
\newblock {Generalized continua and non-homogeneous boundary conditions in
  homogenisation methods}.
\newblock \emph{Journal of Applied Mathematics and Mechanics}, 91\penalty0
  (2):\penalty0 90--109, 2011.

\bibitem[Ganghoffer and Reda(2021)]{GanRed:2021:ava}
J.-F. Ganghoffer and H.~Reda.
\newblock A variational approach of homogenization of heterogeneous materials
  towards second gradient continua.
\newblock \emph{Mechanics of Materials}, 158:\penalty0 103743, 2021.
\newblock ISSN 0167-6636.
\newblock \doi{https://doi.org/10.1016/j.mechmat.2021.103743}.
\newblock URL
  \url{https://www.sciencedirect.com/science/article/pii/S0167663621000028}.

\bibitem[Ghiba et~al.(2015)Ghiba, Neff, Madeo, Placidi, and
  Rosi]{GhiNefMadPlaRos:2015:trl}
I.-D. Ghiba, P.~Neff, A.~Madeo, L.~Placidi, and G.~Rosi.
\newblock The relaxed linear micromorphic continuum: existence, uniqueness and
  continuous dependence in dynamics.
\newblock \emph{Mathematics and Mechanics of Solids}, 20\penalty0
  (10):\penalty0 1171--1197, 2015.

\bibitem[Ghiba et~al.(2022)Ghiba, Rizzi, Madeo, and
  Neff]{GhiRizMadNef:2022:cme}
I.-D. Ghiba, G.~Rizzi, A.~Madeo, and P.~Neff.
\newblock {Cosserat micropolar elasticity: classical Eringen vs. dislocation
  form}, 2022.
\newblock URL \url{https://arxiv.org/abs/2206.02473. To appear in Journal of
  Mechanics of Materials and Structures}.

\bibitem[Glaesener et~al.(2021)Glaesener, Bastek, Gonon, Kannan, Telgen,
  Sp\"ottling, Steiner, and Kochmann]{GlaBasConetal:2021:ctm}
R.~N. Glaesener, J.-H. Bastek, F.~Gonon, V.~Kannan, B.~Telgen, B.~Sp\"ottling,
  S.~Steiner, and D.~M. Kochmann.
\newblock Viscoelastic truss metamaterials as time-dependent generalized
  continua.
\newblock \emph{Journal of the Mechanics and Physics of Solids}, 156:\penalty0
  104569, 2021.

\bibitem[Goda and Ganghoffer(2016)]{GodGan:2016:cof}
I.~Goda and J.-F. Ganghoffer.
\newblock Construction of first and second order grade anisotropic continuum
  media for {3D} porous and textile composite structures.
\newblock \emph{Composite Structures}, 141:\penalty0 292--327, 2016.

\bibitem[Golaszewski et~al.(2019)Golaszewski, Grygoruk, Giorgio, Laudato, and
  Cosmo]{GolGryGioLauCos:2019:mwr}
M.~Golaszewski, R.~Grygoruk, I.~Giorgio, M.~Laudato, and F.~D. Cosmo.
\newblock Metamaterials with relative displacements in their microstructure:
  technological challenges in 3d printing, experiments and numerical
  predictions.
\newblock \emph{Continuum Mechanics and Thermodynamics}, 31:\penalty0
  1015--1034, 2019.

\bibitem[Hosseini and Niiranen(2022)]{HosNii:2022:3sg}
S.~B. Hosseini and J.~Niiranen.
\newblock 3{D} strain gradient elasticity: Variational formulations,
  isogeometric analysis and model peculiarities.
\newblock \emph{Computer Methods in Applied Mechanics and Engineering},
  389:\penalty0 114324, 2022.

\bibitem[H\"utter(2017)]{Hue:2017:hoa}
G.~H\"utter.
\newblock Homogenization of a {C}auchy continuum towards a micromorphic
  continuum.
\newblock \emph{Journal of the Mechanics and Physics of Solids}, 99:\penalty0
  394--408, 2017.

\bibitem[H\"utter(2019)]{Hue:2019:otm}
G.~H\"utter.
\newblock On the micro-macro relation for the microdeformation in the
  homogenization towards micromorphic and micropolar continua.
\newblock \emph{Journal of the Mechanics and Physics of Solids}, 127:\penalty0
  62--79, 2019.

\bibitem[H\"utter(2022)]{Hue:2022:iom}
G.~H\"utter.
\newblock Interpretation of micromorphic constitutive relations for porous
  materials at the microscale via harmonic decomposition.
\newblock \emph{Journal of the Mechanics and Physics of Solids}, page 105135,
  2022.

\bibitem[Ju et~al.(2021)Ju, Mahnken, Liang, and Xu]{JuMahLiaXu:2021:CS}
X.~Ju, R.~Mahnken, L.~Liang, and Y.~Xu.
\newblock Goal-oriented mesh adaptivity for inverse problems in linear
  micromorphic elasticity.
\newblock \emph{Computers and Structures}, 257:\penalty0 106671, 2021.

\bibitem[Khakalo and Niiranen(2019)]{KhaNii:2019:lsa}
S.~Khakalo and J.~Niiranen.
\newblock Lattice structures as thermoelastic strain gradient metamaterials:
  Evidence from full-field simulations and applications to functionally
  step-wise-graded beams.
\newblock \emph{Composites Part B: Engineering}, 177:\penalty0 107224, 2019.

\bibitem[Khakalo and Niiranen(2020)]{KhaNii:2020:asg}
S.~Khakalo and J.~Niiranen.
\newblock Anisotropic strain gradient thermoelasticity for cellular structures:
  Plate models, homogenization and isogeometric analysis.
\newblock \emph{Journal of the Mechanics and Physics of Solids}, 134:\penalty0
  103728, 2020.

\bibitem[Khakalo et~al.(2018)Khakalo, Balobanov, and
  Niiranen]{KhaBalNii:2018:msd}
S.~Khakalo, V.~Balobanov, and J.~Niiranen.
\newblock Modelling size-dependent bending, buckling and vibrations of {2D}
  triangular lattices by strain gradient elasticity models: Applications to
  sandwich beams and auxetics.
\newblock \emph{International Journal of Engineering Science}, 127:\penalty0
  33--52, 2018.

\bibitem[Khan et~al.(2022)Khan, Ghiba, Madeo, and Neff]{KhaChiMadNef:2022:eau}
H.~Khan, I.-D. Ghiba, A.~Madeo, and P.~Neff.
\newblock {Existence and uniqueness of Rayleigh waves in isotropic elastic
  Cosserat materials and algorithmic aspects}.
\newblock \emph{Wave Motion}, 110:\penalty0 102898, 2022.

\bibitem[Kirby et~al.(2012)Kirby, Logg, Rognes, and
  Terrel]{KirLogRogTer:2012:cau}
R.~C. Kirby, A.~Logg, M.~E. Rognes, and A.~R. Terrel.
\newblock Common and unusual finite elements.
\newblock In \emph{Automated Solution of Differential Equations by the Finite
  Element Method: The FEniCS Book}, pages 95--119. Springer Berlin Heidelberg,
  Berlin, Heidelberg, 2012.

\bibitem[Knees et~al.(2022)Knees, Owczarek, and Neff]{KneOwcNef:2022:alr}
D.~Knees, S.~Owczarek, and P.~Neff.
\newblock A local regularity result for the relaxed micromorphic model based on
  inner variations.
\newblock \emph{to appear in Journal of Mathematical Analysis and
  Applications}, 2022.
\newblock \doi{10.48550/ARXIV.2208.04821}.
\newblock URL \url{https://arxiv.org/abs/2208.04821}.

\bibitem[Korelc(2009)]{Kor:2009:aof}
J.~Korelc.
\newblock Automation of primal and sensitivity analysis of transient coupled
  problems.
\newblock \emph{Computational Mechanics}, 44\penalty0 (5):\penalty0 631--649,
  2009.

\bibitem[Korelc and Wriggers(2016)]{KorWri:2016:aofem}
J.~Korelc and P.~Wriggers.
\newblock \emph{Automation of Finite Element Methods}.
\newblock Springer International Publishing, 2016.

\bibitem[Kouznetsova et~al.(2002)Kouznetsova, Geers, and
  Brekelmans]{KouGeeBre:2002:msc}
V.~Kouznetsova, M.G.D. Geers, and W.A.M. Brekelmans.
\newblock Multi-scale constitutive modelling of heterogeneous materials with a
  gradient-enhanced computational homogenization scheme.
\newblock \emph{International Journal for Numerical Methods in Engineering},
  54:\penalty0 1235--1260, 2002.

\bibitem[Kouznetsova et~al.(2004)Kouznetsova, Geers, and
  Brekelmans]{KouGeeBre:2004:mss}
V.~Kouznetsova, M.G.D. Geers, and W.A.M. Brekelmans.
\newblock Multi-scale second-order computational homogenization of multi-phase
  materials: a nested finite element solution strategy.
\newblock \emph{Computer Methods in Applied Mechanics and Engineering},
  193:\penalty0 5525--5550, 2004.

\bibitem[Lahbazi et~al.(2022)Lahbazi, Goda, and Ganghoffer]{LahGodGan:2022:sis}
A.~Lahbazi, I.~Goda, and J.-F. Ganghoffer.
\newblock Size-independent strain gradient effective models based on
  homogenization methods: Applications to 3{D} composite materials, pantograph
  and thin walled lattices.
\newblock \emph{Composite Structures}, 284:\penalty0 115065, 2022.

\bibitem[Lakes(2022)]{Lak:2022:cse}
R.~S. Lakes.
\newblock Cosserat shape effects in the bending of foams.
\newblock \emph{Mechanics of Advanced Materials and Structures}, 2022.
\newblock \doi{10.1080/15376494.2022.2086328}.

\bibitem[Lee et~al.(2012)Lee, Singer, and Thomas]{LeeSinTho:2012:mnm}
J.-H. Lee, J.~P. Singer, and E.~L. Thomas.
\newblock Micro-/nanostructured mechanical metamaterials.
\newblock \emph{Advanced Materials}, 24\penalty0 (36):\penalty0 4782--4810,
  2012.

\bibitem[Leismann and Mahnken(2015)]{LeiMah:2015:coh}
T.~Leismann and R.~Mahnken.
\newblock Comparison of hyperelastic micromorphic, micropolar and microstrain
  continua.
\newblock \emph{International Journal of Non-Linear Mechanics}, 77:\penalty0
  115--127, 2015.

\bibitem[Li et~al.(2022)Li, Wang, Song, Chen, Su, Zhou, and
  Wang]{LiWanSonCheSuZhoWan:2022:osg}
A.~Li, Q.~Wang, M.~Song, J.~Chen, W.~Su, S.~Zhou, and L.~Wang.
\newblock On strain gradient theory and its application in bending of beam.
\newblock \emph{Coatings}, 12\penalty0 (9), 2022.

\bibitem[Li and Zhang(2013)]{LiZha:2013:ana}
J.~Li and X.-B. Zhang.
\newblock A numerical approach for the establishment of strain gradient
  constitutive relations in periodic heterogeneous materials.
\newblock \emph{European Journal of Mechanics - A/Solids}, 41:\penalty0 70--85,
  2013.

\bibitem[Liebold and M\"uller(2016)]{LieMue:2016:cof}
C.~Liebold and W.~H. M\"uller.
\newblock Comparison of gradient elasticity models for the bending of
  micromaterials.
\newblock \emph{Computational Materials Science}, 116:\penalty0 52--61, 2016.

\bibitem[Madeo et~al.(2015)Madeo, Neff, Ghiba, Placidi, and
  Rosi]{MadNefGhiPlaRos:2015:bgi}
A.~Madeo, P.~Neff, I.-D. Ghiba, L.~Placidi, and G.~Rosi.
\newblock Band gaps in the relaxed linear micromorphic continuum.
\newblock \emph{Zeitschrift f\"ur angewandte Mathematik und Mechanik},
  95\penalty0 (9):\penalty0 880--887, 2015.

\bibitem[Madeo et~al.(2016{\natexlab{a}})Madeo, Neff, d'Agostino, and
  Barbagallo]{MadNefBar:2016:cbg}
A.~Madeo, P.~Neff, M.~V. d'Agostino, and G.~Barbagallo.
\newblock Complete band gaps including non-local effects occur only in the
  relaxed micromorphic model.
\newblock \emph{Comptes Rendus M{\'e}canique}, 344\penalty0 (11-12):\penalty0
  784--796, 2016{\natexlab{a}}.

\bibitem[Madeo et~al.(2016{\natexlab{b}})Madeo, Neff, Ghiba, and
  Rosi]{MadNefGhiRos:2016:rat}
A.~Madeo, P.~Neff, I.-D. Ghiba, and G.~Rosi.
\newblock {Reflection and transmission of elastic waves in non-local band-gap
  metamaterials: a comprehensive study via the relaxed micromorphic model}.
\newblock \emph{Journal of the Mechanics and Physics of Solids}, 95:\penalty0
  441--479, 2016{\natexlab{b}}.

\bibitem[Madeo et~al.(2017)Madeo, Neff, Barbagallo, d'Agostino, and
  Ghiba]{MadNefBarGhi:2017:aro}
A.~Madeo, P.~Neff, G.~Barbagallo, M.~V. d'Agostino, and I.-D. Ghiba.
\newblock A review on wave propagation modeling in band-gap metamaterials via
  enriched continuum models.
\newblock In \emph{Mathematical {M}odelling in {S}olid {M}echanics}, volume~69
  of \emph{Advanced Structured Materials}, pages 89--105. Springer, Singapore,
  2017.

\bibitem[Mindlin(1964)]{Min:1964:msi}
R.~D. Mindlin.
\newblock Micro-structure in linear elasticity.
\newblock \emph{Archive for Rational Mechanics and Analysis}, 16:\penalty0
  51--78, 1964.

\bibitem[Mindlin and Eshel(1968)]{MinEsh:1968:ofsg}
R.~D. Mindlin and N.~N. Eshel.
\newblock On first strain-gradient theories in linear elasticity.
\newblock \emph{International Journal of Solids and Structures}, 4\penalty0
  (1):\penalty0 109--124, 1968.

\bibitem[Monchiet et~al.(2020)Monchiet, Auffray, and
  Yvonne]{MonAufYvo:2020:sgh}
V.~Monchiet, N.~Auffray, and J.~Yvonne.
\newblock Strain-gradient homogenization: A bridge between the asymptotic
  expansion and quadratic boundary condition methods.
\newblock \emph{Mechanics of Materials}, 143:\penalty0 103309, 2020.

\bibitem[N\'ed\'elec(1980)]{Ned:1980:mfe}
J.~C. N\'ed\'elec.
\newblock Mixed finite elements in {$\R^3$}.
\newblock \emph{Numerische Mathematik}, 35\penalty0 (3):\penalty0 315--341,
  1980.

\bibitem[N\'ed\'elec(1986)]{Ned:1986:anf}
J.~C. N\'ed\'elec.
\newblock A new family of mixed finite elements in {$\R^3$}.
\newblock \emph{Numerische Mathematik}, 50:\penalty0 57--81, 1986.

\bibitem[Neff(2006)]{Nef:2006:tcc}
P.~Neff.
\newblock The {C}osserat couple modulus for continuous solids is zero viz the
  linearized {C}auchy-stress tensor is symmetric.
\newblock \emph{Zeitschrift f\"ur Angewandte Mathematik und Mechanik},
  86\penalty0 (11):\penalty0 892--912, 2006.

\bibitem[Neff and Forest(2007)]{NefFor:2007:age}
P.~Neff and S.~Forest.
\newblock A geometrically exact micromorphic model for elastic metallic foams
  accounting for affine microstructure. modelling, existence of minimizers,
  identification of moduli and computational results.
\newblock \emph{Journal of Elasticity}, 87\penalty0 (2):\penalty0 239--276,
  2007.

\bibitem[Neff et~al.(2010)Neff, Jeong, M\"unch, and
  Ram\'ezani]{NefJeoMueRam:2010:lce}
P.~Neff, J.~Jeong, I.~M\"unch, and H.~Ram\'ezani.
\newblock Linear {C}osserat elasticity, conformal curvature and bounded
  stiffness.
\newblock In \emph{Mechanics of Generalized Continua: One Hundred Years After
  the Cosserats}, pages 55--63. Springer New York, 2010.

\bibitem[Neff et~al.(2014)Neff, Ghiba, Madeo, Placidi, and
  Rosi]{NefGhiMadPlaRos:2014:aup}
P.~Neff, I.-D. Ghiba, A.~Madeo, L.~Placidi, and G.~Rosi.
\newblock A unifying perspective: the relaxed linear micromorphic continuum.
\newblock \emph{Continuum Mechanics and Thermodynamics}, 26\penalty0
  (5):\penalty0 639--681, 2014.

\bibitem[Neff et~al.(2015)Neff, Ghiba, Lazar, and Madeo]{NefGhiLazMad:2015:trl}
P.~Neff, I.~D. Ghiba, M.~Lazar, and A.~Madeo.
\newblock The relaxed linear micromorphic continuum: well-posedness of the
  static problem and relations to the gauge theory of dislocations.
\newblock \emph{The Quarterly Journal of Mechanics and Applied Mathematics},
  68\penalty0 (1):\penalty0 53--84, 2015.

\bibitem[Neff et~al.(2020)Neff, Eidel, d`Agostino, and
  Madeo]{NefEidMad:2019:ios}
P.~Neff, B.~Eidel, M.~V. d`Agostino, and A.~Madeo.
\newblock Identification of scale-independent material parameters in the
  relaxed micromorphic model through model-adapted first order homogenization.
\newblock \emph{Journal of Elasticity}, 139:\penalty0 269--298, 2020.

\bibitem[Owczarek et~al.(2019)Owczarek, Ghiba, d'Agostino, and
  Neff]{OwcGhidAgNef:2019:nmi}
S.~Owczarek, I.-D. Ghiba, M.~V. d'Agostino, and P.~Neff.
\newblock Nonstandard micro-inertia terms in the relaxed micromorphic model:
  well-posedness for dynamics.
\newblock \emph{Mathematics and Mechanics of Solids}, 24\penalty0
  (10):\penalty0 3200--3215, 2019.

\bibitem[Placidi et~al.(2017)Placidi, Barchiesi, and
  Battista]{PlaBarBat:2017:aim}
L.~Placidi, E.~Barchiesi, and A.~Battista.
\newblock An inverse method to get further analytical solutions for a class of
  metamaterials aimed to validate numerical integrations.
\newblock In F.~dell'Isola, M.~Sofonea, and D.~Steigmann, editors,
  \emph{Mathematical Modelling in Solid Mechanics}, pages 193--210. Springer
  Singapore, Singapore, 2017.

\bibitem[Reda et~al.(2021)Reda, Alavi, Nasimsobhan, and
  Ganghoffer]{RedAlaNasGan:2021:htc}
H.~Reda, S.~E. Alavi, M.~Nasimsobhan, and J.-F. Ganghoffer.
\newblock Homogenization towards chiral {C}osserat continua and applications to
  enhanced {T}imoshenko beam theories.
\newblock \emph{Mechanics of Materials}, 155:\penalty0 103728, 2021.

\bibitem[Rizzi et~al.(2021{\natexlab{a}})Rizzi, H\"utter, Madeo, and
  Neff]{RizHueMadNef:2021:aso3}
G.~Rizzi, G.~H\"utter, A.~Madeo, and P.~Neff.
\newblock Analytical solutions of the simple shear problem for micromorphic
  models and other generalized continua.
\newblock \emph{Archive of Applied Mechanics}, 91:\penalty0 2237--2254,
  2021{\natexlab{a}}.

\bibitem[Rizzi et~al.(2021{\natexlab{b}})Rizzi, H\"utter, Madeo, and
  Neff]{RizHueMadNef:2021:aso4}
G.~Rizzi, G.~H\"utter, A.~Madeo, and P.~Neff.
\newblock Analytical solutions of the cylindrical bending problem for the
  relaxed micromorphic continuum and other generalized continua.
\newblock \emph{Continuum Mechanics and Thermodynamics}, 33:\penalty0
  1505--1539, 2021{\natexlab{b}}.

\bibitem[Rizzi et~al.(2021{\natexlab{c}})Rizzi, Khan, Ghiba, Madeo, and
  Neff]{RizKhaGhiMadNef:2021:aso2}
G.~Rizzi, H.~Khan, I.-D. Ghiba, A.~Madeo, and P.~Neff.
\newblock Analytical solution of the uniaxial extension problem for the relaxed
  micromorphic continuum and other generalized continua (including full
  derivations).
\newblock \emph{Archive of Applied Mechanics}, 2021{\natexlab{c}}.
\newblock \doi{10.1007/s00419-021-02064-3}.

\bibitem[Rizzi et~al.(2022{\natexlab{a}})Rizzi, d'Agostino, Neff, and
  Madeo]{RizdAgNefMad:2022:mam}
G.~Rizzi, M.~V. d'Agostino, P.~Neff, and A.~Madeo.
\newblock Boundary and interface conditions in the relaxed micromorphic model:
  Exploring finite-size metastructures for elastic wave control.
\newblock \emph{Mathematics and Mechanics of Solids}, 27\penalty0 (6):\penalty0
  1053--1068, 2022{\natexlab{a}}.

\bibitem[Rizzi et~al.(2022{\natexlab{b}})Rizzi, H\"utter, Khan, Ghiba, Madeo,
  and Neff]{RizHueKhaGhiMadNef:2021:aso1}
G.~Rizzi, G.~H\"utter, H.~Khan, I.-D. Ghiba, A.~Madeo, and P.~Neff.
\newblock Analytical solution of the cylindrical torsion problem for the
  relaxed micromorphic continuum and other generalized continua (including full
  derivations).
\newblock \emph{Mathematics and Mechanics of Solids}, 27\penalty0 (3):\penalty0
  507--553, 2022{\natexlab{b}}.

\bibitem[Rizzi et~al.(2022{\natexlab{c}})Rizzi, Neff, and
  Madeo]{RizNefMad:2022:msf}
G.~Rizzi, P.~Neff, and A.~Madeo.
\newblock Metamaterial for inner protection and outer tuning through a relaxed
  micromorphic approach.
\newblock \emph{Philosophical Transactions of the Royal Society A:
  Mathematical, Physical and Engineering Sciences}, 380\penalty0
  (2231):\penalty0 20210400, 2022{\natexlab{c}}.

\bibitem[Rizzi et~al.(2022{\natexlab{d}})Rizzi, Tallarico, Neff, and
  Madeo]{RizTalNefMad:2022:ttc}
G.~Rizzi, D.~Tallarico, P.~Neff, and A.~Madeo.
\newblock Towards the conception of complex engineering meta-structures:
  Relaxed-micromorphic modelling of low-frequency mechanical
  diodes/high-frequency screens.
\newblock \emph{Wave Motion}, 113:\penalty0 102920, 2022{\natexlab{d}}.

\bibitem[Rognes et~al.(2009)Rognes, Kirby, and Logg]{RogKirAnd:2009:eao}
M.~E. Rognes, R.~C. Kirby, and A.~Logg.
\newblock Efficient assembly of {$H(\mathrm{div})$ and $H(\mathrm{curl})$}
  conforming finite elements.
\newblock \emph{SIAM Journal on Scientific Computing}, 31\penalty0
  (9):\penalty0 4130--4151, 2009.

\bibitem[Roko\v{s} et~al.(2019)Roko\v{s}, Ameen, Peerlings, and
  Geers]{RokAmePeeGee:2019:mch}
O.~Roko\v{s}, M.~M. Ameen, R.~H.~J. Peerlings, and M.~G.~D. Geers.
\newblock Micromorphic computational homogenization for mechanical
  metamaterials with patterning fluctuation fields.
\newblock \emph{Journal of the Mechanics and Physics of Solids}, 123:\penalty0
  119--137, 2019.

\bibitem[Roko\v{s} et~al.(2020{\natexlab{a}})Roko\v{s}, Ameen, Peerlings, and
  Geers]{RokAmePeeGee:2020:emc}
O.~Roko\v{s}, M.~M. Ameen, R.~H.~J. Peerlings, and M.~G.~D. Geers.
\newblock Extended micromorphic computational homogenization for mechanical
  metamaterials exhibiting multiple geometric pattern transformations.
\newblock \emph{Extreme Mechanics Letters}, 37:\penalty0 100708,
  2020{\natexlab{a}}.

\bibitem[Roko\v{s} et~al.(2020{\natexlab{b}})Roko\v{s}, Zeman,
  Do\v{s}k\'{a}\v{r}, and Krysl]{RokZemDosKry:2020:ris}
O.~Roko\v{s}, J.~Zeman, M.~Do\v{s}k\'{a}\v{r}, and P.~Krysl.
\newblock Reduced integration schemes in micromorphic computational
  homogenization of elastomeric mechanical metamaterials.
\newblock \emph{Advanced Modeling and Simulation in Engineering Sciences}, 7,
  2020{\natexlab{b}}.

\bibitem[Rueger et~al.(2019)Rueger, Ha, and Lakes]{RugHaLak:2019:cel}
Z.~Rueger, C.~S. Ha, and R.~S. Lakes.
\newblock Cosserat elastic lattices.
\newblock \emph{Meccanica}, 54:\penalty0 1983--1999, 2019.

\bibitem[Sarhil et~al.(2021)Sarhil, Scheunemann, Neff, and
  Schr\"oder]{SarSchNefSch:2021:oat}
M.~Sarhil, L.~Scheunemann, P.~Neff, and J.~Schr\"oder.
\newblock On a tangential-conforming finite element formulation for the relaxed
  micromorphic model in 2{D}.
\newblock \emph{Proceedings in Applied Mathematics and Mechanics}, 21\penalty0
  (1):\penalty0 e202100187, 2021.

\bibitem[Sarhil et~al.(2023)Sarhil, Scheunemann, Schr\"oder, and
  Neff]{SarSchSchNef:2023:mts}
M.~Sarhil, L.~Scheunemann, J.~Schr\"oder, and P.~Neff.
\newblock Modeling the size-effect of metamaterial beams under bending via the
  relaxed micromorphic continuum.
\newblock \emph{Proceedings in Applied Mathematics and Mechanics}, 22\penalty0
  (1):\penalty0 e202200033, 2023.

\bibitem[Schmidt et~al.(2022)Schmidt, Kr\"uger, Keip, and
  Hesch]{SchKruKeiHes:2022:cho}
F.~Schmidt, M.~Kr\"uger, M.-A. Keip, and C.~Hesch.
\newblock Computational homogenization of higher-order continua.
\newblock \emph{International Journal for Numerical Methods in Engineering},
  123\penalty0 (11):\penalty0 2499--2529, 2022.

\bibitem[Schr\"oder et~al.(2022)Schr\"oder, Sarhil, Scheunemann, and
  Neff]{SchSarSchNef:2022:lhb}
J.~Schr\"oder, M.~Sarhil, L.~Scheunemann, and P.~Neff.
\newblock {Lagrange and $H(\curl,\cal{B})$ based Finite Element formulations
  for the relaxed micromorphic model}.
\newblock \emph{Computational Mechanics}, 70\penalty0 (6):\penalty0 1309--1333,
  2022.
\newblock \doi{10.1007/s00466-022-02198-3}.

\bibitem[Shekarchizadeh et~al.(2021)Shekarchizadeh, Abali, Barchiesi, and
  Bersani]{SheAbaBarBer:2021:iao}
N.~Shekarchizadeh, B.~E. Abali, E.~Barchiesi, and A.~M. Bersani.
\newblock Inverse analysis of metamaterials and parameter determination by
  means of an automatized optimization problem.
\newblock \emph{Zeitschrift f\"ur Angewandte Mathematik und Mechanik},
  101\penalty0 (8):\penalty0 e202000277, 2021.

\bibitem[Shekarchizadeh et~al.(2022)Shekarchizadeh, Abali, and
  Bersani]{SheAbaBer:2022:abs}
N.~Shekarchizadeh, B.~E. Abali, and A.~M. Bersani.
\newblock A benchmark strain gradient elasticity solution in two-dimensions for
  verifying computational approaches by means of the finite element method.
\newblock \emph{Mathematics and Mechanics of Solids}, 27\penalty0
  (10):\penalty0 2218--2238, 2022.

\bibitem[Shi et~al.(2022)Shi, Fantuzzi, Trovalusci, Li, and
  Wei]{ShiFanTroLiWei:2022:sfe}
F.~Shi, N.~Fantuzzi, P.~Trovalusci, Y.~Li, and Z.~Wei.
\newblock Stress field evaluation in orthotropic microstructured composites
  with holes as {C}osserat continuum.
\newblock \emph{Materials}, 15\penalty0 (18), 2022.

\bibitem[Skrzat and Eremeyev(2020)]{SkrEre:2020:ote}
A.~Skrzat and V.~A. Eremeyev.
\newblock On the effective properties of foams in the framework of the couple
  stress theory.
\newblock \emph{Continuum Mechanics and Thermodynamics}, 32:\penalty0
  1779--1801, 2020.

\bibitem[Sky et~al.(2021)Sky, Neunteufel, M\"unch, Sch\"oberl, and
  Neff]{SkyNeuMueSchNef:2021:CM}
A.~Sky, M.~Neunteufel, I.~M\"unch, J.~Sch\"oberl, and P.~Neff.
\newblock A hybrid {$H^1\times H(\mathrm{curl})$} finite element formulation
  for a relaxed micromorphic continuum model of antiplane shear.
\newblock \emph{Computational Mechanics}, 68:\penalty0 1--24, 2021.

\bibitem[Sky et~al.(2022)Sky, Neunteufel, M\"unch, Sch\"oberl, and
  Neff]{SkyNeuMueSchNef:2022:pam}
A.~Sky, M.~Neunteufel, I.~M\"unch, J.~Sch\"oberl, and P.~Neff.
\newblock Primal and mixed finite element formulations for the relaxed
  micromorphic model.
\newblock \emph{Computer Methods in Applied Mechanics and Engineering},
  399:\penalty0 115298, 2022.

\bibitem[Sky et~al.(2023)Sky, M\"unch, Rizzi, and Neff]{SkyMueRizNeff:2023:hob}
A.~Sky, I.~M\"unch, G.~Rizzi, and P.~Neff.
\newblock Higher order {B}ernstein-{B}\'ezier and {N}\'ed\'elec finite elements
  for the relaxed micromorphic model.
\newblock \emph{https://arxiv.org/abs/2301.01491}, 2023.
\newblock \doi{10.48550/ARXIV.2301.01491}.

\bibitem[Sridhar et~al.(2016)Sridhar, Kouznetsova, and
  Geers]{SirKouGee:2016:hol}
A.~Sridhar, V.~G. Kouznetsova, and M.~G.~D. Geers.
\newblock Homogenization of locally resonant acoustic metamaterials towards an
  emergent enriched continuum.
\newblock \emph{Computational Mechanics}, 57:\penalty0 423--435, 2016.

\bibitem[Sridhar et~al.(2018)Sridhar, Liu, Kouznetsova, and
  Geers]{SirLiuKouGee:2018:hec}
A.~Sridhar, L.~Liu, V.~G. Kouznetsova, and M.~G.~D. Geers.
\newblock Homogenized enriched continuum analysis of acoustic metamaterials
  with negative stiffness and double negative effects.
\newblock \emph{Journal of the Mechanics and Physics of Solids}, 119:\penalty0
  104--117, 2018.

\bibitem[Suhubl and Eringen(1964)]{SuhEri:1964:nto}
E.~S. Suhubl and A.~C. Eringen.
\newblock Nonlinear theory of micro-elastic solids-{II}.
\newblock \emph{International Journal of Engineering Science}, 2\penalty0
  (4):\penalty0 389--404, 1964.

\bibitem[Surjadi et~al.(2019)Surjadi, Gao, Du, Li, Xiong, Fang, and
  Lu]{SurGaoDuLiXioFanLu:2019:AEM}
J.~U. Surjadi, L.~Gao, H.~Du, X.~Li, X.~Xiong, N.~X. Fang, and Y.~Lu.
\newblock Mechanical metamaterials and their engineering applications.
\newblock \emph{Advanced Engineering Materials}, 21\penalty0 (3):\penalty0
  1800864, 2019.

\bibitem[Tran et~al.(2012)Tran, Monchiet, and Bonnet]{TraMonBon:2012:amb}
T.-H. Tran, V.~Monchiet, and G.~Bonnet.
\newblock A micromechanics-based approach for the derivation of constitutive
  elastic coefficients of strain-gradient media.
\newblock \emph{International Journal of Solids and Structures}, 49\penalty0
  (5):\penalty0 783--792, 2012.

\bibitem[Trinh et~al.(2012)Trinh, J\"anicke, Auffray, Diebels, and
  Forest]{TriJaeAufDieFor:2012:eog}
D.~K. Trinh, R.~J\"anicke, N.~Auffray, S.~Diebels, and S.~Forest.
\newblock Evaluation of generalized continuum substitution models for
  heterogeneous materials.
\newblock \emph{International Journal for Multiscale Computational
  Engineering}, 10(6):\penalty0 527--549, 2012.

\bibitem[Trovalusci and Pau(2014)]{TroPau:2014:dom}
P.~Trovalusci and A.~Pau.
\newblock Derivation of microstructured continua from lattice systems via
  principle of virtual works: the case of masonry-like materials as micropolar,
  second gradient and classical continua.
\newblock \emph{Acta Mechanica}, 225:\penalty0 157--177, 2014.
\newblock \doi{10.1007/s00707-013-0936-9}.

\bibitem[Waseem et~al.(2020)Waseem, Heuz\'e, Stainier, Geers, and
  Kouznetsova]{WasHeuStaGeeKou:2020:ecf}
A.~Waseem, T.~Heuz\'e, L.~Stainier, M.G.D. Geers, and V.G. Kouznetsova.
\newblock Enriched continuum for multi-scale transient diffusion coupled to
  mechanics.
\newblock \emph{Advanced Modeling and Simulation in Engineering Sciences},
  7:\penalty0 1--32, 2020.

\bibitem[Weeger(2021)]{Wee:2021:nho}
O.~Weeger.
\newblock Numerical homogenization of second gradient, linear elastic
  constitutive models for cubic {3D} beam-lattice metamaterials.
\newblock \emph{International Journal of Solids and Structures}, 224:\penalty0
  111037, 2021.

\bibitem[Yang and M\"uller(2021)]{YanMue:2021:seo}
H.~Yang and W.~H. M\"uller.
\newblock Size effects of mechanical metamaterials: a computational study based
  on a second-order asymptotic homogenization method.
\newblock \emph{Archive of Applied Mechanics}, 91:\penalty0 1037--1053, 2021.

\bibitem[Yang et~al.(2020)Yang, Abali, Timofeev, and
  M\"uller]{YanAbaTimMue:2020:dom}
H.~Yang, B.~E. Abali, D.~Timofeev, and W.~H. M\"uller.
\newblock Determination of metamaterial parameters by means of a homogenization
  approach based on asymptotic analysis.
\newblock \emph{Continuum Mechanics and Thermodynamics}, 32:\penalty0
  1251--1270, 2020.

\bibitem[Yang et~al.(2021)Yang, Timofeev, Abali, Li, and
  M\"uller]{YanTimAbaLiMue:2021:vos}
H.~Yang, D.~Timofeev, B.~E. Abali, B.~Li, and W.~H. M\"uller.
\newblock Verification of strain gradient elasticity computation by analytical
  solutions.
\newblock \emph{Zeitschrift f\"ur Angewandte Mathematik und Mechanik},
  101\penalty0 (12):\penalty0 e202100023, 2021.

\bibitem[Yang et~al.(2022)Yang, Abali, M\"uller, Barboura, and
  Li]{YanAbaMueetal:2022:voa}
H.~Yang, B.~E. Abali, W.~H. M\"uller, S.~Barboura, and J.~Li.
\newblock Verification of asymptotic homogenization method developed for
  periodic architected materials in strain gradient continuum.
\newblock \emph{International Journal of Solids and Structures}, 238:\penalty0
  111386, 2022.

\bibitem[Yin et~al.(2021)Yin, Xiao, Deng, Zhang, Liu, and
  Gu]{YinXiaDenetal:2021:iao}
S.~Yin, Z.~Xiao, Y.~Deng, G.~Zhang, J.~Liu, and S.~Gu.
\newblock Isogeometric analysis of size-dependent {Bernoulli-Euler} beam based
  on a reformulated strain gradient elasticity theory.
\newblock \emph{Computers and Structures}, 253:\penalty0 106577, 2021.

\bibitem[Yu et~al.(2018)Yu, Zhou, Liang, Jiang, and Wu]{YaZhoLiaJiaWu:2018:mma}
X.~Yu, J.~Zhou, H.~Liang, Z.~Jiang, and L.~Wu.
\newblock Mechanical metamaterials associated with stiffness, rigidity and
  compressibility: A brief review.
\newblock \emph{Progress in Materials Science}, 94:\penalty0 114--173, 2018.

\bibitem[Yvonnet et~al.(2020)Yvonnet, Auffray, and
  Monchiet]{YvoAufMon:2020:cso}
J.~Yvonnet, N.~Auffray, and V.~Monchiet.
\newblock Computational second-order homogenization of materials with effective
  anisotropic strain-gradient behavior.
\newblock \emph{International Journal of Solids and Structures},
  191-192:\penalty0 434--448, 2020.

\bibitem[Zadpoor(2016)]{Zad:2016:mmm}
A.~A. Zadpoor.
\newblock Mechanical meta-materials.
\newblock \emph{Materials Horizons}, 3:\penalty0 371--381, 2016.

\bibitem[Zhi et~al.(2022)Zhi, Poh, Tay, and Tan]{ZhiPohTayTan:2022:dfm}
J.~Zhi, L.~H. Poh, T.-E. Tay, and V.~B.~C. Tan.
\newblock Direct {FE2} modeling of heterogeneous materials with a micromorphic
  computational homogenization framework.
\newblock \emph{Computer Methods in Applied Mechanics and Engineering},
  393:\penalty0 114837, 2022.

\bibitem[Zohdi and Wriggers(2005)]{ZohWri:2005:ait}
T.~I. Zohdi and P.~Wriggers.
\newblock \emph{An Introduction to Computational Micromechanics}.
\newblock Springer Berlin, Heidelberg, 2005.

\end{thebibliography}
}
%=== appendix

\begin{appendix}
\end{appendix}

\end{document}